\documentclass[MNSC]{informs}
\OneAndAHalfSpacedXII 
\usepackage{endnotes}
\usepackage{subcaption}
\usepackage{hyperref}
\usepackage{lipsum}
\let\footnote=\endnote

\newcommand{\R}{\mathbb{R}} 
\newcommand{\Q}{\mathbb{Q}} 
\newcommand{\Z}{\mathbb{Z}} 
\newcommand{\N}{\mathbb{N}} 
\newcommand{\E}{\mathbb{E}} 

\newcommand{\numsym}[2]{\stackrel{\scriptscriptstyle(\mkern-1.5mu#1\mkern-1.5mu)}{#2}}
\newenvironment{cproof}[2] {\noindent \textit{Proof of {#1} {#2}.}\\\noindent}{\strut\hfill$\square$\\} 

\usepackage{dsfont}
\usepackage{comment}
\usepackage{hyperref}
\usepackage{xcolor}
\usepackage{booktabs}
\usepackage{float}
\usepackage{algorithm}
\usepackage{algpseudocode}
\usepackage{algorithmicx}
\usepackage{multirow}
\usepackage{pifont}
\algdef{SE}[DOWHILE]{Do}{doWhile}{\algorithmicdo}[1]{\algorithmicwhile\ #1}%
\hypersetup{
    colorlinks=true,
    linkcolor=blue,
    filecolor=magenta,      
    urlcolor=cyan,
    pdftitle={Overleaf Example},
    pdfpagemode=FullScreen,
    citecolor=blue
    }
\usepackage{dirtytalk}
\usepackage{svg}
\usepackage{tikz}
\usetikzlibrary{decorations.pathreplacing}
\usepackage{tabularx}
\usepackage{makecell}
\newcolumntype{Y}{>{\centering\arraybackslash}X}

\usepackage{natbib}
\bibpunct[, ]{(}{)}{,}{a}{}{,}%
\def\bibfont{\small}%

\TheoremsNumberedThrough   
\ECRepeatTheorems

\EquationsNumberedThrough
\usepackage{bm}
\usepackage{stmaryrd}

\newcommand{\ear}{\lambda}
\newcommand{\arate}{\mu}
\newcommand{\rv}{X}

\RUNAUTHOR{Bouscary, Danaf, Amin}

\RUNTITLE{Integrated Bundling and Pricing of Unique Items}

\TITLE{Integrated Bundling and Pricing of Unique Items}
\ARTICLEAUTHORS{%
\AUTHOR{Maxime Bouscary}
\AFF{Operations Research Center, Massachusetts Institute of Technology, Cambridge, MA 02139, \EMAIL{mbscry@mit.edu}}
\AUTHOR{Mazen Danaf}
\AFF{Uber Freight, 433 W. Van Buren St. Chicago, IL 60607,  \EMAIL{mdanaf@uberfreight.com}}
\AUTHOR{Saurabh Amin}
\AFF{Laboratory for Information and Decision Systems, Massachusetts Institute of Technology, Cambridge, MA 02139, \EMAIL{amins@mit.edu}}
} 

\ABSTRACT{
Retailers have significant potential to improve recommendations through strategic bundling and pricing. By taking into account different types of customers and their purchasing decisions, retailers can better accommodate customer preferences and increase revenues while reducing unsold items. We consider a retailer seeking to maximize its expected revenue by selling unique and non-replenishable items over a finite horizon. The retailer may offer each item individually or as part of a bundle.

Our approach provides tractable bounds on expected revenue that are tailored to unique items and suitable for a rich class of choice models. We leverage these bounds to propose a bundling algorithm that efficiently selects bundles in a column-generation fashion. Under the multinomial logit model, our bounds are asymptotically optimal as the expected number of arrivals $\ear$ grows, yielding a performance bound in $O(1/\ear)$. In contrast, we show that both the static and fluid approximations are not asymptotically optimal. Moreover, we propose a greedy algorithm that allows for tractable dynamic bundling.

We show through numerical experiments that third-party logistics providers (3PL) in particular can benefit from improved load recommendations and pricing decisions. By not only minimizing shipping costs, but also reducing empty miles and suggesting options that better match carrier preferences, our methodology benefits digital brokerage platforms while contributing to the transportation industry's sustainability efforts. In a Texas Triangle simulation powered by Uber Freight data, dynamic bundling and pricing reduces costs by 6\%, empty miles by 25\%, and more than halves the number of unmatched loads over dynamic pricing alone.
}


\begin{document}

\maketitle
\section{Introduction}

Online retailing has revolutionized the way companies operate their sales. Among the key features of online retailing is the ability to easily implement dynamic pricing, which is particularly valuable when inventory is limited. Alternatively, offering discounted bundles of items has traditionally been a strategy to increase sales and clear inventories. Employed in similar contexts, these two mechanisms influence each other when used simultaneously. Indeed, bundling certain items affects the optimal price of all other items, and the gains from bundling can be negated by mispricing.
We propose an innovative approach that integrates bundling with dynamic pricing to unify both mechanisms and alleviate adverse interactions.

We consider a retailer owning a set of \textit{unique} items, that is, initially available in quantity one and cannot be replenished. The retailer can decide to offer each item as an individual item or as part of a bundle of two or more items, and we refer to each such offer as an \textit{option}. Importantly, we assume no relationship between the quality of a bundle and the quality of the items it contains. As such, a bundle formed by two high-quality items can be of low quality, and low-quality items can form a high-quality bundle. This allows to account for the complementarity and substitutability of items. Item uniqueness allows for fine-grained pricing and bundling. For example, considering each seat on an airplane as a unique item to be sold allows for price differentiation based on various features (window access, legroom, distance to the engine) and bundling based on the relative position of the seats. Options are static in the sense that they are determined at the first time period and remain unchanged thereafter. The approach we develop is extensible to dynamic bundling where options are updated throughout the selling season by frequently re-optimizing.

Customers arrive sequentially and we refer to the expected number of arrivals throughout the selling season as the \textit{demand}. Arriving customers are segmented into different customer \textit{types}. The retailer can observe the types, and customers of the same type perceive the quality of options identically. Different types may correspond to distinct geographical regions or socioeconomic characteristics such as age, gender, or income level. 
For each arrival, the retailer observes the type of customer and adjusts the price of the options offered. Based on the quality of the options and the prices displayed, customers can decide to accept an option or leave the marketplace empty-handed. At the end of the selling season, items are salvaged at an item-specific value.

Our objective is to maximize the retailer's expected total revenue by forming a set of options that partition the set of items and by pricing these options dynamically. Specifically, we focus on the question: \textit{Can bundling improve expected revenue alongside dynamic pricing, and if so, how to determine which bundles to offer to maximize expected revenue?} Although our performance guarantees focus on the multinomial logit (MNL) model, we provide a general bundling framework that is amenable to a large class of choice models, including but not limited to multinomial logit and nested logit models. 

\begin{figure}[H]
\centering
\begin{minipage}{0.35\textwidth}
\small
\vfill
\raggedright
Quality of options:
\vspace{2mm}
\begin{center}
\begin{tikzpicture}[scale=0.8]
\tikzset{every node/.style={font=\small}}
\pgfmathsetmacro{\v}{2}
\pgfmathsetmacro{\r}{0.55}
\pgfmathsetmacro{\h}{0.4}

\node (A) at (0,0) {Item $a$};
\node (B) at (\v,0) {Item $b$};
\node (C) at (2*\v,0) {Item $c$};

\draw[thick] (-\r,\h) -- (-\r,\h+0.1) -- (\r,\h+0.1) -- (\r,\h);
\draw[thick] (\v-\r,\h) -- (\v-\r,\h+0.1) -- (\v+\r,\h+0.1) -- (\v+\r,\h);
\draw[thick] (2*\v-\r,\h) -- (2*\v-\r,\h+0.1) -- (2*\v+\r,\h+0.1) -- (2*\v+\r,\h);
\node at (0,2*\h) {1};
\node at (\v,2*\h) {1.5};
\node at (2*\v,2*\h) {2.5};

\draw[thick] (-\r,-\h) -- (-\r,-\h-0.1) -- (\v+\r,-\h-0.1) -- (\v+\r,-\h);
\draw[thick] (\v-\r,-2*\h) -- (\v-\r,-2*\h-0.1) -- (2*\v+\r,-2*\h-0.1) -- (2*\v+\r,-2*\h);
\node at (\v/2,-2*\h) {3};
\node at (3*\v/2,-3*\h) {4};
\end{tikzpicture}
\end{center}
\vfill
Considered sets of options:
\vspace{-3mm}
\begingroup
\setlength{\jot}{5pt}
\begin{align*}
    S_0 &= \{a,b,c\}\\ \noalign{\vskip-2mm}
    S_1 &= \{(a,b),c\}\\ \noalign{\vskip-2mm}
    S_2 &= \{a,(b,c)\}
\end{align*}
\endgroup
\vfill
\end{minipage}
\begin{minipage}{0.6\textwidth}
    \centering    \includegraphics[width=\textwidth]{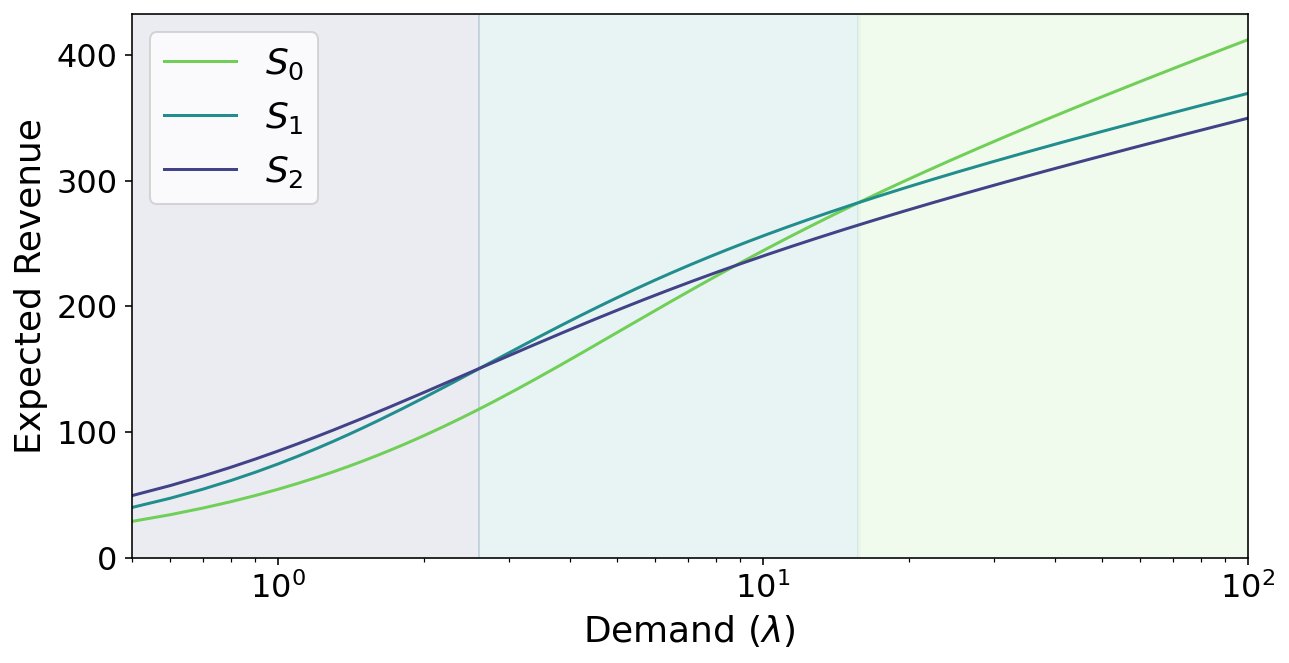}
\end{minipage}
\caption{Expected revenue under optimal pricing, function of the demand and set of options.}
\label{ex:introduction}
\end{figure}
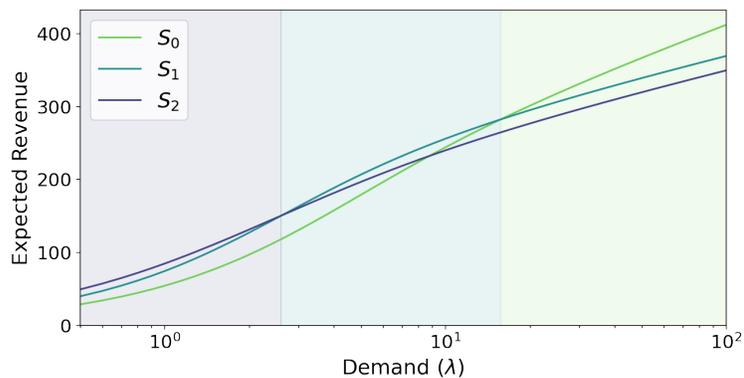

Let us consider a generic example of the expected revenue under optimal pricing policy, summarized in Figure \ref{ex:introduction}. In this example with a single type of customer, the retailer has $L=3$ items: $a$, $b$, and $c$. Note that $a$ and $b$ are complementary since the quality of bundle $(a,b)$ is greater than the sum of the qualities of $a$ and $b$, whereas $b$ and $c$ are neither complementary nor substitutable. We show the expected revenue of 3 key sets of options as the demand varies. The set of options $S_0$ offers each item individually, $S_1$ offers bundle $(a,b)$ and item $c$, and $S_2$ offers bundle $(b,c)$ and item $a$. For each considered set of options, the optimal expected revenue and pricing policy (given in the Appendix [\ref{apx:intro_example}]) are obtained by solving the dynamic program \eqref{eq:bellman_equation}. When the demand, $\ear$, is low, the retailer maximizes expected revenue by offering set $S_2$ as it contains the option with the highest quality (bundle $(b,c)$), the other options having little effect on expected revenue. As the demand increases, the retailer can expect to sell all items. In this case, it is best for the retailer to offer the set $S_1$ exploiting the complementarity of $a$ and $b$ to maximize the sum of the quality of the options. Finally, as the demand grows large, the retailer should opt for the set of options with the largest number of options, $S_0$, as it can leverage dynamic pricing and the variability in customers' utility to make the most of each option. This example illustrates how the optimal set of options varies with the demand under optimal prices. It also points out that it may be optimal not to bundle, and suggests that the retailer may particularly benefit from bundling when demand is low or when the complementarity of items allows for an increase in the total quality of options. Our framework, which is not restricted to a single customer type, tractably selects a set of options from all feasible sets in order to maximize expected revenue under optimal pricing.\\

Dynamic pricing lies at the core of our problem. As explained in 
\cite{Inventory_Competition_Dynamic_Consumer_Choice}, \cite{Stocking_Assortments_Consumer_Substitution}, and \cite{Static_Approximation_Dynamic_Demand_Substitution} considering a similar multi-item and multi-period pricing problem with dynamic demand substitution, the availability of an item is jointly determined by the inventory levels of all items, resulting in a dynamic program subject to the curse of dimensionality. Solving such a pricing problem to optimality requires enumerating a number of states that is exponential in the number of items. For this reason, although the single-period pricing problem admits a near-closed form solution given in \cite{Optimizing_Multinomial_Logit_Profit_Functions}, solving the multi-period problem to optimality remains intractable as the number of items increases. Due to the uniqueness of items and the absence of replenishment, our pricing problem is all the more affected by dynamic substitution as the sale of an item has a significant effect on future sales, while the sale of an item in a context of high inventories does not substantially affect the marginal value of its replicas and other items. The uncertain number of arrivals and the dynamic interplay between the retailer's offers and customers' decisions are other key aspects of the complexity of our problem. Specifically, since customers' decisions depend on the retailer's offer, the use of techniques relying on hindsight optimization, such as the framework developed in \cite{online_allocation_pricing} or \cite{Approximations_Stochastic_DP_Information_Relaxation_Duality}, is prohibited.

\subsection{Main Contributions}

Our approach to solving the joint bundling and pricing of unique items involves analyzing the optimal solution of the pricing problem (expected revenues, prices, choice probabilities) as demand varies. We develop tractable bounds on the optimal revenue for a given set of options by reducing the state space of dynamic pricing to a single state for each time period. We then develop bundling algorithms leveraging these bounds and evaluate the performance of our bounds and algorithms on real-world data. Our main contributions consist of the following:
\begin{enumerate}
    \item \textbf{Analysis of the pricing problem under the MNL model.} We introduce the unique item pricing problem and provide a closed-form solution under the MNL model as the demand, $\ear$, tends to infinity in Theorem \ref{theorem:asymptotic_costs}, allowing us to develop insights into the optimal revenue improvement achievable via bundling. In particular, depending on the demand, we find that bundling substitutable items can increase expected revenue while bundling complementary items can decrease it. We then show that optimal bundling jointly maximizes expected revenue and the attractiveness of the resulting set of options in Theorem \ref{theorem:revenue_utility}.
    \item \textbf{Two classes of approximations tailored to the unique items pricing problem and suitable for a rich class of choice functions.} We first develop bounds on the optimal revenue considering best and worst-case scenarios in terms of option availability. Unlike the fluid approximation, our upper bound shows asymptotic optimality. We then develop an approximation framework based on a probabilistic approach to item availability estimation. Under the class of choice functions we consider, which requires choice functions to admit a non-increasing and convex continuous extension with respect to item availability, this approach yields a lower bound on the optimal revenue. This bound significantly outperforms the static lower bound in numerical experiments, and we show in Theorem \ref{theorem:performance_bound} that it is asymptotically optimal under the MNL model as its relative difference with our upper bound is $O\left(1/\ear\right)$.
    \item \textbf{Tractable bundling algorithms.} Using the established bounds, we develop a bundling algorithm compatible with a large class of choice models that efficiently explores and evaluates candidate sets of options to maximize our lower bound on expected revenue. We then provide a greedy heuristic specific to the MNL model enabling dynamic bundling of hundreds of items.
    \item \textbf{Numerical experiments on a case study with real data from the freight industry.} Using historical data from Uber Freight in the Texas Triangle (including freight shipments between San Antonio, Austin, Dallas, and Houston), we benchmark our bounds against the static and fluid approximations and analyze the revenue improvement induced by our bundling algorithms. In a 2-year simulation, we find that dynamic bundling and pricing reduces costs by 6\%, empty miles by 25\%, and more than halves the number of unmatched loads over dynamic pricing alone. The findings suggest that offering bundles results in a win-win situation for digital brokerage platforms, provided that the bundles and prices offered are judiciously selected by the platform.\\
\end{enumerate}



\subsection{Related Work}

The two papers most closely related to our work are \cite{Static_Approximation_Dynamic_Demand_Substitution}, considering dynamic substitution in a multi-item multi-period pricing problem, and \cite{Data_Driven_Personalized_Bundle_Pricing}, dynamically offering a bundle of items to arriving customers. While these two papers emphasize how the performance of their approach improves as the problem becomes less inventory constrained, our unique item pricing problem lies on the other side of the spectrum. We consider three areas of the literature that are most closely related to our work: dynamic pricing, assortment optimization, and bundling.

\textbf{Dynamic pricing}. The dynamic pricing literature, which originally focused primarily on a single item, is well summarized by \cite{overview_pricing_models_revenue_management}. However, there is a large body of literature on subjects covering a wide variety of consumer utility models and multiple items.
\cite{Dynamic_Pricing_Strategies_Multiproduct} considers a multiproduct, single-resource setting and proposes a fluid approximation that reduces the problem to a single-product problem. In general, fluid approximations are used for tractability \cite{Clearance_Pricing_Zara} or to derive an upper bound on the expected profit \cite{Posted_Price_Auction_Mechanisms_Freight_Transportation_Marketplaces}. The quality of this bound improves as inventory and demand scale proportionally, but need not be tight if only demand scales for a fixed inventory level. In fact, we show that this bound is not asymptotically optimal when only demand grows [\ref{prop:fluid_approx_asymptotical}]. In the case of network-based problems with relocating resources, \cite{Pricing_Optimization_Shared_Vehicle_Systems} and \cite{Dynamic_pricing_relocating_resources_large_networks} develop tractable solutions that maximize the steady-state rate of reward accumulation by ignoring demand thinning due to resource unavailability, or using Lagrangian relaxation on the positive inventory constraints. Those approximations are natural when resource levels are high, but become inappropriate when inventories are sparse. Finally, the substitution effect, which is ubiquitous when inventories are sparse, is considered in a number of papers: \cite{DP_Inventory_Control_Substitute_Products} provides an exact dynamic pricing policy along with an analysis of optimal pricing patterns, while \cite{Static_Approximation_Dynamic_Demand_Substitution} develops a static approximation framework based on a fluid network model and a service-inventory mapping. The focus of these works is either to derive insights into optimal solutions or to develop approximations that prove to perform well when inventories are large. In contrast, our work considers non-replenishable inventories whose initial level are set to one, in which case common approximations such as the static and fluid approximations perform poorly.

\textbf{Assortment Optimization}. \cite{Kok2009} provides an overview of the early work on single-period assortment planning under inventory or budget constraints, including extensions to other consumer choice models and various dynamic substitution effects. Several contributions address variations of the assortment optimization problem. \cite{Dynamic_Assortment_Optimization_MNL_Capacity_Constraint} and \cite{Assortment_Optimization_Consider_Then_Choose} provide efficient practical solutions but consider unlimited inventory, whereas constrained inventory is a key aspect of our problem. Other works, such as \cite{Revenue_Management_General_Discrete_Choice_Model_Consumer_Behavior} and
\cite{Approximation_Algorithms_Dynamic_Assortment_Optimization_Models} consider cardinality-constrained assortment under limited inventories but assume customer preferences independent of time or price. As a result, they cannot assess the impact of dynamic pricing on the evolution of the optimal assortment over time. \cite{Assortment_Short_Lived_Products} and \cite{Nonparametric_Joint_Assortment_Price} study a multi-period problem with exogenous prices resulting in optimal assortments potentially significantly different from assortments considering integrated pricing. Finally, \cite{Optimal_Assortments_Endogenous_Pricing} provides some structural properties of the single-period optimal assortment with endogenous prices, but does not specify a practical method for constructing such an assortment. In contrast, we provide structural properties of the optimal assortment both for the single-period problem and for the multi-period problem in the asymptotic regime as demand grows, as well as algorithms to build an assortment.

\textbf{Bundling}. A related line of work studies the use of bundles of items with potential discounts or premiums. In a single-period setup, \cite{bakos1999bundling} shows the value of large bundles of items with zero marginal cost and claims that pure bundling can be unprofitable when marginal costs are significant. We find similar results under dynamic pricing and express the value of bundling as a function of the demand, as discussed in Section \ref{section:DP}. Other papers address multi-period bundling: \cite{Data_Driven_Personalized_Bundle_Pricing} focuses on bundle recommendation and pricing, but is limited to recommending a single bundle in addition to other items offered individually. In contrast, our work extends to the selection and pricing of sets of options potentially containing a large number of bundles. \cite{On_Pricing_Composition_Bundles} derives an efficient method for identifying the profit-maximizing bundle in the single-period case with linear bundle attractiveness, whereas we consider a finite inventory from which we want to extract the maximum revenue over a finite horizon. Finally, \cite{Optimal_Bundling_Pricing_Independently_Valued_Products} and \cite{armstrong2013more} study the mechanisms of pure and mixed bundling in a two-item market, but do not provide a practical way to bundle a large number of items. In the freight industry, and more specifically in third-party logistics, our approach complements studies that already suggest that bundling can increase truck utilization and reduce freight emissions such as \cite{Combinatorial_Bid_Generation}, \cite{Decision_TL_Combinatorial_Auctions}, \cite{InVehicleConsolidationLessThenTruckload}, and \cite{heilmann2020information}. By offering low empty mile options in line with carrier preferences, we show that freight platforms can drive the decarbonization of the transportation industry while improving expected revenues and the attractiveness of their offerings.\hfill\\

\section{Problem Formulation} \label{sec:formulation}\hfill

Consider a retailer seeking to maximize its expected revenue by offering a set of unique and
non-replenishable items over a finite horizon. The retailer owns $L$ items, indexed by $\mathcal{L} = \{1,...,L\}$. All items are initially available in quantity one, and item $j$ has a salvage value $\xi_j$. Options can either consist of a single item or a sequence, also called \textit{bundle}, of two or more items. An option $i$ is described by two quantities: a salvage value $\xi_i$, and a quality vector $q_i$ where $q_i^\omega$ is the perceived quality of the $i$-th option for a customer of type $\omega$. The salvage value of a bundle is the sum of the salvage values of the items included in that bundle. Note however that we do not make any assumption on the relation between the quality of a bundle and the quality of the items it contains, allowing us to model complementary and substitutable items. Moreover, we allow the order of items in a bundle to influence the quality of that bundle. This aspect of our formulation is essential when items in a bundle are displayed or consumed in a particular order.

To ensure that the formed options meet the retailer's commercial requirements, we require sets of options to include at most $K_s$ options of size at most $K_b$. We define $\mathcal{O}$ as the set of all options of size at most $K_b$. We consider pure bundling which requires that sets of options contain non-overlapping options i.e. each item belongs to exactly one option. This ensures that the revenue does not increase as a result of customers considering additional options, and allows us to isolate the effect of bundling on revenue. Indeed, offering bundles in addition to all individual items would always yield higher revenue provided the options are priced optimally. We let $\mathcal{P}(\mathcal{L})$ be the collection of all sets formed by options included in $\mathcal{O}$ and partitioning items $\mathcal{L}$.


Let $\ear$ denote the demand i.e. the expected number of arrivals throughout the selling season. A season is discretized by time periods indexed by an integer $t \in \{T(\ear),T(\ear)-1,...,0\}$ with smaller values of $t$ representing later points in time. We assume that the time discretization is sufficiently fine-grained such that we can ignore the probability that more than one arrival occurs during a period. The probability of an arrival in each period is denoted by $\arate > 0$. Time $t=0$ corresponds to the end of the season and $T(\ear) = \left\lfloor \ear / \arate \right\rfloor$ is the total number of periods during the season. The retailer selects a set of options $S \in \mathcal{P}(\mathcal{L})$ before the season begins and cannot change it afterward. This set constitutes the retailer's offer, and the options it contains are displayed to the customers until sold or season ends.

Arriving customers can be one of the $k$ different customer types $\Omega = \{\omega_1,...,\omega_k\}$. We assume that the type of arriving customers are independent and identically distributed random variables $\rv_{T(\ear)},...,\rv_1$ characterized by a probability mass function $p_{\rv}(\omega) = \mathbb{P}(\rv = \omega)$ for all $\omega \in \Omega$. Hence, at each period, a customer of type $\omega$ arrives with probability $\arate \cdot p_{\rv}(\omega)$, and no customer arrives with probability $1-\arate$. For each arrival, the retailer observes the type of customer and sets a price on each available option.

The analysis of the optimal solution structure (bundles and prices) in Section \ref{section:DP} and the asymptotic performance guarantees of our bounds as $\ear$ increases in Section \ref{section:bounds} are specific to MNL models whose price sensitivity is constant and identical across customer types.
 Under the MNL model, each customer chooses the option that maximizes her utility. The utility of option $i$ at price $p_i^\omega$ for customer of type $\omega \in \Omega$ is $u_i^\omega = q_i^\omega + \beta_p p_i^\omega+ \epsilon_{i}$ where $\beta_p<0$ is the price sensitivity and $\epsilon_{i}$ is an additive random component. Under the usual assumption that the random component, which represents errors in the modeler’s ability to represent all of the elements which influence the utility of an option to an individual, is independently and identically Gumbel distributed across options, the probability $\rho_i^\omega$ that a customer of type $\omega$ selects option $i$ among the set $S$ becomes:
\begin{align}
    \rho_i^\omega(p, S) = \frac{e^{q_i^\omega + \beta_p p_i^\omega}}{1+\sum_{j \in S} e^{q_j^\omega + \beta_p p_j^\omega}} ,\label{eq:rho_definition}
\end{align}
and the probability that no option is selected, denoted $\rho_0^\omega$, is given by:
\begin{align*}
    \rho_0^\omega(p, S) = \frac{1}{1+\sum_{j \in S} e^{q_j^\omega + \beta_p p_j^\omega}}.
\end{align*}

In the following, we introduce the class of choice functions, which includes notably MNL models, latent class models, and nested logit models. The bounds of Section \ref{section:bounds} and the bundling framework of Section \ref{section:Bunlding} are suitable for choice models that fall into this class.

\vspace{0.25cm}
\begin{definition}[Choice Function] \label{def:choice_funciton}
A choice function $f$ is such that for all $S \in \mathcal{P}(\mathcal{L})$, $i\in S$ and $\omega \in \Omega$, $f_i^\omega(p, S)$ denote the probability that a customer of type $\omega$ accepts option $i$ when the set of options $S$ is displayed at price vector $p$, with $f^\omega(\cdot, S):{\R}^{|S|} \to [0,1]^{|S|}$ differentiable. A choice function $f$ must satisfy, for all $S \in \mathcal{P}(\mathcal{L})$, $\omega \in \Omega$:
\begin{enumerate}
    \item[\textit{p.1}] $\sum_{i \in S} f_i^\omega(p, S) \leq 1$, $\forall p \in \R^{|S|}$,
    \item[\textit{p.2}] $f_i^\omega(p,S)=0$, $\forall i \notin S$, $\forall p \in \R^{|S|}$,
    \item[\textit{p.3}] $\lim_{p_i \to \infty} p_i \cdot f_i^\omega(p,S)=0$, $\forall i \in S$,
    \item[\textit{p.4}] There exists $\tilde{f}^\omega(p, a): {\R}^{|S|} \times [0,1]^{|S|} \to [0,1]^{|S|}$ such that $\tilde{f}_i^\omega$ is non-increasing and convex in $a$, and $\tilde{f}_i^\omega(p,a(S')) = f_i^\omega(p, S')$ for all $S' \subseteq S$ such $i \in S'$, where $a_j(S) = \mathds{1}_{\{j \in S\}}$ for all $j \in \mathcal{L}$,
        \item[\textit{p.5}] The mapping $p \mapsto p + \frac{f_i^\omega(p,\{i\})}{\frac{\partial}{\partial p_i}f_i^\omega(p,\{i\})}$ is well-defined and non-decreasing in $p_i$ for all $i \in S$.
\end{enumerate}
\end{definition}\vspace{0.25cm}\hfill

The first three properties guarantee that the sum of options in $S$ does not exceed 1, only options in $S$ can be accepted by customers, and the revenue of any option tends to 0 as its price tends to infinity. The last two are technical properties ensuring that $f$ admits a non-increasing and convex continuous extension with respect to $S$, and that the revenue-maximizing price of an option is increasing with its marginal value.
We have that $\rho$ is a valid choice function. Indeed, properties \textit{p.1}, \textit{p.2}, \textit{p.3} and \textit{p.5} can easily be verified, and property \textit{p.4} is proved in Lemma \ref{lemma:convexity_tilde_rho}. Throughout this paper, we denote by $f$ any choice function and by $\rho$ the choice function corresponding to the MNL model.

If option $i$ is selected by a customer at price $p_i^\omega$, the revenue of the retailer increases by $p_i^\omega$, all items included in option $i$ are sold, and option $i$ is no longer available in future time periods. Items cannot be replenishment during the season. As a result, sold items remain out of stock until the end of the season. When the season ends, unsold items are salvaged for a revenue equal to the sum of their salvage values.

Our objective is to maximize the retailer's total expected revenue throughout the season by forming a set of options and pricing them dynamically.
Given a set of options $S \in \mathcal{P}(\mathcal{L})$, the optimal expected revenue can be obtained by solving a dynamic program. Let the value function, denoted $V^*_t(S)$, be defined as the optimal expected revenue obtainable from time $t$ to the terminal time $0$ given that the set of available options at time $t$ is $S$. The Bellman equation is then:
\begin{align}
t\in\{0,1,\hdots,T(\ear)\}, \quad V^*_t(S) = \left\{\begin{tabular}{ll}
    $\displaystyle \sum_{i \in S} \xi_i$ & if $t=0$ \\
    $\displaystyle V^*_{t-1}(S) + \arate \cdot \E_\rv\left[ R\left(f^\rv, S, \Delta V^*_t(S)\right)\right]$ & otherwise \tag{$\mathcal{P}_\text{UIP}$} \label{eq:bellman_equation}
\end{tabular}\right.
\end{align}
where $\Delta V^*_t(S)$ is the vector of marginal values in state $S$, satisfying $\Delta_{i} V^*_t(S) \triangleq V^*_{t-1}(S) - V^*_{t-1}(S\backslash\{i\})$, and
\begin{equation} \label{eq:inner_problem}
    R(g, S, \Delta) = \max_{p \in \R^{|S|}} \sum_{i \in S} g_i(p, S)(p_i - \Delta_i) \tag{$\mathcal{P}_\mathcal{I}$}
\end{equation}
is the the single-stage single-type pricing problem.
We refer to problem \eqref{eq:bellman_equation} as the \textit{Unique Items Pricing} (UIP) problem. Compared to most stochastic and multi-item dynamic pricing problems in the literature, such as the airline seat management problem described in \cite{DP_Airline_Seat_Management_Multiple_Legs}, the initial state of the UIP problem has small inventory levels, namely one unit per item. The UIP problem is hard in itself as it suffers from the curse of dimensionality. Indeed for any $S \in \mathcal{P}(\mathcal{L})$ we have $|S| \geq L/K_b$ so the number of states of the dynamic program at each time period is at least $2^{L/K_b}$, which grows exponentially in $L$. Furthermore, as indicated in Proposition \ref{prop:fluid_approx_asymptotical} and Table \ref{tab:bounds_perf}, the fluid approximation of the UIP problem can be weak due to small inventory levels. Finally, the set of options maximizing the total expected revenue over the season is obtained by solving:
\begin{equation}
   \max_{S \in \mathcal{P}(\mathcal{L})} V^*_{T(\ear)}(S) \tag{$\mathcal{P}_\mathcal{B}$}\label{eq:first_stage} 
\end{equation}
The fact that the UIP problem is hard to solve makes \eqref{eq:first_stage} a highly complex combinatorial problem. Our bundling and pricing problem is thus a bi-level optimization problem where the first level selects a set of options and the second level dynamically assigns prices to available options at each time period. 

Section \ref{section:DP} studies the optimal solution of the single-period pricing problem, the asymptotic optimal revenue and prices for the UIP problem as the demand grows large, and the structural properties of the optimal set of options for various demand levels. In Section \ref{section:bounds}, we present an upper and a lower bound of the UIP problem by considering best and worst-case scenarios in terms of option availability. We then introduce a forward availability approximation framework that enables a tractable lower bound to the UIP problem that is asymptotically optimal in the case of the MNL. Section \ref{section:Bunlding} leverages those bounds to develop a column-generation based bundling algorithm and a fast bundling heuristic designed for large sets of items. We compare the performance of our bounds and bundling algorithms in Section \ref{section:ExperimentalResults}. In addition, we use our heuristic in a dynamic setting where items with distinct item expiration times are successively supplied to the retailer and the set of options is updated during the season. We provide final thoughts and conclude in Section \ref{section:conclusion}. All proofs are provided in the Appendix.\hfill\\
\section{UIP Problem Under MNL Model} \label{section:DP}

In this section, we gain insights into the structure of optimal prices and set of options in the specific case where the choice function is the one derived from the MNL: $f = \rho$.

\subsection{Optimal Pricing and Bundling} \label{subsec:opt_pricing_bundling}
Problem \eqref{eq:inner_problem} is at the core of the UIP problem. As explained in \cite{Optimizing_Multinomial_Logit_Profit_Functions}, problems of the form of \eqref{eq:inner_problem} with $g = \rho^\omega$ are not concave in the price vector $p$. However, under the MNL model, the mapping from a price vector $p$ to a probability vector $\rho$ is bijective as shown in Lemma \ref{lemma:price_proba_bijection} in the Appendix. \cite{Demand_Management_Inventory_Control_Substitutable_Products} point out that this problem in concave in $\rho$. We take advantage of this bijection by maximizing over the probability vector $\rho$ instead of the price vector $p$. Lemma \ref{lemma:inner_problem_solution} derive a near closed-form solution to the optimal expected revenue, prices, and choice probabilities of \eqref{eq:inner_problem}. These expressions use the Lambert $W$ function introduced in \cite{Lambert_W_Function} and defined as the reciprocal function of the mapping $z \mapsto ze^z$.

\begin{lemma}\label{lemma:inner_problem_solution}
For any $S \in \mathcal{P}(\mathcal{L})$ and $\Delta \in \R^{|S|}$, the optimal solution to \eqref{eq:inner_problem} with $g = \rho^\omega$ is given by:
\begin{equation*}
    R(\rho^\omega, S, \Delta) = -\frac{\Gamma}{\beta_p},
\end{equation*}
where $\Gamma \triangleq  W\left(\sum_{i \in S} e^{q_i^\omega + \beta_p \Delta_i - 1}\right)$.\\
Moreover, the associated optimal prices and choice probabilities are:
\begin{align*}
    p^*_i = \Delta_i - \frac{1+\Gamma}{\beta_p}, \quad \forall i \in S
    \hspace{1cm} \text{and} \hspace{1cm}
    \rho^*_i = \frac{\Gamma}{1+\Gamma} \frac{e^{q_i^\omega + \beta_p \Delta_i}}{\sum_{j\in S} e^{q_j^\omega + \beta_p \Delta_j}}, \quad \forall i \in S.
\end{align*}
\end{lemma}

Using the fact that the problem \eqref{eq:inner_problem} is concave in $\rho$ and equating the partial derivatives with respect to $\rho$ to 0 gives the above results. Thus, the optimal price of an option is equal to its marginal value marked up by a quantity that is identical for every option. Further, this near closed-form solution to \eqref{eq:inner_problem} allows us to efficiently solve the maximization problem for each state of the UIP. While this closed-form solution isn't sufficient to solve the UIP problem to optimality due to the curse of dimensionality (see Section \ref{sec:formulation}), it will prove useful in developing approximate solutions to this problem in Section \ref{section:bounds}.

Before exploring these solutions, we are interested in studying the asymptotic dynamics of the exact dynamic program. Specifically, for a fixed set of options $S \in \mathcal{P}(\mathcal{L})$, we study the asymptotic behavior of the value function, optimal prices, and choice probabilities as the demand $\ear$ grows. Let us define the aggregated quality $\kappa_i \triangleq \ln \E_{\rv}\left[e^{q_i^\rv} \right]$ for all $i \in S$ and the mapping $\kappa: S \mapsto \sum_{i \in S} \kappa_i$. The aggregated quality can be interpreted as a weighted RealSoftMax of the perceived quality of items in $S$ with weights corresponding to the probability mass of arrival types. Theorem \ref{theorem:asymptotic_costs} gives a closed-form expression of the asymptotic dynamics of the UIP problem as the demand grows.

\begin{theorem} \label{theorem:asymptotic_costs}
Let $S \in \mathcal{P}(\mathcal{L})$. Then, under the MNL model and as $\ear$ increases:
$$V^*_{T(\ear)}(S) - v_{T(\ear)}(S) = O\left(\frac{\ln^{|S|}(\ear)}{\ear}\right),$$
where $v_{T(\ear)}(S) \triangleq \frac{-1}{\beta_p}\left( |S|(\ln (\ear) - 1) + 
\kappa(S) \right)$ is the asymptotic optimal revenue.\\
In addition, for all $i \in S$ and $\omega \in \Omega$, the optimal price and probability of choice respectively converge to:
\begin{align*}
     p_{T(\ear),i}^{*,\omega} &= \frac{-1}{\beta_p} \left(\ln(\ear) + \kappa_i \right) + O\left(\frac{\ln^{|S|}(\ear)}{\ear}\right),\\
     \rho_{T(\ear),i}^{*,\omega} &= \frac{1}{\ear}e^{q_i^\omega - \kappa_i}+O\left(\frac{\ln^{|S|}(\ear)}{\ear^2}\right).
\end{align*}
as $\lambda$ tends to $\infty$.
\end{theorem}

This theorem follows because each operator $F_t: \epsilon_t \to \epsilon_{t+1}$ (where $\epsilon_t = V^*_t(S) - v_t(S)$ is the difference between the exact and asymptotic optimal revenue) are nonexpansions with fixed points $\epsilon^*_t$ converging to 0 and such that $F_t(\epsilon) - \epsilon$ is lower bounded for all $t$ by a constant independent of $t$. This implies that the sequence $(\epsilon_t)_{t\in \N}$ is bounded. As each $F_t$ restricted to the finite range of $(\epsilon_t)_{t\in \N}$ is a contraction with a factor $\gamma_t$ small enough such that $\prod_{t=1}^\infty \gamma_t = 0$, we obtain $\lim_{t \to \infty} \epsilon_t = 0$.

This theorem reveals that the asymptotic optimal expected revenue is logarithmic in the demand and do not depend on the salvage values. When arrivals (demand) far exceeds the number of options, that is $\ear \gg |S|$, Theorem \ref{theorem:asymptotic_costs} provides a closed-form expression of near-optimal prices, which remarkably does not depend on the customer type. Moreover, since $\mathbb{E}_\rv\left[e^{q_i^\rv - \kappa_i}\right] = 1$, we observe that the asymptotic per period probability of choice is, in expectation, identical for each item and equal to $1/\ear$. This is consistent with the run-out rate described in \cite{Dynamic_Pricing_Strategies_Multiproduct} with a capacity of one. Finally, note that $\kappa_i$ is a weighted RealSoftMax across customer types of the quality vector $q_i$. The fact that the asymptotic optimal expected revenue is expressed in terms of $\kappa(S)$ shows that the seller can benefit greatly from customer types with diverse preferences such that different items are best perceived by different customer types.

Let us consider the difference in revenue between set $S_1$ and $S_2$ be $I_t(S_1, S_2) \triangleq V_{t}(S_1) - V_{t}(S_2)$. We then have:
$$I_{T(\ear)}(S_1, S_2) = \frac{-1}{\beta_p}\left[\left(|S_1|-|S_2|\right) (\ln \ear-1) + \kappa(S_1) - \kappa(S_2)\right] + o(1),$$
implying that for a large enough demand, the expected revenue of $S_1$ exceeds the one of $S_2$ if and only if $|S_1| > |S_2|$, or $|S_1| = |S_1|$ and $\kappa(S_1) > \kappa(S_2)$. This allows us to order sets of options with respect to their asymptotic performances. It also implies that the revenue improvements from bundling become negative as the demand tends to infinity. Indeed, let $S_0$ be the set of options containing single items only. Then $S_0$ has a greater cardinality than any set in $\mathcal{P}(\mathcal{L})$ containing at least one bundle and it follows that there exists $\ear$ for which the expected revenue with $S_{0}$ exceeds those of any set including one or multiple bundles. Specifically, the revenue improvement of set $S \in \mathcal{P}(\mathcal{L})$ over $S_0$ is:
\begin{equation} \label{eq:revenue_improvement}
I_{T(\ear)}(S, S_0) = \frac{-1}{\beta_p}\left[\kappa(S) - \kappa(S_0) -\left(\sum_{i\in S}(\#i-1)\right) \cdot (\ln \ear -1)\right] + o(1),
\end{equation}
where $\#i$ denote the cardinality of option $i$. Consequently, the revenue improvement of any set containing a bundle becomes negative for $\ear$ large enough, as each bundle of size $n$ included in $S$ incurs a penalty of $(n-1)\cdot\left(\ln \lambda - 1\right)/\beta_p$ on the expected revenue. Interestingly, if $\kappa(S)$, the total aggregated quality of set $S$, is significantly larger than $\kappa(S_0)$, the total aggregated quality of $S_0$, then set $S$ can yield greater revenue than set $S_0$ even for $\ear$ relatively large. 

On the other hand, when $T(\ear)=1$ i.e. when the demand is extremely small, the marginal value of each option equates to its salvage value and we obtain with Lemma \ref{lemma:inner_problem_solution}:
$$I_1(S,S_0) = - \frac{\arate}{\beta_p}\mathbb{E}_\rv\left[W\left(\sum_{i \in S} e^{q_i^\rv + \beta_p \xi_i-1}\right)-W\left(\sum_{i \in S_0} e^{q_i^\rv + \beta_p \xi_i-1}\right)\right].$$
Hence the set maximizing $I_1(\cdot,S_0)$ can be drastically different to the one maximizing $I_{T(\ear)}(\cdot,S_0)$ for $\ear$ large. When $\ear$ is small, it is better to favor options with small salvage values and highest quality as they are the key contributors to large expected revenue. In contrast, when $\ear$ is large, the expected revenue is marginally affected by the salvage values and the optimal set of options prioritizes a large total aggregated quality and a small number of bundles.

Let us define $\Delta \kappa(S) \triangleq \kappa(S) - \kappa(S_0)$ as a measure of the complementarity of a set of options, given that $\Delta \kappa(S)>0$ implies that $S$ contains options whose total aggregated quality is greater than that of unbundled items. Equation \eqref{eq:revenue_improvement} gives us a first-order condition for a set $S$ to yield higher revenue than the set of individual items $S_0$:
\begin{equation}\label{eq:condition}
    \Delta \kappa(S) \geq (\ln \ear -1) \cdot \sum_{i\in S}(\#i-1).
\end{equation}

\begin{figure}[h]
    \centering
    \begin{subfigure}{0.49\linewidth}
        \centering
        \includegraphics[width=\linewidth]{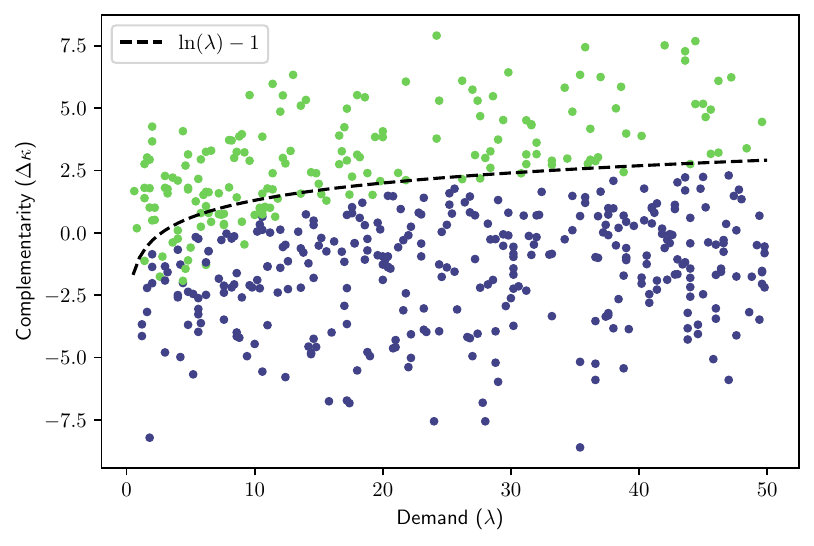}
        \caption{Bundle of size 2}\label{fig:improvement_condition_2}
    \end{subfigure}
    \hfill
    \begin{subfigure}{0.49\linewidth}
        \centering
        \includegraphics[width=\linewidth]{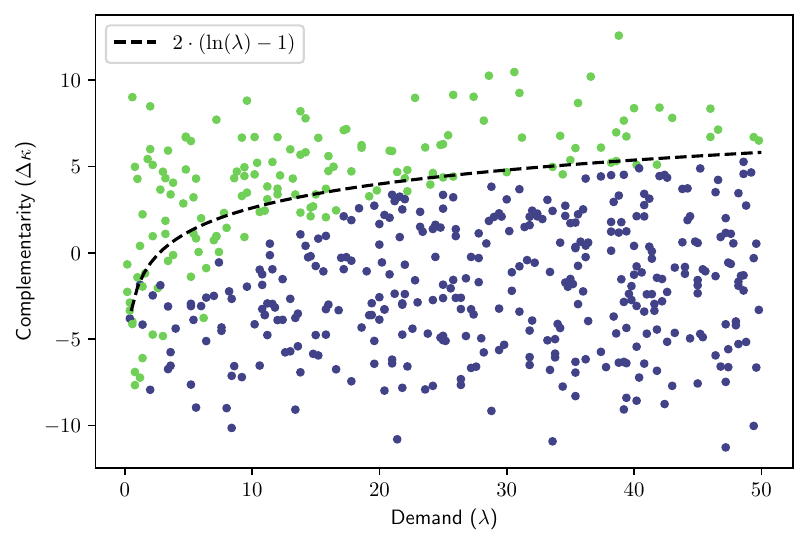}
        \caption{Bundle of size 3}\label{fig:improvement_condition_3}
    \end{subfigure}
    \caption{First-order condition and positive (green)/negative (blue) revenue improvement of sampled sets.}
    \label{fig:improvement_condition}
\end{figure}

Figure $\ref{fig:improvement_condition}$ compares condition \eqref{eq:condition} with the observed revenue improvements on synthetic data: each point corresponds to a randomly sampled set of 5 items, one bundle, and a demand $\ear$. Here, we define the perceived of a bundle as the sum of the perceived quality of each item it contains, multiplied by a random variable following $\mathcal{U}(\frac{1}{2}, \frac{3}{2})$ to reflect substitutability and complementarity. We compare the performance of the set $S$ containing the bundle and the individual items not included in it to the performance of set $S_0$. We consider sets of options containing a single bundle of size $n=2$ or 3, in which case condition \eqref{eq:condition} simplifies to $\Delta \kappa \geq (n-1)\cdot \left(\ln(\lambda)-1\right)$.

As expected, the accuracy of condition \eqref{eq:condition} improves as $\ear$ increases. Since condition \eqref{eq:condition} ignores the expected revenue from salvaging that is included in the big O of equation \eqref{eq:revenue_improvement}, this condition tends to be conservative in the sense that a set satisfying it with equality has a positive revenue improvement with high probability.

We note that it is possible to improve the expected revenue by bundling substitute items ($\Delta \kappa < 0$) when the demand is small. Even if they have a negative impact on the quality of the options offered, bundles can reduce the expected number of unsold items at the end of the season, resulting in improved expected revenue. Conversely, one can deteriorate the expected revenue by bundling complementary items ($\Delta \kappa > 0$) when the demand is large. In this case, reducing the number of options via bundling may overshadow the increased quality and lead to lower expected revenues, as high demand allows the retailer to fully exploit each option through dynamic pricing.

\subsection{Relation Between Revenue and Utility}

Let $\displaystyle U_{T(\ear)}(S) = \ln \sum_{t=1}^{T(\ear)} \sum_{i \in S} \mathbb{E}_\rv\left[e^{u_{t,i}^{\rv,*}}\right]$ be the cumulative aggregated utility, where $u_{t,i}^{\omega, *} = q_i^\omega + \beta_p p_{t,i}^{\omega,*}$ is the utility of option $i$ for customer type $\omega$ under optimal prices when the set of options $S$ is offered. This quantity is the RealSoftMax (defined as the logarithm of the sum of the exponentials of the arguments) of all the utilities $u_{t,i}^{\omega,*}$ offset by $\ln p_\rv(\omega)$ with $p_\rv(\omega)$ being the probability of arrival of type $\omega$: $U_{T(\ear)}(S) = \text{RealSoftMax}\left(\{ u_{t,i}^{\omega,*} + \ln p_\rv(\omega) \;|\;t\in [T(\ear)],\; i\in S,\; \omega \in \Omega \}\right)$. It can be interpreted as a measure of the attractiveness of a set of options over the entire season.

\begin{theorem} \label{theorem:revenue_utility} For any $\ear\geq 0$, optimal solutions to \eqref{eq:first_stage} also maximize the cumulative aggregated utility:
$$\argmax_{S \in \mathcal{P}(\mathcal{L})} V^*_{T(\ear)}(S) = \argmax_{S \in \mathcal{P}(\mathcal{L})} U_{T(\ear)}(S)$$
\end{theorem}

That is, under an optimal pricing policy, maximizing revenue through bundling is equivalent to maximizing the cumulative aggregated utility of customers. As discussed in subsection \ref{subsec:opt_pricing_bundling}, the revenue improvement from bundling becomes negative as the demand tends to infinity, resulting in $S_0$ being optimal for $\ear$ large enough. Consequently, the set of options maximizing the cumulative aggregated utility is $S_0$ for $\ear$ large enough. Therefore, revenue maximization through bundling does not come at the expense of reduced attractiveness of options under optimal pricing.

Intuitively, bundling can improve expected revenue in two ways. First, by increasing the quality of options. If the quality of a bundle is greater than the sum of the qualities of the items it contains, bundling can increase the total quality of options. This higher quality leads to higher prices without increasing the expected number of unsold items. Second, by reducing the number of options. This results in increased demand per option allowing a price increase small enough to reduce the expected number of unsold items. The impact of this latter mechanism increases as demand decreases. Both of these mechanisms allow optimal prices to increase expected revenue and the attractiveness of options simultaneously, leading to the conclusion stated in Theorem \ref{theorem:revenue_utility}.

Since the revenue improvement from bundling need not be positive, we want to obtain tractable and accurate approximations of the expected revenue for a given set of options under a finite horizon in order to make informed bundling decisions. In the next section, we derive bounds on the optimal objective of \eqref{eq:bellman_equation} for a given demand, set of options, and choice function. We will then use these bounds to select sets of options: a lower bound allows us to select the set that maximizes the guaranteed expected revenue, and an upper bound allows us to compute the corresponding optimality gap.\hfill\\
\section{Bounds on Optimal Expected Revenue} \label{section:bounds}

In this section, we develop bounds on the expected revenue under an arbitrary choice function $f$. Specifically, we develop three bounds on $V^*_{T(\ear)}(\bar{S})$ for a fixed demand $\ear$ and a given set of options $\bar{S} \in \mathcal{P}(\mathcal{L})$ with cardinality $N$. We first derive two bounds based on backward recursion. In addition to providing a price trajectory with desirable properties, these simple bounds exhibit better performance than the fluid and static approximations discussed in the Appendix [\ref{apx:fluid_approx} \& \ref{apx:static_approx}]. We then present a lower bound based on a forward approximation of option availability throughout the season.

\subsection{Backward Approximations} \label{subsection:backward}

Since choice functions are decreasing monotone with respect to $S$ (\textit{p.4} in Definition \ref{def:choice_funciton}), we have that for a fixed price vector $p$, the probability of choosing option $i$ is highest when the set of available options is $\{i\}$. Hence, estimating the probability that every option is selected as if each option was the only available one overestimates the per-step expected revenue. This best-case scenario corresponds to replacing each probability $f_i^\omega(p,S)$ with $f_i^\omega(p,\{i\})$ in problem \eqref{eq:bellman_equation}. We obtain a modified dynamic program that is separable in each option, and whose value function is expressed as the sum of the contribution of each option.\\
Formally, let $r \in \R^{T(\ear) \times N}$ be such that for all $i \in \bar{S}$:
\begin{align}
    r_{0,i} &= \xi_i \nonumber\\
    r_{t,i} &= r_{t-1,i} + \arate \cdot \E_\rv \left[ \max_{p \in \R} f_i^\rv(p, \{i\}) \left(p - r_{t-1,i}\right) \right], \quad \forall t \in \{1,\hdots,T(\ear)\}  \label{eq:u_recursion}
\end{align}
We call $r_{t,i}$ the individual upper bound of option $i$. Intuitively, $r_{t,i}$ corresponds to the optimal expected revenue if option $i$ was the only option available during a season of length $t$. When other options are available, customers are no more likely to accept this option (assuming its price remains the same), making $r_{t,i}$ an upper bound on the expected revenue induced by option $i$ when multiple options are available. We define $V^U_t\left(S\right) = \sum_{i \in S} r_{t,i}$ for all $t \in \{0,\hdots,T(\ear)\}$. Computing $V^U_{T(\ear)}\left(S\right)$ requires solving $|\Omega| \cdot |S| \cdot T(\ear)$ problems of the form of $\eqref{eq:inner_problem}$. Note however that the sequences $(r_{t,i})_{t \in [T(\ear)]},\; i=1,...,|\bar{S}|$ are independent and thus can be computed in parallel. The optimal value of this modified dynamic program yields an upper bound on the optimal value of the original problem that is asymptotically optimal in the case of the MNL model.

\begin{proposition} \label{prop:upper_bound}
For all $S \subseteq \bar{S}$ and $t \in \{0,\hdots,T(\ear)\}$, $V_t^*(S) \leq V^U_t(S)$. In addition, $\lim_{\ear \to \infty}  V_{T(\ear)}^*(S) - V^U_{T(\ear)}(S) = 0$ when $f = \rho$.
\end{proposition}

The first result of this proposition is obtained by an induction argument. The second follows from applying Theorem \ref{theorem:asymptotic_costs} to each sequence $(r_{t,i})_{t \in [T(\ear)]}$ to obtain $\lim_{\ear \to \infty} r_{{T(\ear)},i} - v_{T(\ear)}(\{i\}) = 0$. We obtain the desired result by summing over the options included in $S$.

This result contrasts with Proposition \ref{prop:fluid_approx_asymptotical} given in the Appendix, which states that the fluid approximation of the UIP problem is not asymptotically optimal as $\lambda$ tends to infinity. This supports the idea that $V^U$ is better suited to the UIP problem than the fluid approximation.

We say a price trajectory $\tau \in \R^{T(\ear) \times N \times |\Omega|}$ is \textit{homogeneous} if prices are independent of the customer type i.e. $\tau_{t,i}^\omega = \tau_{t,i}^{\omega'}$ for all $(\omega, \omega') \in \Omega^2$, in which case we drop the superscript for brevity. In addition to being a valuable technical property, homogeneous prices may be desirable or even required by the retailer for reasons of fairness or feasibility. While solving recursion \eqref{eq:u_recursion}, we obtain a price trajectory consisting of the optimal price for each triple $(t,i,\omega)$. Taking the maximum over $\omega$ gives a non-increasing and homogeneous price trajectory. 

\begin{proposition} \label{prop:prices_upper_bound}
The homogeneous price trajectory $\tau^U$ defined as:
\begin{align}
    \tau^U_{t,i} &\triangleq \max_{\omega \in \Omega}\left\{\argmax_{p \in \R} f_i^\omega(p, \{i\}) \left(p - r_{t,i}\right)\right\},\quad \forall t\in \{0,\hdots,T(\ear)\}, \quad \forall i \in \bar{S} \label{eq:tau^U}
\end{align}
is non-increasing as $t$ decreases for all $i \in \bar{S}$.
\end{proposition}
This proposition follows since $r_{t,i}$ is increasing in $t$ and properties $\textit{p.3}$ and $\textit{p.5}$ in definition \ref{def:choice_funciton} ensure that $\partial \tau^U_{t,i}/ \partial r_{t,i} \geq 0$. In the specific case of the MNL model, we obtain with Theorem \ref{theorem:asymptotic_costs} that $\tau_{T(\ear)}^U$ converges to the optimal price vector as $\ear$ tends to infinity.

Analogously, one can obtain a lower bound $V^L$ on the expected revenue by underestimating the probability that an option is selected in problem \eqref{eq:bellman_equation}. To obtain such a bound we replace $f_i^\omega(p,S)$ with $f_i^\omega(p, \bar{S})$. This corresponds to always estimating the selection probabilities as if all other options were available. The full description of this bound is available in the Appendix [\ref{apx:backward_lower_bound}]. However, contrary to the upper bound, the price trajectory induced by the lower bound need not be non-increasing as $t$ decreases when $N > 1$. Indeed, as explained and illustrated by \cite{DP_Inventory_Control_Substitute_Products}, optimal prices under the MNL model exhibit a complex behavior when multiple options are offered to arriving customers: the main effects driving prices are the difference in quality between available options, inventory scarcity, and the interplay of these two effects. In our case, low inventory levels lead to increasing prices for low-quality options when the demand is small and other higher-quality options are available. Furthermore, although more flexible than the static lower bound [\ref{apx:static_approx}], this bound may be loose as it considers an extreme worst case in terms of option availability. In the next section, we develop a much more useful lower bound that considers real-valued availabilities.

\subsection{Deterministic Forward Approximation} \label{subsection:DFA}

The two previous bounds maintain an estimate of the marginal value of each option from the final to the initial state, presuming the availability of other options. We will now build on this idea and maintain an estimate of the availability of each option from the initial to the final state under a fixed trajectory of prices.

In each state of the dynamic program \eqref{eq:bellman_equation}, only a subset of the initial set of options $\bar{S}$ is available. We represent the probability of these availabilities with a real-valued vector $a \in [0,1]^N$. Let $\tilde{f}:\R^{|\Omega|\times N} \times [0,1]^{N} \to [0,1]^{|\Omega|\times N}$ be an extended choice function as in definition \ref{def:choice_funciton}. The mapping $\Tilde{f}_i^\omega$ gives the probability that option $i$ is selected by a customer of type $\omega$ with price vector $p$ conditioned to the availability of option $i$ and the availability of other options given by vector $a$. In the specific case of the MNL model, a natural extension of $\rho$ is $\tilde{\rho}$ satisfying:
\begin{equation}
    \Tilde{\rho}_i^\omega(p, a) = \frac{e^{q_i^\omega + \beta_p p_i}}{\displaystyle 1 + e^{q_i^\omega + \beta_p p_i} + \sum_{j \in \bar{S} \backslash\{i\}} a_j e^{q_j^\omega + \beta_p p_j}}, \quad \forall \omega \in \Omega, \quad \forall i \in \bar{S}. \label{eq:tilde_rho_def}
\end{equation}
For any set of options $S \subseteq \bar{S}$, let $a(S) \in [0,1]^N$ be the availability vector of set $S$ i.e. $a_i(S) = \mathds{1}_{\{i \in S\}}$ for all $i \in \bar{S}$. Using this specific vector, we have $\rho_i^\omega(\cdot, S) = \Tilde{\rho}_i^\omega(\cdot, a(S))$ for all $\omega \in \Omega$ and $i \in S$. Hence, $\Tilde{\rho}$ is indeed an extension of $\rho$ that allows for real-valued availabilities.

The approximation we construct below uses $\Tilde{f}_i^\omega$ to compute a trajectory of real-valued availability vectors, given a fixed price trajectory $\tau$. 
Let $A^* \in [0,1]^{T(\ear) \times N}$ be the matrix containing the probabilities of options' availability throughout the season under the price trajectory $\tau$. The actual availability of the options can then be modeled by the matrix of random variables $A$ satisfying $A_{t,i} \sim Ber(A^*_{t,i})$ for all $t$ and $i$. In particular, $A_t$ is a vector of random variables whose realization represents the availability of options in $\bar{S}$ at time $t$. Let $Y_{t,i}$ be the event \say{option $i$ is selected at time $t$, conditioned to the arrival of a customer at time $t$}. We can express the expected revenue obtained under the price trajectory $\tau$ as:
\begin{equation}
    V_{T(\ear)}(\bar{S}, \tau) \triangleq \sum_{i \in \bar{S}} \underbrace{\mathbb{P}(A_{0,i} = 1) \xi_i}_{(a_i)} + \sum_{1<t\leq T(\ear),i\in \bar{S}} \underbrace{\arate \mathbb{P}(A_{t,i} = 1) \E_\rv \left[\mathbb{P}(Y_{t,i} = 1 \;|\; \rv, A_{t,i} = 1)\right] \tau_{t,i}}_{(b_{t,i})} \label{eq:expected_revenue_tau}
\end{equation}
where $(a_i)$ is the expected revenue from salvation of option $i$, and $(b_{t,i})$ is the expected revenue from selling option $i$ at time $t$.
Conditioned to the availability of option $i$ and the arrival of a customer of type $\omega$ at time $t$, the probability that option $i$ is selected at time $t$ with price vector $\tau_{t}$ is:
\begin{equation*}
\mathbb{P}(Y_{t,i} = 1 \;|\; \omega, A_{t,i} = 1) = \E_{A}\left[\tilde{f}_i^\omega(\tau_t, A_t)\right]
\end{equation*}
We can then reformulate the total expected revenue as:
\begin{equation} \label{eq:total_costs}
    V_{T(\ear)}(\bar{S}, \tau) = \sum_{i\in \bar{S}} A^*_{0,i} \xi_i + \arate \sum_{1<t\leq T(\ear),i\in \bar{S}} A^*_{t,i}  \E_{\rv,A}\left[\tilde{f}_i^\rv(\tau_t, A_t)\right] \tau_{t,i}
\end{equation}
Initially, every option is available so $A^*_{T(\ear),i} = 1$ for all $i \in \bar{S}$. The availability in the following periods is then computed in a forward manner. Indeed, option $i$ is available at time $t-1$ if and only if option $i$ is available at time $t$ and either no customer arrives at time $t$ or a customer arrives but does not choose option $i$. We thus obtain the following recursion:
\begin{align}
    A^*_{t-1,i} &= \mathbb{P}(A_{t,i} = 1) \left((1-\arate)+\arate\E_{\rv}\left[\mathbb{P}(Y_{t,i} = 0 \;|\; \rv, A_{t,i} = 1)\right]\right) \nonumber\\
    &= A^*_{t,i} \cdot \left(1-\arate \E_{\rv,A}\left[\tilde{f}_i^\rv(\tau_t, A_t)\right]\right) \label{eq:availability_recursion}
\end{align}
Note however that the random variable $A_t$ can take $2^N$ values so computing the expectation over $A$ is costly. Instead, we will recursively compute an estimate $\hat{A}$ of $A^*$ by replacing $\E_A\left[\tilde{f}_i^\omega(\tau_t, A_t)\right]$ with $\tilde{f}_i^\omega(\tau_t, \E_A\left[A_t\right])$. Since $\tilde{f}_i^\omega(\tau_t, \cdot)$ is non-increasing and convex (property \textit{p.4} in definition \ref{def:choice_funciton}), Jensen's inequality tells us that the estimate $\hat{A}$ is an overestimate of $A^*$. Consequently, the approximation overestimates the expected time needed to sell options. If $\tau$ is homogeneous and non-increasing, substituting $A^*$ with $\hat{A}$ in $\eqref{eq:total_costs}$ yields a lower bound on the exact total revenue.\\
Formally, we define the Deterministic Forward Approximation (DFA) for a choice function $f$ and price trajectory $\tau$ as:
\begin{equation}
    V_{T(\ear)}^\text{DFA}\left(\bar{S}, \tau\right) \triangleq \sum_{i\in \bar{S}} \hat{A}_{0,i}\xi_i + \arate \sum_{1<t\leq T(\ear),i\in \bar{S}} \hat{A}_{t,i} \E_\rv\left[\tilde{f}_i^\rv\left(\tau_t, \hat{A}_{t}\right)\right] \tau_{t,i} \label{eq:def_DFA}
\end{equation}
where $\hat{A} \in [0,1]^{T(\ear) \times N}$ is such that for all $i \in \bar{S}$:
\begin{align}
    \hat{A}_{T(\ear),i} &= 1 \nonumber\\
    \hat{A}_{t-1,i} &= \hat{A}_{t,i} \cdot \left(1-\arate\E_\rv\left[\tilde{f}_i^\rv\left(\tau_t, \hat{A}_{t}\right)\right]\right), \quad \forall t \in \{1,\hdots,T(\ear)\} \label{eq:estimated_availability_recursion}
\end{align}
As intuitively explained above, the DFA provides a lower bound on the optimal total revenue if $\tau$ is a homogeneous and non-increasing price trajectory. We note that the DFA need not be a lower bound if $\tau$ is not homogeneous, even if each customer type is given a non-increasing price trajectory.
\begin{theorem}\label{theorem:DFA}
Given any homogeneous and non-increasing price trajectory $\tau$, the DFA lower bounds the optimal expected revenue i.e. $V_{T(\ear)}^\text{DFA}\left(\bar{S},\tau\right) \leq V_{T(\ear)}^*\left(\bar{S}\right)$.
\end{theorem}
The proof first requires showing that $A^* \leq \hat{A}$, which can be done using induction. By optimality of $V_{T(\ear)}^*\left(\bar{S}\right)$, we have $V_{T(\ear)}^*\left(\bar{S}\right) \geq V_{T(\ear)}\left(\bar{S}, \tau\right)$ for any $\tau$. We then remark that the term $(b_{t,i})$ in equation \eqref{eq:expected_revenue_tau} is equal to $\left(A_{t,i}^* - A_{t-1,i}^*\right)\tau_{t,i}$. As $A^* \leq \hat{A}$ and since $\tau$ is non-increasing as $t$ decreases, we add $\sum_{i \in \bar{S}} \sum_{t=0}^{{T(\ear)}-1} \left(\hat{A}_{t,i}-A^*_{t,i}\right)\left(\tau_{t,i}-\tau_{t+1,i}\right) \leq 0$ to $V_{T(\ear)}\left(\bar{S}, \tau\right)$ and obtain $V_{T(\ear)}^\text{DFA}\left(\bar{S},\tau\right)$ after canceling multiple terms, finally establishing that \eqref{eq:def_DFA} is a lower bound.

To obtain the tightest lower bound, we want to find the homogeneous and non-increasing price trajectory $\tau$ that maximizes $V_{T(\ear)}^\text{DFA}(\bar{S}, \cdot)$. A good candidate is $\tau^U$ defined in \eqref{eq:tau^U}. In addition to being homogeneous and non-increasing as shown in Proposition \ref{prop:prices_upper_bound}, this price trajectory has the convenience of being readily available once the upper bound $V_{T(\ear)}^U(\bar{S})$ has been calculated. Further, $\tau^U_{T(\ear)}$ converges to the optimal prices as $\ear$ tends to infinity under the MNL model. Since $V_{T(\ear)}^\text{DFA}(\bar{S}, \cdot)$ is differentiable, one can use gradient descent with $\tau^U$ as starting point to find a local minimum $\tau^*$. However, this procedure is relatively computationally expensive, making it suitable for refining a bound, but not for bounding a large number of sets of options. For this reason, we will mostly use $\tau^U$ to obtain DFA bounds, and we now refer to $V_{T(\ear)}^\text{DFA}\left(S,\tau^U\right)$ as $V_{T(\ear)}^\text{DFA}(S)$ unless stated otherwise.

The following theorem states that, under the MNL model, the DFA approximation is asymptotically optimal with this choice of price trajectory, and gives an asymptotic performance bound for $V^\text{DFA}$ and $V^U$.
\begin{theorem} \label{theorem:performance_bound}
For $f = \rho$, the DFA evaluated with $\tau^U$ is asymptotically optimal as $\ear \to \infty$ and we have for all $S \in \mathcal{P}\left(\mathcal{L}\right)$:
    \begin{align*}
   \frac{V_{T(\ear)}^\text{DFA}(S) - V_{T(\ear)}^*(S)}{V_{T(\ear)}^*(S)} = O\left(\frac{1}{\ear}\right),
   \quad \text{and} \quad
   \frac{V_{T(\ear)}^U(S) - V_{T(\ear)}^*(S)}{V_{T(\ear)}^*(S)} = O\left(\frac{1}{\ear}\right)
\end{align*}
\end{theorem}
Let $\hat{A}(T)$ be the availability trajectory satisfying $A_{T,i} = 1$ for all $i \in \bar{S}$ and the recursion in \eqref{eq:estimated_availability_recursion}. Theorem \ref{theorem:performance_bound} leverages the fact that $\hat{A}(T)/\mathbb{P}(Y_{T,i} = 0) \leq \hat{A}(T-1)$ with Lemma \ref{lemma:telescopic_inequality} to lower bound $V_T^\text{DFA}(S)$ as a function of $V_{T-1}^\text{DFA}(S)$. The difference $V_{T(\ear)}^\text{DFA}(S) - V_{T(\ear)}^U(S)$ is then examined, and we obtain using the results of Theorem \ref{theorem:asymptotic_costs} that $V_{T(\ear)}^\text{DFA}(S) - V_{T(\ear)}^U(S) = O(\ln(\ear)/\ear)$ and $V_{T(\ear)}^*(S) = \Theta(\ln(\ear))$. Combined with the inequalities of Proposition \ref{prop:upper_bound} and Theorem \ref{theorem:DFA} we obtain the desired performance bounds, implying that the DFA is asymptotically optimal when evaluated with $\tau^U$.

This result makes the DFA an appropriate choice of lower bound for the UIP problem. In contrast, the static lower bound need not be asymptotically optimal when $\ear$ tends to infinity, as shown in Proposition \ref{prop:static_approx_asymptotical}.

\begin{figure}[h]
    \centering
    \begin{subfigure}{0.49\linewidth}
        \centering
        \includegraphics[width=\linewidth]{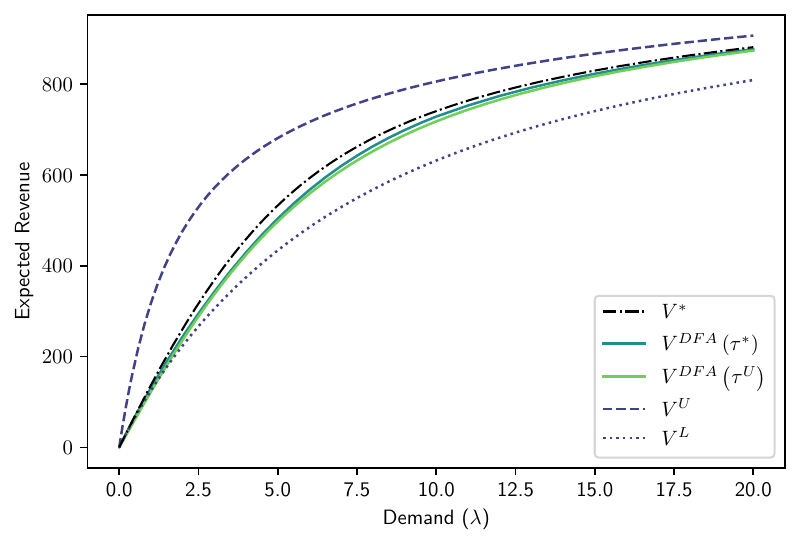}
    \end{subfigure}
    \hfill
    \begin{subfigure}{0.49\linewidth}
        \centering
        \includegraphics[width=\linewidth]{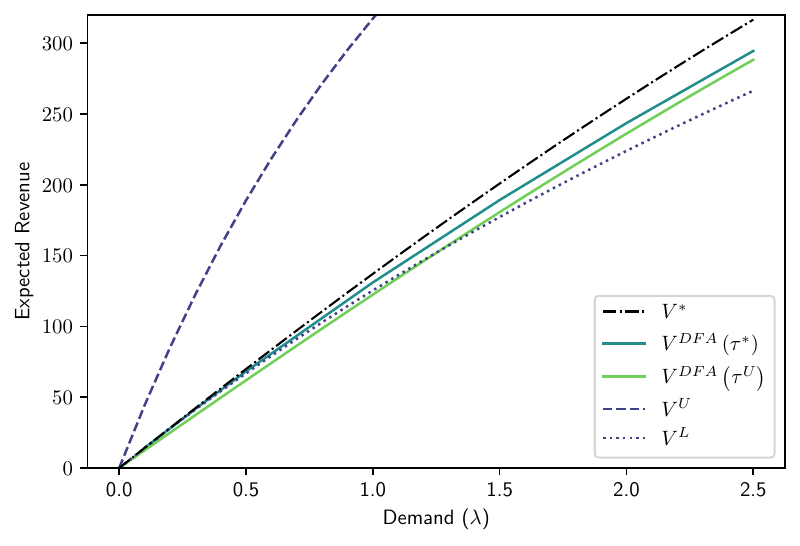}
    \end{subfigure}
    \caption{Performance comparison of the bounds for $L = 5$. (b) is a magnified version of (a) highlighting the performance for $\ear$ taking small values.}
    \label{fig:DFA_illustration}
\end{figure}

Figure \ref{fig:DFA_illustration} illustrates the behavior of the bounds introduced in this section on an example with 5 options, a single type of customer, and no salvage value. We observe that $V^\text{DFA}\left(\tau^U\right)$ is a tighter lower bound than $V^L$ except for small values of $\ear$. In particular, the DFA converges to the exact expected value much faster than the backward lower bounds as $\ear$ increases. Further, we see that $V^\text{DFA}\left(\tau^U\right)$ comes quiteclose to $V^\text{DFA}\left(\tau^*\right)$ while incurring considerably less computational cost.

\vspace{-1.5cm}\hfill\\
\section{Bundling Algorithms} \label{section:Bunlding}

Now that we have established some tools for evaluating the performance of a given set of options, we turn our attention to problem \eqref{eq:first_stage}. Namely, we want to find the set of options $S \in \mathcal{P}(\mathcal{L})$ that maximizes $V_{T(\ear)}^*(S)$ for a given $\ear$. 
To achieve this, we first present a column generation-based algorithm suitable for any arbitrary choice function $f$, building on the bounds developed in Section \ref{section:bounds}. We then provide a simpler and faster heuristic in the specific case of the MNL choice function. 

\subsection{Column-Generation Algorithm}

Leveraging the DFA, we can efficiently and accurately lower bound the expected revenue of a set of options. Consequently, we focus our efforts on identifying the set of options maximizing the DFA instead of the set maximizing the exact expected revenue.

The main challenge in finding the set of options maximizing the DFA is that the cardinality of $\mathcal{P}(\mathcal{L})$ is exponential in $\mathcal{L}$. In addition, $V^\text{DFA}_{T(\ear)}$ is a non-linear set function, and the cardinality of the set of feasible options $\mathcal{O}$ is in the order of $L^{K_b}$, making $|\mathcal{O}|$ large enough to prevent evaluating the DFA of each option in $\mathcal{O}$. We address these challenges by developing an algorithm that efficiently generates bundles and forms sets of options to be evaluated with the DFA. This algorithm leverages the fact that $V^U$ is cost-efficient and that individual upper bounds $r$ are parallelizable. It involves three major steps:
\begin{enumerate}
    \item Compute the individual upper bound $r_{T(\ear),i}$ for all $i \in \mathcal{O}$.
    \item Linearize $V^\text{DFA}_{T(\ear)}$ at $S_0$ (set of options containing individual items only) and form a set partitioning problem whose binary variables denote whether an option is included in the set of options. Relax the binary constraints of this problem and solve it in a column-generation fashion using perturbed reduced costs involving the individual upper bounds $r$.
    \item Obtain up to $N_\text{eval}$ solutions to the set partitioning problem with binary constraints, evaluate their DFA, and select the one achieving the largest value.\\
\end{enumerate}

Recall that the backward upper bound $V^U_{T(\ear)}(S)$ of a set of options $S$ is the sum of the contributions $r_{T(\ear),i}$ of each option it contains, which can be computed in parallel. Therefore, it is sufficient to determine the individual upper bound $r_{T(\ear),i}$ of each option $i$ to retrieve the upper bound of any set of options easily.

Let $S_i$ be the set containing bundle $i$ and all individual items not included in bundle $i$. We linearize the DFA at $S_0$ by letting $\delta_i = V_{T(\ear)}^\text{DFA}\left(S_i\right) - V_{T(\ear)}^\text{DFA}\left(S_0\right)$ be the DFA improvement of option $i$. We then consider the following set partitioning problem:
\begin{align} \label{eq:P_IO}
    \max_{z} & \sum_{i \in \mathcal{O}} z_i \delta_i \tag{$\mathcal{P}_\text{IO}$}\\
    \text{s.t.} & \sum_{i \in \mathcal{O}\;:\;l \in i} z_i = 1, \quad \forall l \in \mathcal{L}\nonumber\\
    &\sum_{i \in \mathcal{O}\;:\;\#i > 1} z_i \leq K_s, \nonumber\\
    &z_i \in \{0,1\}, \quad \forall i \in \mathcal{O} \nonumber
\end{align}
where $\#i$ is the cardinality of option $i$ i.e. the number of items it contains.

Problem \eqref{eq:P_IO} consists of selecting the partition of $\mathcal{L}$ maximizing the total DFA improvement. The first two constraints ensures that the solution of \eqref{eq:P_IO} is a partition of $\mathcal{L}$ with at most $K_s$ bundles, and the last one ensures that the solution is binary. By relaxing the latter, we obtain a linear program whose primal and dual formulations are given below:

\vspace{-0.7cm}

\noindent
\makebox[0.95\linewidth]{%
\begin{minipage}[t]{.4\linewidth}
\begin{align} \label{eq:P_LP}
    \max_{z} & \sum_{i \in \mathcal{O}} z_i \delta_i \tag{$\mathcal{P}_\text{LP}$}\\
    \text{s.t.} & \sum_{i \in \mathcal{O}\;:\;l \in i} z_i = 1, \quad \forall l \in \mathcal{L} \nonumber\\
    &\sum_{i \in \mathcal{O}\;:\;\#i > 1} z_i \leq K_s, \nonumber\\
    &z_i \in [0,1], \quad \forall i \in \mathcal{O} \nonumber
\end{align}
\end{minipage}
\hfill
\raisebox{-0.7cm}{\vrule width 0.4pt depth 4cm}
\hfill
\begin{minipage}[t]{.45\linewidth}
\begin{align} \label{eq:P_D}
    \min_{\mu} & \sum_{l \in \mathcal{L}} \mu_l + \sum_{i \in \mathcal{O}} \mu_i + \mu_{K_s}\tag{$\mathcal{P}_\text{D}$}\\
    \text{s.t.} & \sum_{l\in i} \mu_l + \mu_i + \mu_{K_s} \mathds{1}_{\{\#i>1\}} \geq \delta_i, \quad \forall i \in \mathcal{O} \nonumber\\
    &\mu_i \geq 0, \quad \forall i \in \mathcal{O}\nonumber \\
    &\mu_{K_s} \geq 0 \nonumber
\end{align}
\end{minipage}
}
\vspace{0.7cm}

\noindent where $\mu_l$ is the dual variable of the equality constraint associated with item $l$ in the primal problem, $\mu_i$ the dual variable of the inequality constraint $z_i \leq 1$, and $\mu_{K_s}$ is the dual variable of the constraint $\sum_{i \in \mathcal{O}\;:\;\#i > 1} z_i \leq K_s$. This relaxation allows us to use column-generation to form high-potential bundles. The master problem is defined as \eqref{eq:P_LP} except that the set $\mathcal{O}$ is restrained to the set $\overline{\mathcal{O}}$ initially containing individual items only.

At each iteration, we solve the master problem to obtain the primal and dual optimal solutions $(z^*, \mu^*)$. Ideally, we then add to the master problem the option not already included in $\overline{\mathcal{O}}$ with the largest reduced cost $\overline{\delta}_i = \delta_i - \sum_{l \in i} \mu^*_l - \mu^*_{K_s}$. However, by the definition of $\delta_i$, this would require computing $V_{T(\ear)}^\text{DFA}\left(S_i\right)$ for all $i \in \mathcal{O}$ which is too computationally expensive. Instead, we consider the perturbed reduced cost $\tilde{\delta}_i \triangleq V^U_{T(\ear)}(S_i) - V_{T(\ear)}^\text{DFA}\left(S_0\right) - \sum_{l \in i} \mu^*_l - \mu^*_{K_s}$ for each option, using the backward upper bound instead of the DFA of the set of options $S_i$. Note that $\overline{\delta}_i \leq \tilde{\delta}_i$. Let $i^* = \argmax_{i} \tilde{\delta}_i$. If $\tilde{\delta}_{i^*} > 0$, option $i$ is added to $\overline{\mathcal{O}}$. Otherwise, as $\tilde{\delta}_{i^*} \leq 0$ implies $\max_{i} \overline{\delta}_i \leq 0$, no option has a positive reduced cost and the generation of options stops. Since the number of iterations to obtain non-positive reduced costs can be large, we limit the number of generated options to $N_\text{gen}$.

We then restore the binary constraints and solve \eqref{eq:P_IO} using the restrained set of options $\overline{\mathcal{O}}$. Note that $\delta_i$ simply compares the DFA of $S_i$ and $S_0$. Therefore, solving \eqref{eq:P_IO} over $\overline{\mathcal{O}}$ can produce solutions that suffer from neglecting the interactions between options not included in $S_0$, or, in other words, from treating the DFA as a linear set function. Instead, we generate a set $\mathcal{S}$ of up to $N_\text{eval}$ solutions to \eqref{eq:P_IO} over $\overline{\mathcal{O}}$ by sequentially solving and constraining problem \eqref{eq:P_IO} to prevent selecting any previous solution. We then evaluate $V_{T(\ear)}^\text{DFA}\left(S\right)$ for all $S \in \mathcal{S}$, and select the set achieving the largest value. A pseudocode of the bundling algorithm is given in Algorithm \ref{alg:bld_alg}.

\begin{algorithm}
\caption{Column Generation Algorithm}\label{alg:bld_alg}
\begin{algorithmic} \small
\Require $\mathcal{L}$: set of items, $K_s$: maximum number of bundles, $K_b$: maximum bundle size, $N_\text{gen}$: maximum number of options generated, $N_\text{eval}$: maximum number of sets of options evaluated
\State $\mathcal{O} \leftarrow \left\{(l_1, \hdots, l_k) \in \mathcal{L}^k\;|\; k\leq K_b \text{ and } l_a \neq l_b, \text{ for } a\neq b\right\}$  \Comment{$\mathcal{O}$ is the set of all options of size at most $K_b$}
\For{$i \in \mathcal{O}$}
    \State Compute $\left(r_{t,i} , \tau^U_{t,i}\right)_{t \in [T(\ear)]}$
\EndFor
\State $\overline{\mathcal{O}} \leftarrow \left\{i \in \mathcal{O} \;|\;\#i=1\right\}$
\Do
\State Solve master problem to obtain $(z^*,\mu^*)$
\State $i^* \leftarrow \argmax_{i \in \mathcal{O}\backslash\overline{\mathcal{O}}} \left\{ r_{T(\ear),i} + \sum_{l \notin i} r_{T(\ear),l} - V_{T(\ear)}^\text{DFA}\left(S_0\right) - \sum_{l \in i} \mu_l^* - \mu^*_{K_s} \right\}$
\If{$\tilde{\delta}_{i^*} > 0$}
\State $\delta_{i^*} \leftarrow V_{T(\ear)}^\text{DFA}\left(S_{i^*}\right) - V_{T(\ear)}^\text{DFA}\left(S_{0}\right)$
\State $\overline{\mathcal{O}} \leftarrow \overline{\mathcal{O}} \cup \left\{i^*\right\}$
\EndIf
\doWhile{$\tilde{\delta}_i^* > 0$ and $|\overline{\mathcal{O}}| < |\mathcal{L}| + N_\text{gen}$}
\State Compute $\mathcal{S} =$ top $N_\text{eval}$ solutions of $\left(\mathcal{P}_{IO}\right)$ restrained to $\overline{\mathcal{O}}$
\State \Return $\argmax_{S \in \mathcal{S}} V_{T(\ear)}^\text{DFA}(S)$
\end{algorithmic}
\end{algorithm}

Moreover, since the upper bound $V^U_{T(\ear)}(S)$ of a set $S$ is linear in the individual upper bounds $r_{T(\ear),i}$ of each option it contains, we can find the highest achievable upper bound by solving the following MILP:
\begin{align*}
    Z^* = \max_{z} & \sum_{i \in \mathcal{O}} z_i \cdot r_{T(\ear),i}\\
    \text{s.t.} & \sum_{i \in \mathcal{O}\;:\;l \in i} z_i = 1, \quad \forall l \in \mathcal{L}\\
    &\sum_{i \in \mathcal{O}\;:\;\#i > 1} z_i \leq K_s,\\
    &z_i \in \{0,1\}, \quad \forall i \in \mathcal{O}
\end{align*}
We use the optimal objective $Z^*$ to compute the optimality gap of any set of options.

The complexity of the first loop of Algorithm \ref{alg:bld_alg} is $O \left(L^{K_b} \cdot |\Omega| \cdot T(\ear)\right)$. Being easily vectorized, this loop is much less computationally expensive than the rest of the algorithm. The master problem in the while loop is easy to solve since it is restrained to $\overline{\mathcal{O}}$ and can use the solution of the previous loop iteration as a warm start. Finally, the practical bottleneck of the algorithm is the computation of $V_{T(\ear)}^\text{DFA}\left(S\right)$ whose complexity is $O \left(L^2 \cdot |\Omega| \cdot T(\ear)\right)$. Thus, the worst-case complexity of the algorithm is $O \left(L^{K_b+2} \cdot |\Omega| \cdot T(\ear)\right)$ if there is no limit to the number of options generated i.e. $N_\text{gen} = \infty$. Hence the versatility and performance of Algorithm \ref{alg:bld_alg} comes at the cost of non-negligible computing time as $L$ increases.

We note that the quality of the bounds increases with $\ear$, making $\tilde{\delta}$ a better approximation of $\overline{\delta}$, which in turn improves the quality of the bundles generated. Therefore, although the complexity of each iteration increases linearly in $\ear$, Algorithm \ref{alg:bld_alg} generally requires fewer iterations to produce high-quality bundles as $\ear$ increases. As a result, Algorithm \ref{alg:bld_alg} is able to generate high-quality sets of options in a tractable manner for a wide range of values of $\ear$, as shown in Table \ref{tab:bdl_perf}.





\subsection{Greedy Algorithm}

We now provide a fast heuristic to select a set of options under the MNL choice function. The key idea of this heuristic is to use a value function approximation to estimate the marginal costs in the first period of the season, and select the set $S$ that maximizes the expected revenue at this time period. From Lemma \ref{lemma:inner_problem_solution}, this quantity can be expressed as:
$$R_{T(\ear)}(S) = \frac{-1}{\beta_p} \cdot \mathbb{E}_\rv \left[W\left(\sum_{i \in S} \varphi_i^\rv(\ear)\right)\right]$$
where $\varphi_i^\omega(\ear) = e^{q_i^\omega + \beta_p \Delta_i V^*_{T(\ear)} -1}$. From Theorem \ref{theorem:asymptotic_costs} we have $\varphi_i^\omega(\ear) = O(1/\ear)$. Since $W(0) = 0$ and $W'(0)=1$, the Taylor series of $W$ at $0$ gives $R_{T(\ear)}(S) = \frac{-1}{\beta_p} \cdot \sum_{i \in S} \mathbb{E}_\rv \left[\varphi_i^\rv(\ear)\right] + O(1/\ear^2)$. Thus the revenue at time $T(\ear)$ can be approximately maximized by selecting $S$ maximizing $\sum_{i \in S} \mathbb{E}_\rv \left[\varphi_i^\rv(\ear)\right]$. Let $V$ be an approximation of the value function $V^*$, such as $V^\text{DFA}$ or $V^U$. Ignoring the interaction between options the marginal value of option $i$ becomes 
\begin{equation*}
\Delta_i V_{T(\ear)} \simeq \sum_{l \in i}\Delta_l V_{T(\ear)} = \sum_{l \in i} \left(V_{T(\ear)-1}(S_0)-V_{T(\ear)-1}(S_0\backslash\{l\})\right).
\end{equation*}
Our greedy algorithm works as follows:
\begin{enumerate}
    \item Initialize $S = \emptyset$.
    \item Ranks feasible options in decreasing $\mathbb{E}_\rv \left[e^{q_i^\rv + \beta_p \Delta_i V_{T(\ear)}}\right]$.
    \item Enumerate all options starting from the top-ranked one until $S$ forms a partition of $\mathcal{L}$, and add an option to $S$ if:
    \begin{enumerate}
        \item the option does not overlap with options already included in $S$.
        \item adding the option does not violate the bundle limit $K_s$.
    \end{enumerate}
\end{enumerate}

A pseudocode of this greedy algorithm, Algorithm \ref{alg:greedy_alg}, is provided in the Appendix. Note that to compute every marginal value, one only requires to call the value function approximation $L+1$ times. In contrast, the column generation algorithm may require one call per option, resulting in up to $L^{K_b}$ calls. This makes the greedy algorithm much faster as the number of items or maximum bundle size increases.\hfill\\
\section{Computational Results and Industry-Driven Case Study} \label{section:ExperimentalResults}


To implement the algorithms presented in Section \ref{section:Bunlding}, we partnered with the third-party logistics provider (3PL) Uber Freight, a growing segment within Uber Technologies that provides a digital platform for shippers and carriers to facilitate the movement of freight. This collaboration allowed us access to real-world data and a carrier choice model. In 2023, Uber Freight moved 26 billion pounds of freight across 9,500 cities in the United States, for a total of 306 million miles driven by carriers on its platform. Already leveraging dynamic pricing, freight brokerage platforms managing long-haul full truckloads have a significant opportunity to improve load recommendations and encourage carriers to smartly string loads together through bundle suggestions.

\begin{figure}[h]
    \centering
    \begin{subfigure}{0.65\linewidth}
        \centering
        \includegraphics[width=\linewidth]{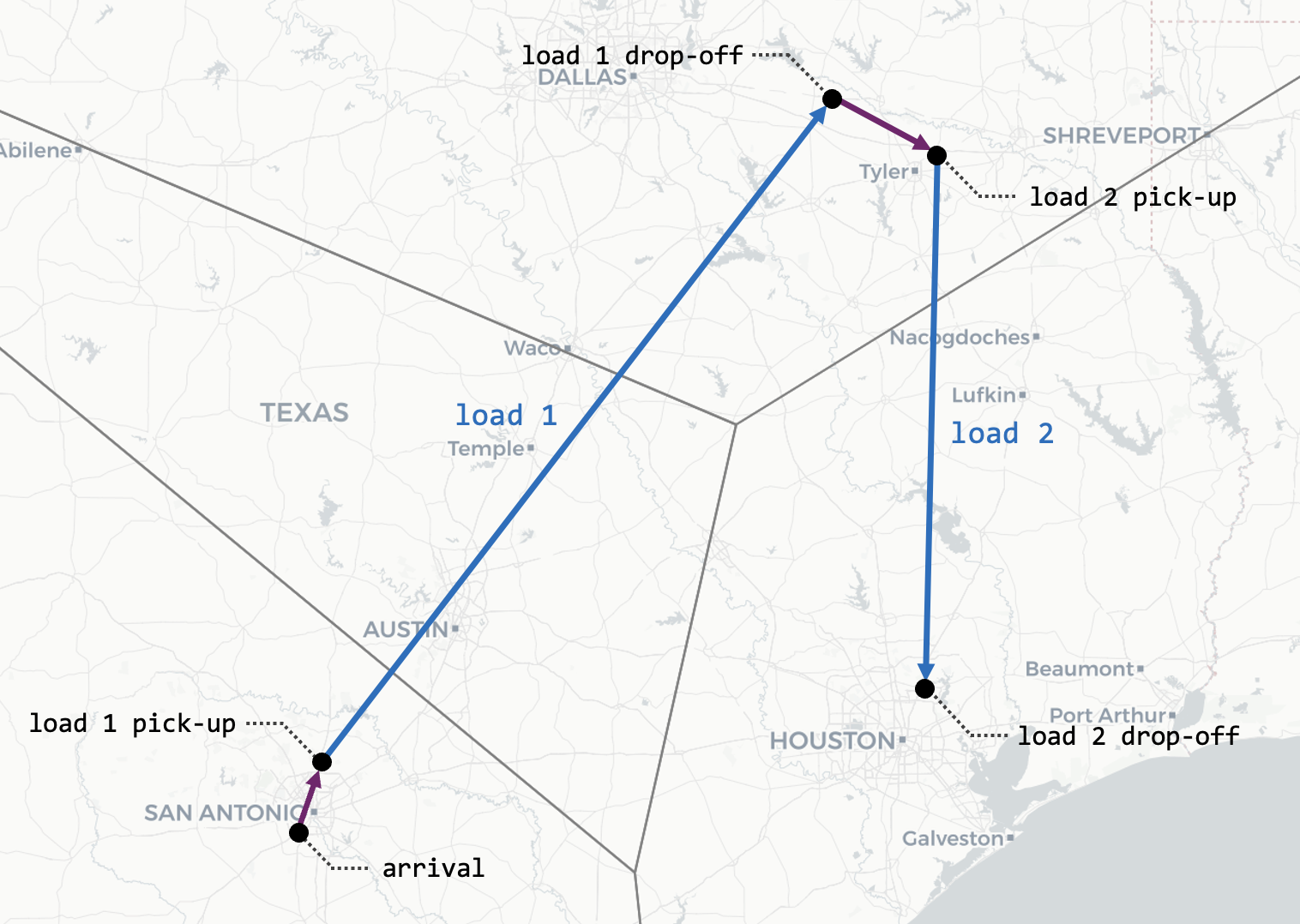}
    \end{subfigure}
    \caption{Example of a bundle delivery. In-load and empty miles are denoted in blue and purple respectively.}
    \label{fig:ex_texas_triangle}
\end{figure}

On one side, shippers provide \textit{loads} (items) consisting of a single shipment associated with an origin-destination pair and an expiration time. On the other side, \textit{carriers} (customers) arrive on the platform, observe the option-price pairs displayed by the platform, and decide to book one of the options or leave the platform. By bundle, we refer here to two loads to be picked up and dropped off sequentially as illustrated in figure \ref{fig:ex_texas_triangle}, and offered together by the 3PL provider as a single option. Specifically, carriers cannot transport multiple loads simultaneously. Note that the order of loads in a bundle can drastically affect empty miles by influencing both the distance to first pickup and the distance between loads. The objective of the platform is to minimize its operational costs, or equivalently, maximize its revenue, by selecting and pricing options to offer to arriving carriers. In this case study, we analyze not only revenue, but also empty miles and timely load delivery as markers of efficient platform operations.

The data includes 959 inter-region loads in the Texas Triangle, with specified pick-up and drop-off location pairs, supplied to the platform over 2.5 months. We partition the space into 4 regions (San Antonio, Austin, Dallas, Houston) using K-means clustering on the drop-off location of the loads and let $\Omega = \{\omega_\text{S},\omega_\text{A},\omega_\text{D},\omega_\text{H}\}$ be the set of regions with corresponding centroids $\bar{\omega}_\text{S}, \bar{\omega}_\text{A}, \bar{\omega}_\text{D}, \bar{\omega}_\text{H} \in \R^2$.
We have access to a logit model that predicts the probability of a carrier accepting an option, with parameters calibrated using historical data from our industry partner. Note that this model does not suffer from the \say{spiral-down effect} described in \cite{Spirals} as it is fitted on historical data that is independent of it. We use here a MNL model with these parameters. In the Appendix [\ref{apx:case_study}], we show that simulating carrier choices as a sequence of binary decisions yields similar results and insights. Further, the simulation without bundles yields a cost-per-mile distribution similar to the historical one, indicating that the MNL model accurately captures carrier preferences within the scope of this simulation. The model includes several parameters: $\beta_0$, the MNL constant; $\beta_p$, $\beta_d$, $\beta_{e}$, and $\beta_{b}$, the sensitivity to price, in-load miles, empty miles, and bundle respectively; $\beta_{\texttt{org}}(\omega)$ and $\beta_{\texttt{dst}}(\omega)$ for all $\omega \in \Omega$, denoting origin and destination region preferences. The perceived quality of option $i$ for a carrier arriving in region $\omega \in \Omega$ is then computed as follows:
$$q_i^\omega = \beta_0 + \beta_d \cdot d_i + \beta_{e} \cdot e_i^\omega + \beta_{b} \cdot \mathds{1}_{\{\#i=2\}} + \beta_{\texttt{org}}(\omega^\texttt{org}_i) + \beta_{\texttt{dst}}(\omega^\texttt{dst}_i),$$
where $d_i$ is the total in-load distance of option $i$, $e_i^\omega$ is the total empty miles from delivering option $i$ starting from location $\bar{\omega}$, and $\omega^\texttt{org}_i$ and $\omega^\texttt{dst}_i$ are respectively the origin and destination regions of option $i$ i.e. regions whose centroids are closest to the pick-up and drop-off location of option $i$. Here, a bundle is complementary if it results in low empty miles and the markets of first load pickup and last load drop-off are attractive, and substitutable otherwise.

Salvage values are difficult to estimate, even for 3PLs, since the freight industry is volatile and rates can vary significantly within days, most loads are booked well ahead of their expiration time, and observations of these salvage values are subject to survival bias since underpriced loads are more likely to reach their expiration time. If a load is not booked by its expiration time, we assume that a proprietary truck will deliver it from the centroid of the load's origin region at a cost per mile resulting in a proportion of loads delivered on time in the dynamic setting without bundles that match the historical proportion.

We point out that, unlike in the general retail setting, here price sensitivity of carriers is positive. This makes the price trajectories non-decreasing instead of decreasing, but the results of Section \ref{section:DP} and Section \ref{section:bounds} still hold in this context. Further, we only allow bundles of up to two loads. As the sensitivity to larger bundles is very negative, allowing bundles of 3 or more loads does not significantly affect the maximum revenue improvement here (unless $\ear \ll L$) while substantially increasing the computational cost in the dynamic setting.


\subsection{Static setting}

In the static setting, we consider a finite set of loads $\mathcal{L}$ that all have the same expiration time which marks the end of the season. Carriers only arrive at the centroid of each region, with the distribution of arrivals across regions following the distribution of historical loads drop-off region. A set of options is formed at the beginning of the season and cannot be changed subsequently. Consequently, the platform has to solve a problem of the form \eqref{eq:first_stage}.

We first compare the performance of the bounds developed in Section \ref{section:bounds} with the fluid and static approximations for varying numbers of loads generated $L$ and demand\footnotemark{}\footnotetext{In the context of 3PLs, demand usually refers to freight, but we choose to refer to carriers here to remain consistent with previous sections.} levels $\ear$. All loads are offered individually, and $L$ is small enough that we can calculate the exact expected revenue. Table \ref{tab:bounds_perf} reports the relative error in expected revenue $(V_{T(\ear)}-V_{T(\ear)}^*)/V_{T(\ear)}^*$ for different value approximations $V$. We report the relative error, averaged over 50 runs with loads randomly sampled from historical loads, of the fluid approximation $V^\text{Fluid}$ [\ref{apx:fluid_approx}], the backward upper bound $V^U$ [\ref{subsection:backward}], the Deterministic Forward Approximation $V^\text{DFA}$ [\ref{subsection:DFA}], the static approximation $V^\text{Static}$ [\ref{apx:static_approx}], and the backward lower bound $V^L$ [\ref{apx:backward_lower_bound}]. We first note that $V^\text{Static}$ and $V^L$ perform similarly, with $V^L$ gaining the edge as the demand per load increases. Between two upper bounds, $V^\text{Fluid}$ and $V^U$, we observe that $V^\text{Fluid}$ performs better for $\ear/L$ small, but the relative error of the fluid approximation increases much more than that of $V^U$ as $\ear/L$ increases. Finally, the DFA lower bound outperforms all other approximations. In particular, the gap in performance between the DFA and the fluid and static approximations widens significantly as $\ear$ increases.

\begin{table}[H]
\centering
{
\fontsize{8}{12}\selectfont
\begin{tabular}{|l|ccc|ccc|ccc|}
\hline
$\ear$ & \multicolumn{3}{c|}{5} & \multicolumn{3}{c|}{20} & \multicolumn{3}{c|}{50} \\
\hline
$L$ & 3 & 5 & 8 & 3 & 5 & 8 & 3 & 5 & 8\\
\hline
$V^\text{Fluid}$     &   0.0540 &   0.0365 &   0.0247 &   0.1612 &   0.1536 &   0.1443 &   0.2114 &   0.2078 &   0.2036 \\
$V^U$     &   0.0267 &   0.0428 &   0.0572 &   0.0185 &   0.0317 &   0.0456 &   0.0109 &   0.0188 &   0.0274 \\
$V^\text{DFA}$ &  \textbf{-0.0068} &  \textbf{-0.0070} &  \textbf{-0.0073} &  \textbf{-0.0052} &  \textbf{-0.0071} &  \textbf{-0.0096} &  \textbf{-0.0023} &  \textbf{-0.0032} &  \textbf{-0.0046} \\
$V^L$     &  -0.0104 &  -0.0121 &  -0.0121 &  -0.0212 &  -0.0322 &  -0.0404 &  -0.0200 &  -0.0326 &  -0.0442 \\
$V^\text{Static}$   &  -0.0115 &  -0.0126 &  -0.0123 &  -0.0336 &  -0.0410 &  -0.0467 &  -0.0526 &  -0.0603 &  -0.0678 \\
\hline
\end{tabular}
}
\caption{Relative error in expected revenue.}
\label{tab:bounds_perf}
\end{table}

Since the relative error of the DFA does not exceed $1\%$ in our experiments, and since the computation time of the exact expected revenue explodes with $L$, we use the DFA to evaluate the revenue improvement of different bundling algorithms compared to the \say{no bundle} set of options. Table \ref{tab:bdl_perf} shows the expected revenue improvement for the set of options returned by Algorithm \ref{alg:bld_alg} (column generation), the sets of options returned by Algorithm \ref{alg:greedy_alg} (greedy algorithm) with value function approximations $V^\text{Static}$, $V^\text{Fluid}$, and $V^\text{DFA}$, and the set of options that minimizes expected empty miles following load drop-offs. We use the latter set of options as a benchmark given 3PLs' efforts to reduce driving distances as a means of jointly increasing revenue and reducing emissions. It is obtained by solving a MILP using historical deadhead miles and is provided in the Appendix [\ref{subsection:min_empty_miles}]. The results are averaged over 50 runs with loads randomly sampled from the historical loads. Although computationally expensive for large values of $L$, the column generation algorithm outperforms all other approaches, yielding the best revenue for all $\lambda$ and $L$ considered. The greedy algorithm performs similarly for all three value-function approximations for small $\lambda$. However, as $\lambda$ increases, the DFA enables the greedy algorithm to perform comparably to the column generation algorithm, while the static and fluid approximations result in negligible or even negative revenue improvement. Intuitively, the fluid approximation overestimates expected revenue by an amount proportional to the number of options when demand is high (as shown in Proposition \ref{prop:fluid_approx_asymptotical}), thus overvaluing single item options. Conversely, the static approximation undervalues each option, resulting in overly aggressive bundling. Finally, while the set of options minimizing empty miles outperforms the set of options without bundles for $\lambda$ small, it fails to capture carrier preferences and does not account for the penalizing effects of a reduced number of options when demand is high, yielding significantly lower revenues compared to other bundling methods and a substantially negative revenue improvement for $\lambda$ large. Thus, the column generation algorithm appears to be the best algorithm, while the greedy algorithm with DFA value function approximation emerges as the best alternative if the number of loads makes column generation too expensive for the application considered.

\begin{table}[H]
\centering
{
\fontsize{8}{12}\selectfont
\begin{tabular}{|l|ccc|ccc|ccc|}
\hline
Demand $(\ear)$ & \multicolumn{3}{c|}{10} & \multicolumn{3}{c|}{30} & \multicolumn{3}{c|}{50} \\
\hline
Number of loads $(L)$  &       5 &      10 &      20 &      5 &      10 &      20 &  5 &      10 &      20 \\
\hline
Column Generation   &  +5.005\% &  +7.974\% &  +9.655\% &  +1.449\% &  +3.401\% &  +7.621\% &  +0.000\% &  +0.467\% &  +1.928\% \\
Greedy $<$DFA$>$    &  +4.900\% &  +7.550\% &  +9.220\% &  +1.386\% &  +3.052\% &  +7.017\% &  +0.000\% &  +0.337\% &  +1.543\% \\
Greedy $<$Fluid$>$  &  +4.641\% &  +7.371\% &  +9.218\% &  +0.019\% &  +1.229\% &  +7.006\% &  +0.000\% &  +0.000\% &  +0.221\% \\
Greedy $<$Static$>$ &  +4.817\% &  +7.531\% &  +9.180\% &  +1.249\% &  +2.923\% &  +7.018\% &  -0.272\% &  -0.793\% &  +0.375\% \\
Min. Empty Miles    &  +3.448\% &  +5.993\% &  +7.349\% &  -0.449\% &  +1.621\% &  +5.065\% &  -5.351\% &  -5.272\% &  -2.942\% \\
\hline
\end{tabular}
}
\caption{Expected revenue improvement relative to the “no bundle” set of options.}
\label{tab:bdl_perf}
\end{table}

\subsection{Dynamic setting}

The dynamic setting provides a more realistic environment where loads have different expiration times $\bar{t}_l$, are successively supplied and booked by carriers, and bundles are updated regularly. Moreover, carriers' arrival locations are more fine-grained in that they follow the exact distribution of historical load drop-off locations. This is as accurate as we can model arrivals based on our data, since carriers tend to log into the platform after completing a delivery. 

We incorporate two additional considerations to enhance the simulation's accuracy in reflecting real-world conditions. First, carriers only observe a subset of the options offered by the platform. When a carrier arrives at location $\omega$, options in $S$ are sorted in decreasing $q_i^\omega + \beta_p \Delta_i V_t$ such that the expected revenue is maximized under optimal prices according to Lemma \ref{lemma:inner_problem_solution}. The carrier then observes the first $k$ options of $S$, with $k$ following a distribution provided by our industry partner. Second, we require prices to be homogeneous (independent of carriers' location) notably to increase carrier satisfaction with fair pricing and prevent carriers from abusing the pricing system by manipulating their virtual position.

While real operations can benefit from prices obtained using the backward upper bound and marginal costs estimated using DFA (see Section \ref{section:bounds}), here we consider a closed-form approximation derived next to speed up the simulation. This allows us to simulate trucking operations over a 2-year horizon in less than 3 hours, with statistically significant results.

The pricing policy we use for individual loads is derived from the optimal and asymptotic prices given in Lemma \ref{lemma:inner_problem_solution} and Theorem \ref{theorem:asymptotic_costs}. The \textit{expiration price} $\bar{p}_l$ of load $l$ (price at time $\bar{t}_l-1$, the period preceding the expiration of $l$) is set to be the expected optimal price across regions assuming that $l$ is the only available option at that time:
\begin{equation*}
    \bar{p}_l \triangleq \xi_l - \frac{1}{\beta_p}\left(1+\mathbb{E}_\rv\left[W\left(e^{q_l^\rv+\beta_p \cdot \xi_l - 1}\right)\right]\right).
\end{equation*}

Moreover, we know from Theorem \ref{theorem:asymptotic_costs} that prices are asymptotically logarithmic in demand, which lets us consider the following pricing policy for individual loads:
\begin{equation*}
p_{t,l} \triangleq \frac{-1}{\beta_p}\left(\ln\left(e^{-(\beta_p \cdot \bar{p}_l + \kappa_l)} + \alpha \cdot \arate \cdot (\bar{t}_l - 1 - t) \right)+ \kappa_l\right), \quad t = \{0,\hdots,\bar{t}_l-1\},
\end{equation*}
where $\alpha$ is a parameter controlling the curvature and aggressiveness of the price trajectory. This parameter is determined empirically to maximize revenue in the \say{no bundle} framework. Note that the larger $\beta_p \cdot \bar{p}_l + \kappa_l$, the steeper the curvature of the price trajectory near the expiration time. As such, a discrepancy between expiration prices and asymptotic quality results in more volatile prices, which is a desirable behavior observed on optimal price trajectories of \eqref{eq:bellman_equation}.

We now turn our attention to bundle pricing. Our partner offers personalized bundles with a linear pricing strategy: the price of any bundle equates to the sum of the prices of the loads it contains. Namely, the \textit{linear price} of a bundle $b$ at time $t$ is $p_{t,b} \triangleq \sum_{l \in b} p_{t,l}$. This pricing policy is popular in practice due to its simplicity and because it reduces the joint bundling and pricing problem to an assortment problem if the prices of loads are fixed. Our analysis in Section \ref{section:DP} suggests that the platform can benefit from offering custom (bundle-specific) prices that depend on the quality and marginal value of bundles instead of using linear prices. We estimate the marginal value of load $l$ by assuming that the price trajectory $(p_{t,l})_{0\leq t \leq \bar{t}_l-1}$ is optimal, in which case the relation between optimal prices and marginal value given in Lemma \ref{lemma:inner_problem_solution} gives:
\begin{equation*}
\Delta_l V_t = p_{t,l} + \frac{1}{\beta_p}\left(1+\mathbb{E}_\rv\left[e^{q_l^\rv+\beta_p \cdot p_{t,l}}\right]\right).\end{equation*}
Ignoring the interactions between loads, the marginal value of a bundle $b$ becomes $\Delta_b V_t = \sum_{l \in b} \Delta_l V_t$.
We define the \textit{custom price} of 
$b$ at time $t$ as $p_{t,b}^\dagger \triangleq \Delta_b V_t -\frac{1}{\beta_p}\left(1+\mathbb{E}_\rv \left[ W\left(e^{q_b^\rv+\beta_p \cdot \Delta_b V_t - 1}\right]\right)\right)$, the optimal price given in Lemma \ref{lemma:inner_problem_solution} function of the quality of bundle $b$ and its marginal value. Crucially, custom pricing takes into account the substitutability or complementarity of bundles, whereas linear pricing is agnostic to the quality of bundles and only accounts for the quality of the loads composing them.

In parallel to the pricing policies, we compare three frameworks: no bundling, \textit{rolling horizon} bundling (the set of options $S$ is recomputed every $T$ time steps and supplied loads are offered individually until options are recomputed), and \textit{personalized} bundling (the set of options $S$ is recomputed for each arrival based on the exact carrier location).

We simulate 2 years of operations based on the historical load distribution. The main results of the simulation are shown in Table $\ref{tab:simulation_results}$. The costs per loaded mile, the average empty miles per load, and the proportion of loads that failed to be delivered on time are reported for every bundling and pricing method with an associated 99\% confidence interval.

We first note that taking into account carriers' preferences and expected demand (arrivals) is crucial for improving the platform's revenue. Indeed, selecting bundles based solely on expected empty miles following load drop-offs fails to form attractive bundles or bundles that align with the spatial distribution of carriers. In contrast, the greedy algorithm is able to generate bundles resulting in improved revenue, reduced empty miles, and reduced proportion of unmet deadlines compared to the \say{no bundle} strategy. Furthermore, accounting for the distribution of carrier arrivals allows the greedy algorithm to achieve lower empty miles than the sets of options that minimize expected empty miles following load drop-offs. We also observe, especially with linear pricing, that poorly designed or mispriced bundles can lead to a cannibalism effect, whereby unattractive options presented to carriers reduce the number of attractive options observed by carriers. This cannibalism effect leads to an accumulation of loads on the platform as shown in Figure \ref{fig:acceptance_distrib}, resulting in both higher costs and fewer loads delivered on time. Here, linear pricing is a driver of cannibalization as it tends to underprice bundles. Among underpriced bundles, only high quality ones (with low empty miles while meeting other carrier preferences) are likely to be accepted, which slightly reduces empty miles. On the other hand, custom prices lead to early acceptances of options, increased number of loads delivered on time, and reduced costs by matching options to a carrier before their prices spike. Early acceptances also prevent loads to accumulate on the platform which improves the proportion of options observed in future arrivals and further decreases the likelihood of cannibalism effect.

Finally, we find that while personalized bundles perform better than bundles recomputed on a rolling basis with our greedy algorithm, the latter framework is able to capture most of the improvement across all three reported performance metrics despite being less computationally intensive since far fewer options need to be priced and ranked for each arrival. This can be explained by the fact that, despite the large number of feasible options, only a few high-quality options need to be offered in each region to achieve close to optimal revenue, which our greedy algorithm successfully identifies by taking into account the distribution of arrival locations.

\begin{table}
\centering
{
\small
\begin{tabular}{|l|c|c||c|c|c|}
\hline
Framework & \makecell{Bunlding \\ method} & Pricing & \makecell{Cost/loaded mile \\ (\$/mi.)}& \makecell{Avg. empty miles \\ (mi.)} & Unmet deadlines\\
\hline
No bundle & - & - & $2.8181\pm0.0135$ & $33.1618\pm0.6271$ & $5.81\%\pm0.41$\\
\hline
\multirow{4}{*}{Rolling horizon}
& \multirow{2}{*}{Empty miles} 
& Linear & $3.0151\pm0.0152$ & $28.7050\pm0.5202$ & $10.82\%\pm0.54$\\
    \cline{3-6}
& & Custom & $2.8668\pm0.0120$ & $28.9775\pm0.5098$ & $3.66\%\pm0.33$\\
\cline{2-6}
& \multirow{2}{*}{Greedy Alg.} 
& Linear & $2.7018\pm0.0111$ & $25.8897\pm0.4105$ & $4.24\%\pm0.35$\\
    \cline{3-6}
& & Custom & $2.6860\pm0.0121$ & $26.1519\pm0.4268$ & $2.49\%\pm0.27$\\
\hline
\multirow{2}{*}{Personalized}
& \multirow{2}{*}{-}
& Linear & $2.6904\pm0.0093$ & $\bm{23.6795\pm0.3264}$ & $4.75\%\pm0.37$\\
    \cline{3-6}
& & Custom & $\bm{2.6485\pm0.0107}$ & $24.6380\pm0.3508$ & $\bm{2.27\%\pm0.26}$\\
\hline
\end{tabular}
}
\caption{Two-year simulation: operational results.}
\label{tab:simulation_results}
\end{table}

Overall, we find that joint bundling and pricing can improve revenue over dynamic pricing alone, provided that the bundles are carefully chosen. The performance of the rolling horizon framework shows that our bundling approach for the static setting can be successfully extrapolated to a dynamic setting through periodic re-optimization. Compared to the \say{no bundle} strategy, personalized bundling with custom pricing reduces costs by 6\%, empty miles by 25\%, and more than halves the number of unmatched loads.\hfill\\
\section{Concluding Remarks} \label{section:conclusion}

In this article, we studied the bundling and dynamic pricing of nonreplenishable and unique items. Our methodology (i) provides insights on the optimal bundling structure by obtaining near closed-form solutions under MNL choice function to both the single-period pricing problem and the Unique Item Pricing problem as the demand tends to infinity; (2) enables computation of tractable bounds on the expected revenue over a finite horizon; and (3) leverages these bound to determine what bundles, if any, to offer.

Notably, the bounds and bundling approach we provide are amenable to a variety of choice functions, enabling bundling under a variety of choice models, including latent class and nested logit models. Under the MNL model, we prove that these bounds are asymptotically optimal as the horizon grows, and yield a performance bound in $O\left(1/\ear\right)$ as the demand $\ear$ increases. We proposed a column generation-based bundling algorithm that exploits these bounds to efficiently form bundles that maximize our lower bound on expected revenue. We then presented a greedy algorithm under the MNL model that enables bundling large numbers of items in a dynamic bundling setting. Finally, we conducted a numerical study to assess the performance of our bounds against the static and fluid approximations, the efficiency of our bundling algorithms, and the potential of bundling to improve expected revenue when combined with dynamic pricing. Our industry-driven case study on third-party logistics shows that carefully selected bundles can increase the platform's revenue while reducing empty miles and the proportion of unmatched loads.

Our model could be further refined to account for non-stationary arrivals, non-stationary item quality, or item-dependent expiration time. Despite the increased complexity of the theoretical analysis introduced by these relaxations, we remain convinced that the approach outlined in this paper would be suitable for these variations of the problem we have covered. Moreover, the DFA could be extended to non-unique items: instead of recursively computing a scalar denoting the availability of an item at a given time period, one could propagate a vector whose $k$-th entry denotes the probability that $k$ units of that item are in stock at that time period. Finally, in the context of freight operations, it would be interesting to explore the potential of bundling to achieve desirable long-term carrier positioning and further reduce costs and empty miles.

We are hopeful that our framework provide an appropriate model for dynamic pricing problems in practice. In many cases, dynamic pricing problems can be studied in a large inventory setting by grouping similar items together. Although this simplification makes the fluid approximation well-suited, it comes at the expense of treating potentially dissimilar items as identical and interchangeable. Our framework allows for a more fine-grained approach, thereby enabling more accurate solutions and explicit bundles.\hfill\\

\bibliographystyle{informs2014}
\bibliography{refs}

\newpage

\newcounter{storedlemma}
\setcounter{storedlemma}{\value{lemma}}
\renewcommand{\thelemma}{A-\the\numexpr\value{lemma}-\value{storedlemma}\relax}

\newcounter{storedequation}
\setcounter{storedequation}{\value{equation}}
\renewcommand{\theequation}{A.\the\numexpr\value{equation}-\value{storedequation}\relax}

\newcounter{storedproposition}
\setcounter{storedproposition}{\value{proposition}}
\renewcommand{\theproposition}{A-\the\numexpr\value{proposition}-\value{storedproposition}\relax}

\newcounter{storedfigure}
\setcounter{storedfigure}{\value{figure}}
\renewcommand{\thefigure}{A-\the\numexpr\value{figure}-\value{storedfigure}\relax}

\newcounter{storedtable}
\setcounter{storedtable}{\value{table}}
\renewcommand{\thetable}{A-\the\numexpr\value{table}-\value{storedtable}\relax}

\section*{Appendix}\hfill

\renewcommand{\thesubsection}{A.\arabic{subsection}}

\subsection{More on the Introduction's Generic Example} \label{apx:intro_example}

The optimal pricing policies of the introductory example are obtained by solving the dynamic program \eqref{eq:bellman_equation}. The dynamic program can be solved to optimality here, since no more than 3 options are offered in each of the sets considered. As the optimal price of an option depends not only on demand but also on the availability of other options, we plot the optimal price of each option when offered alongside all other options or as the only remaining option.

\begin{figure}[htb]
    \centering
    \begin{subfigure}[b]{0.49\textwidth}
        \centering
        \includegraphics[width=\textwidth]{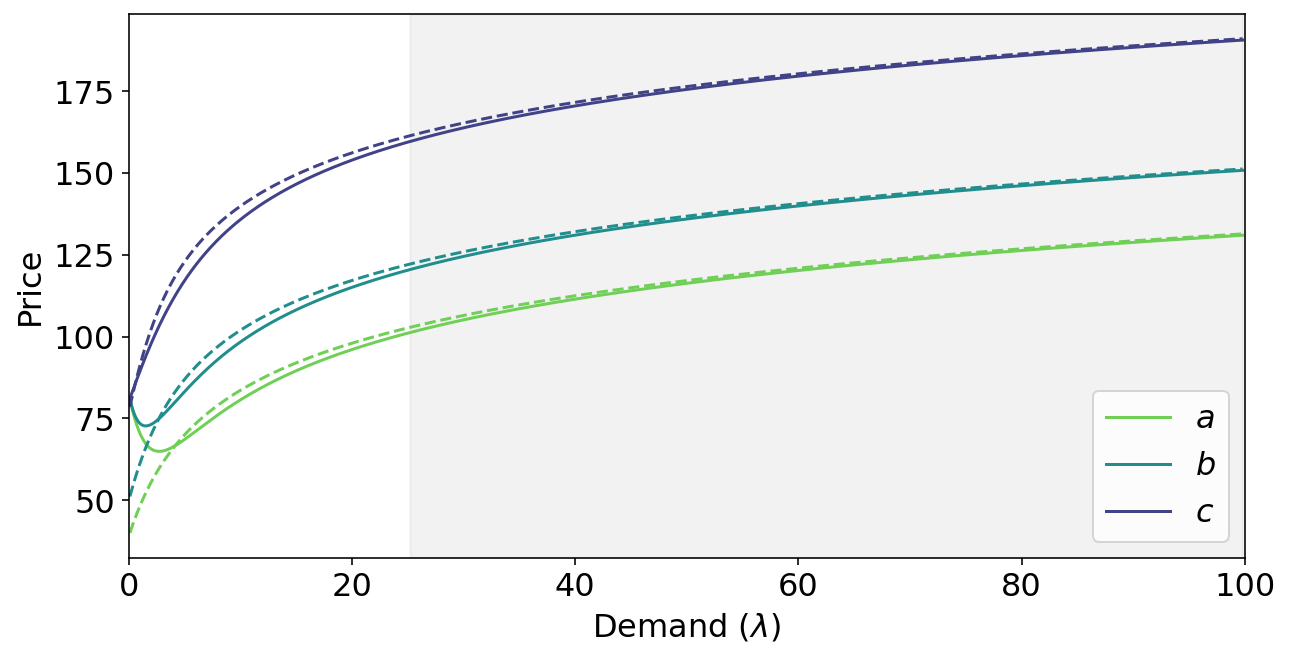}
        \caption{$S_0$}
    \end{subfigure}
    \hfill
    \begin{subfigure}[b]{0.49\textwidth}
        \centering
        \includegraphics[width=\textwidth]{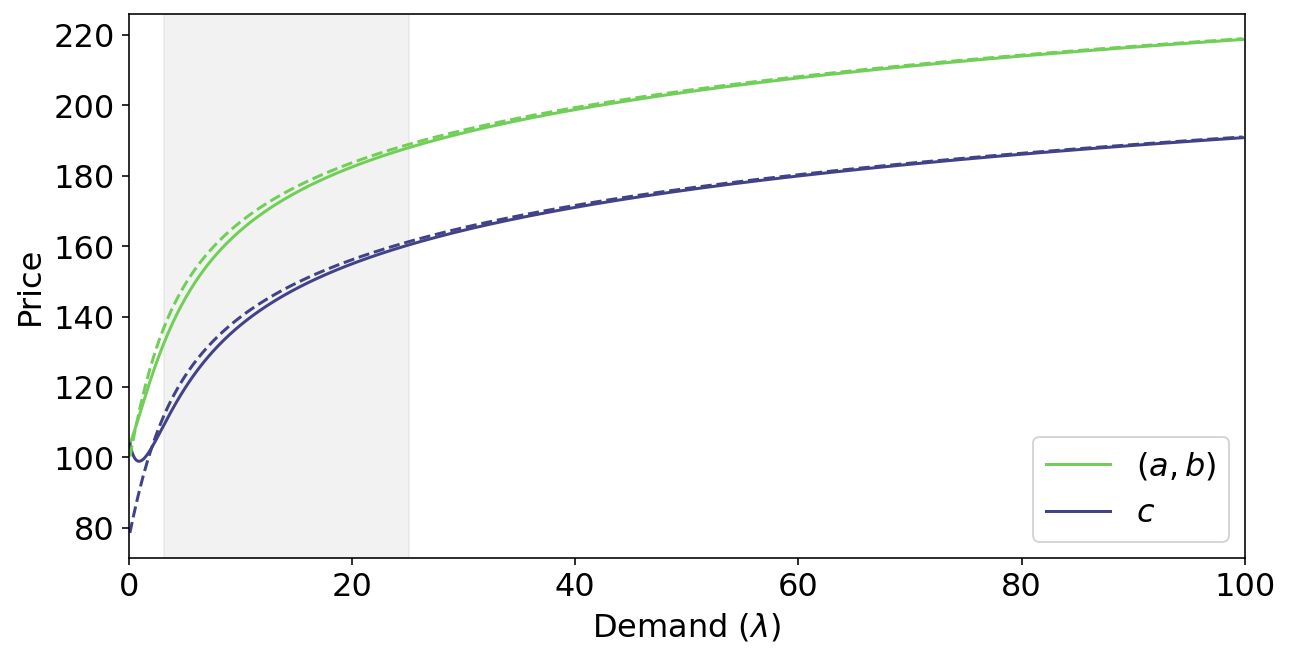}
        \caption{$S_1$}
    \end{subfigure}
    \begin{subfigure}[b]{0.49\textwidth}
        \centering
        \includegraphics[width=\textwidth]{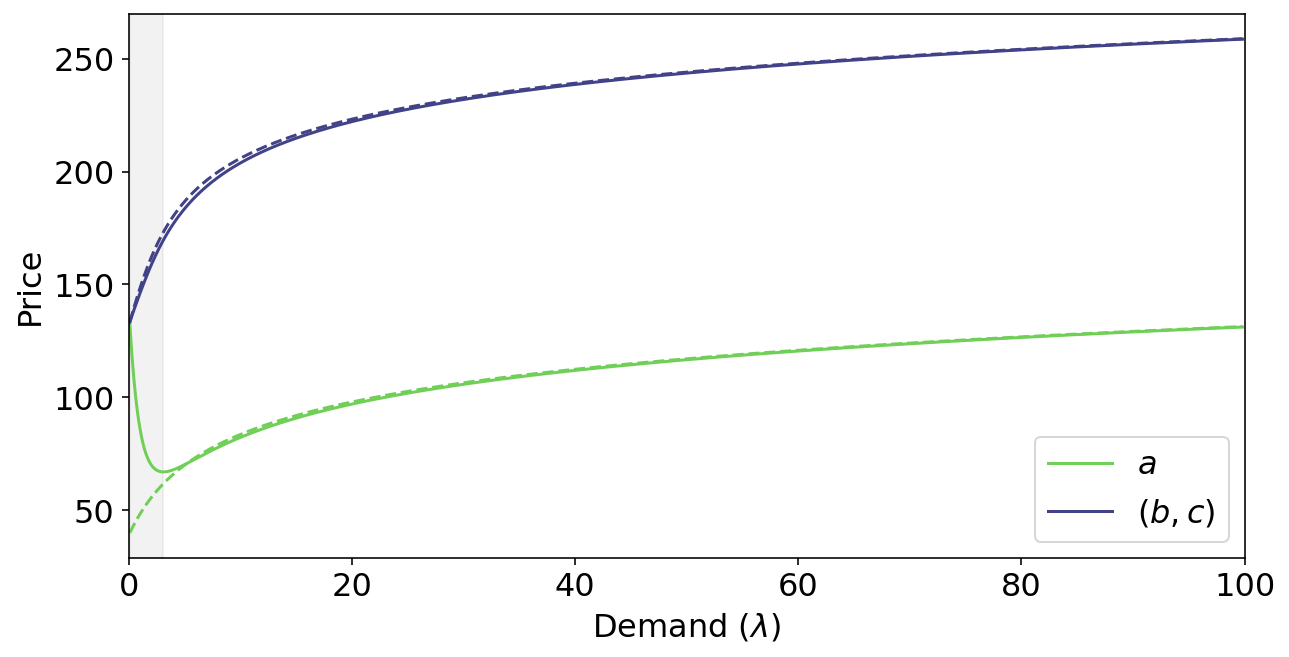} 
        \caption{$S_2$}
    \end{subfigure}
    \caption{Optimal prices of each option function of demand for the 3 sets of options considered. Solid lines denote optimal prices when all options are available and offered. Dashed lines denote optimal prices when an option is the only one offered, i.e. when other options have been purchased and are no longer available. The gray area highlights the demand interval in which that particular set of options is optimal.}
    \label{fig:prices_intro_example}
\end{figure}

\subsection{Fluid Approximation} \label{apx:fluid_approx}
A common approach to simplify stochastic dynamic programs consists of formulating a deterministic problem where inventory levels are continuous, options are sold in continuous amounts, and the demand is deterministic. With initial inventory levels equal to one, the fluid approximation of problem \eqref{eq:bellman_equation} is:
\begin{align*}
    \max_{\rho(\cdot)} \quad & \sum_{i \in \bar{S}}\int_0^T \arate \E_\rv \left[\rho_i^\rv(t) p_i(\rho^\rv(t))\right] dt + \sum_{i \in \bar{S}} \xi_i \left(1-\int_0^T \arate \E_\rv \left[\rho_i^\rv(t)\right] dt\right)\\
    \text{s.t.}\quad & \int_0^T \arate \E_\rv \left[\rho_i^\rv(t)\right] dt \leq 1, \quad \forall i \in \bar{S}\\
    &\sum_{i \in \bar{S}} \rho_i^\omega(t) \leq 1, \quad \forall \omega \in \Omega, \; \forall t\in[0,T]
\end{align*}
where the first constraint says that each item can be sold at most once, and the second constraint ensures that the probability of an item being purchased is less than one for each customer type.
As explained in \cite{Multiproduct_Dynamic_Pricing_Problem_Yield}, this problem reduces to a concave problem with constant solutions when the revenue and demand function are time-invariant. The problem can be reformulated as:
\begin{align}
    \max_{\rho \in (0,1)^{N\times|\Omega|}} \quad & \arate T \sum_{i \in \bar{S}} \E_\rv \left[\rho_i^\rv \cdot \left(p_i^\rv(\rho^\rv) - \xi_i\right)\right] + \sum_{i \in \bar{S}} \xi_i \tag{$\mathcal{P}_\text{Fluid}$} \label{eq:fluid_approx}\\
    \text{s.t.} \quad & \E_\rv \left[\rho_i^\rv\right] \leq \frac{1}{\arate T}, \quad \forall i \in \bar{S} \nonumber\\
    & \sum_{i \in \bar{S}} \rho_i^\omega \leq 1, \quad \forall \omega \in \Omega \nonumber
\end{align}
We denote by $V_T^\text{Fluid}\left(\bar{S}\right)$ the optimal solution to problem \eqref{eq:fluid_approx}. Fluid approximations are often asymptotically optimal upper bound to the original stochastic problem as a common factor scales the initial inventory levels and number of periods. However, given a fixed inventory level, the fluid approximation of the UPP problem is not asymptotically optimal as the number of periods tends to infinity.

\begin{proposition}\label{prop:fluid_approx_bound}
The optimal objective of \eqref{eq:fluid_approx} upper bounds the optimal objective of the original stochastic problem. Namely, $V_T^\text{Fluid}\left(\bar{S}\right) \geq V_T^*\left(\bar{S}\right)$ for any $T \in \N$.
\end{proposition}

\begin{cproof}{Proposition}{\ref{prop:fluid_approx_bound}}
Let $V_T^\text{Fluid}\left(\bar{S},d\right)$ be the expected revenue of the perfect foresight model with $d^\omega$ arrivals of type $\omega$:
\begin{align}
    V_T^\text{Fluid}\left(\bar{S}, d\right) = \max_{y \in \R_+^{N\times|\Omega|}} \quad & \sum_{i \in \bar{S}} \sum_{\omega \in \Omega} \left[y_i^\omega \cdot \left(p_i^\omega(y^\omega/d^\omega) - \xi_i\right)\right] + \sum_{i \in \bar{S}} \xi_i\\
    \text{s.t.} \quad & \sum_{\omega \in \Omega} y^\omega_i \leq 1, \quad \forall i \in \bar{S} \nonumber\\
    & \sum_{i \in \bar{S}} y_i^\omega \leq d^\omega, \quad \forall \omega \in \Omega \nonumber
\end{align}
For every instance of $d$, the expected sales under the optimal dynamic policy obtained solving \eqref{eq:bellman_equation} constitute a feasible solution to the perfect foresight model. Since stationary solutions are optimal under known demand, the corresponding objective value is larger than the one achieved by the optimal dynamic policy in the UPP problem. Thus the revenue obtained by the optimal policy under demand realization $d$ is at most $V_T^\text{Fluid}\left(\bar{S},d\right)$. Taking the expectation we obtain $V_T^*\left(\bar{S}\right) \leq \E_d\left[V_T^\text{Fluid}\left(\bar{S},d\right)\right]$. Note that $V_T^\text{Fluid}\left(\bar{S},d\right)$ is concave in $d$. Indeed, a convex combination of feasible solutions for two values of $d$ is feasible for the convex combination of the values of $d$, so the optimal solution is at least as large as the convex combination of the objectives. Therefore, by Jensen's inequality, we have $\E_d\left[V_T^\text{Fluid}\left(\bar{S},d\right)\right] \leq V_T^\text{Fluid}\left(\bar{S}, \E_d\left[d\right]\right) = V_T^\text{Fluid}\left(\bar{S}\right)$. Combining the two inequalities, we obtain $V_T^\text{Fluid}\left(\bar{S}\right) \geq V_T^*\left(\bar{S}\right)$.
\end{cproof}

\begin{proposition}\label{prop:fluid_approx_asymptotical}
The optimal objective of \eqref{eq:fluid_approx} does not converge to the optimal objective of the stochastic problem as $T$ tends to $\infty$. Moreover, $\lim_{T \to \infty}  V_T^\text{Fluid}(\bar{S}) - V_T^*(\bar{S}) \geq -\frac{|\bar{S}|}{\beta_p} > 0$.
\end{proposition}

\begin{cproof}{Proposition}{\ref{prop:fluid_approx_asymptotical}}
Let $T$ be such that $T \geq \frac{1}{\arate} \max_\omega \sum_{i \in \bar{S}} e^{q_i^\omega - \kappa_i}$. Consider the solution $\hat{\rho}$ defined by $\hat{\rho}_i^\omega = \frac{1}{\arate T} e^{q_i^\omega - \kappa_i}$. Then for any $i \in \bar{S}$ we have $\E_\rv \left[\hat{\rho}_i^\rv\right] = 1/\arate T$ and $\sum_{i \in \bar{S}} \hat{\rho}_i^\omega = \frac{1}{\arate T} \sum_{i \in \bar{S}} e^{q_i^\omega - \kappa_i} \leq 1$ for all $\omega \in \Omega$. Hence $\hat{\rho}$ satisfies the two constraints of \eqref{eq:fluid_approx} and is thus a feasible solution to \eqref{eq:fluid_approx}. Let $V_T^\text{Fluid}(\bar{S})$ be the optimal solution to \eqref{eq:fluid_approx}. By optimality of $V_T^\text{Fluid}(\bar{S})$ we have $V_T^\text{Fluid}(\bar{S}) \geq \arate T \sum_{i \in \bar{S}} \E_\rv \left[\hat{\rho}_i^\rv \left(p_i^\rv(\hat{\rho}) - \xi_i\right)\right] + \sum_{i \in \bar{S}} \xi_i$. We obtain from Lemma \ref{lemma:price_proba_bijection} that:
\begin{align*}
    p_i^\omega(\hat{\rho}) &= \frac{1}{\beta_p}\left(\ln \hat{\rho}_i - \ln \hat{\rho}_0 - q_i^\omega\right)\\
    &= \frac{1}{\beta_p}\left(\ln \frac{\frac{1}{\arate T} e^{q_i^\omega - \kappa_i}}{1-\frac{1}{\arate T} \sum_{i \in \bar{S}}e^{q_i^\omega - \kappa_i}} - q_i^\omega\right)\\
    &= \frac{1}{\beta_p}\left(-\kappa_i - \ln\left(\arate T - \sum_{j \in \bar{S}}e^{q_j^\omega - \kappa_j}\right)\right)
\end{align*}
and it follows, since $\E_\rv \left[\hat{\rho}_i^\rv\right] = 1/\arate T$, that:
\begin{align*}
    V_T^\text{Fluid} &\geq \arate T \sum_{i \in \bar{S}} \E_\rv \left[\hat{\rho}_i^\rv \left(p_i^\rv(\hat{\rho}) - \xi_i\right)\right] + \sum_{i \in \bar{S}} \xi_i\\
    &= \arate T \sum_{i \in \bar{S}} \E_\rv \left[\hat{\rho}_i^\rv p_i^\rv(\hat{\rho})\right]\\
    &= \sum_{i \in \bar{S}} \E_\rv \left[\frac{e^{q_i^\rv - \kappa_i}}{\beta_p}\left(-\kappa_i - \ln\left(\arate T - \sum_{j \in \bar{S}}e^{q_j^\rv - \kappa_j}\right)\right)\right]\\
    &= \frac{-1}{\beta_p}\sum_{i \in \bar{S}}\left(\kappa_i + \ln(\arate T)\right) + \sum_{i \in \bar{S}}\E_\rv \left[\frac{-e^{q_i^\rv - \kappa_i}}{\beta_p} \cdot \ln\left(1 - \frac{1}{\arate T}\sum_{j \in \bar{S}}e^{q_j^\rv - \kappa_j}\right)\right]\\
    &= \frac{-1}{\beta_p}\left(\kappa(\bar{S}) + |\bar{S}|\ln(\arate T)\right) + O\left(\frac{1}{T}\right)
\end{align*}
Introducing the term $v_T$ defined in Theorem \ref{theorem:asymptotic_costs}, we obtain:
$$V_T^\text{Fluid}(\bar{S}) \geq v_T(\bar{S}) - \frac{|\bar{S}|}{\beta_p} + O\left(\frac{1}{T}\right)$$
Recall from Theorem \ref{theorem:asymptotic_costs} that $\lim_{T \to \infty} V_T^*(\bar{S}) - v_T(S) = 0$. Combining the above results we have $\lim_{T \to \infty}  V_T^\text{Fluid}(\bar{S}) - V_T^*(\bar{S}) \geq -\frac{|\bar{S}|}{\beta_p} > 0$.
\end{cproof}

\subsection{Static Approximation} \label{apx:static_approx}
Another common approach to simplify the stochastic dynamic program is to constrain quantities such as the price vector $p$ and probability vector $\rho$ to be time-invariant. With this approximation, problem \eqref{eq:bellman_equation} reduces to:
\begin{align*}
    \max_{\rho \in (0,1)^{N\times|\Omega|}} \quad & \sum_{i \in \bar{S}} \sum_{t=1}^T \left(1-\arate \E_\rv\left[\rho_i^\rv\right]\right)^{T-t}  \E_\rv\left[\arate \rho_i^\rv p_i^\rv(\rho^\rv)\right] + \sum_{i \in \bar{S}} \left(1-\arate \E_\rv\left[\rho_i^\rv\right]\right)^{T} \xi_i\\
    \text{s.t.} \quad & \sum_{i \in \bar{S}} \rho_i^\omega \leq 1, \quad \forall \omega \in \Omega
\end{align*}
which can be further simplified to:
\begin{align}
    \max_{\rho \in (0,1)^{N\times|\Omega|}} \quad & \sum_{i \in \bar{S}} \frac{1-\left(1-\arate \E_\rv\left[\rho_i^\rv\right]\right)^{T}}{\E_\rv\left[\rho_i^\rv\right]} \E_\rv\left[ \rho_i^\rv p_i^\rv(\rho^\rv)\right] + \sum_{i \in \bar{S}} \left(1-\arate \E_\rv\left[\rho_i^\rv\right]\right)^{T} \xi_i \tag{$\mathcal{P}_\text{Static}$} \label{eq:static_approx}\\
    \text{s.t.} \quad & \sum_{i \in \bar{S}} \rho_i^\omega \leq 1, \quad \forall \omega \in \Omega \nonumber
\end{align}
We denote by $V_T^\text{Static}\left(\bar{S}\right)$ the optimal solution to problem \eqref{eq:static_approx}. As a constrained variant of the original problem, the optimal solution of \eqref{eq:static_approx} yields a lower bound for the original stochastic problem.\\

\begin{proposition}\label{prop:static_approx_asymptotical}
The optimal objective of \eqref{eq:static_approx} need not converge to the optimal objective of the stochastic problem as $T$ tends to $\infty$.
\end{proposition}

\begin{cproof}{Proposition}{\ref{prop:static_approx_asymptotical}}
Consider the special case where $|\Omega| = |\bar{S}| = 1$ and $\xi = q = 0$. Let $R^*(T)$ be the optimal objective of \eqref{eq:static_approx}. In this case, the problem reduces to:
\begin{align*}
    R^*(T) = \max_{\rho \in (0,1)} \left(1-(1-\arate \rho)^T\right) p(\rho)
\end{align*}
which, using Lemma \ref{lemma:price_proba_bijection}, further reduces to:
\begin{align*}
    R^*(T) = \max_{\rho \in (0,1)} \frac{-1}{\beta_p}\left(1-(1-\arate \rho)^T\right) \ln\left(\frac{1-\rho}{\rho}\right)
\end{align*}
We first note that $\ln\left(\frac{1-\rho}{\rho}\right) \leq \ln(1/\rho)$ for any $\rho \in (0,1)$. We then let $\alpha = \arate T \rho$ and maximize over $\alpha$ instead:
\begin{align}
    R^*(T) &\leq \max_{\alpha \in (0,\arate T)} \frac{-1}{\beta_p}\left(1-\left(1-\frac{\alpha}{T}\right)^T\right) \ln\left(\frac{\arate T}{\alpha}\right)\nonumber\\
    &\leq \max_{\alpha \in (0,\arate T)} \frac{-1}{\beta_p}\left(1-e^{-\alpha}\right) \left(\ln\left(\arate T\right) - \ln(\alpha)\right)  \label{eq:static_alpha}
\end{align}
where the last inequality follows since $R^*(T) \geq 0$, implying $\ln\left(\arate T\right) - \ln(\alpha) \geq 0$ for any optimal $\alpha$, and $\left(1-\frac{\alpha}{T}\right)^T \leq e^{-\alpha}$ for $T \geq 1$.
Let $\alpha^*(T)$ be the optimal solution to \eqref{eq:static_alpha}.
As $\lim_{T \to \infty} \alpha^*(T) = \infty$, there exists $T_0 \in \N$ such that $\alpha^*(T) \geq 2e$ for all $T \geq T_0$. It follows that for any $T \geq \max(T_0,1/\arate)$:
\begin{align*}
    R^*(T) &\leq \frac{-1}{\beta_p}\left(1-e^{-\alpha^*(T)}\right) \left(\ln\left(\arate T\right) - \ln\left(e\cdot \frac{\alpha^*(T)}{e}\right)\right)\\
    &\leq \frac{-1}{\beta_p} \left(\ln\left(\arate T\right) - 1\right) + \frac{1}{\beta_p}\left(1-e^{-2e}\right) \ln\left(2\right)
\end{align*}
Using Theorem \ref{theorem:asymptotic_costs} we conclude that $\lim_{T \to \infty} V_T^* - R^*(T) > 0$. Thus, we found a counterexample where the optimal objective of \eqref{eq:static_approx} does not converge to the optimal objective of the stochastic problem.
\end{cproof}

\subsection{Backward Lower Bound} \label{apx:backward_lower_bound}

Since selection probabilities decrease as new options are added, replacing $f_i^\omega(S)$ with $f_i^\omega(\bar{S})$ in \eqref{eq:bellman_equation} corresponds to a worst-case scenario. This substitution simplifies the dynamic program to a 1-dimensional recursion, and the total expected revenue can be expressed as the sum of the individual contribution of each option.\\
Let $l \in \R^{T \times N}$ be such that for all $i \in \bar{S}$:
\begin{align*}
    l_{0,i} &= \xi_i\\
    l_{t,i} &= l_{t-1,i} + \arate \cdot \E_\rv \left[ f_i^\rv(\hat{p}_t^\rv, \bar{S})(\hat{p}^\omega_{t,i} - l_{t-1,i}) \right], \quad \forall t \in \{1,\hdots,T(\ear)\}
\end{align*}
where
\begin{align*}
    \hat{p}_t^\omega& = \argmax_{p \in \R^N} \sum_{i \in \bar{S}} f_i^\omega(p, \bar{S}) \left(p_i - l_{t-1,i}\right), \quad \forall t \in \{0,\hdots,T(\ear)\}, \quad \forall \omega \in \Omega
\end{align*}
We define $V^L_t(S) = \sum_{i \in S} l_{t,i}$ for all $t \in \N$. This time, obtaining $V^L_{T(\ear)}(S)$ only requires solving $|\Omega| \cdot T(\ear)$ problems of the form of \eqref{eq:inner_problem}, but each maximization problem is slightly more complex as it involves $|S|$ options instead of one as in the case of the upper bound. The optimal value of the modified dynamic program yields a lower bound on the optimal objective of the original problem.

\begin{proposition} \label{prop:lower_bound}
For all $S \subseteq \bar{S}$ and $t \in \N$, $V^L_t(S) \leq V_t^*(S)$.
\end{proposition}

\subsection{Minimum Expected Empty Miles: Mixed Integer Linear Program}\label{subsection:min_empty_miles}

If a load is the last load of an option, its expected empty miles to the next load are estimated using historical data. Namely, let $\hat{e} \in \R^{|\Omega|}$ be the vector such that $\hat{e}_\omega$ is the average historical empty miles of loads whose destination market is $\omega$. Thus, following the delivery of load $k$ included in option $i$, the expected empty miles to the next load is either $\hat{e}_{k_\texttt{dst}}$ if load $k$ is the last load of the option $i$, or $e_{k,l}$, the distance between the drop-off location of load $k$ and the pick-up location of load $l$, if load $l$ is following load $k$ in option $i$.

We formulate a Mixed Integer Linear Program that minimizes the expected empty miles following load drop-offs by forming bundles of size at most 2:
\begin{align*}
    \min_{x \in \{0,1\}^{L \times L}} \quad& \sum_{k \in \mathcal{L}} \left[\hat{e}_{k_\texttt{dst}} + \sum_{l \in \mathcal{L}\backslash\{k\}} x_{k,l} \left( e_{k,l} - \hat{e}_{k_\texttt{dst}} \right) \right] \tag{$\mathcal{P}_\text{EM}$} \label{eq:minimum_empty_miles}\\
    \text{s.t.} \quad& \sum_{l \in \mathcal{L}\backslash\{k\}} x_{k,l} + x_{l,k} \leq 1, \quad \forall k \in \mathcal{L}
\end{align*}
where the constraint ensures that each load is preceded or succeeded by at most one other load, or equivalently, paired to at most one other load. Given a solution $x^*$ of \eqref{eq:minimum_empty_miles}, the corresponding set $S_\text{EM}$ is formed by bundles $(k,l)$ for all $(k,l) \in \mathcal{L}^2$ if $x^*_{k,l} = 1$ and single loads $k \in \mathcal{L}$ for all $k$ such that $\sum_{l \in \mathcal{L}\backslash\{k\}} x^*_{k,l} + x^*_{l,k} = 0$.

\newpage

\subsection{Greedy Algorithm}

\begin{algorithm}[h]
\caption{Greedy Algorithm}\label{alg:greedy_alg}
\begin{algorithmic} \small
\Require $V$: value function approximation, $\mathcal{L}$: set of products, $K_s$: maximum number of bundles, $K_b$: maximum bundle size
\State $\mathcal{O} \leftarrow \left\{(l_1, \hdots, l_k) \in \mathcal{L}^k\;|\; k\leq K_b \text{ and } l_a \neq l_b, \text{ for } a\neq b\right\}$  \Comment{$\mathcal{O}$ is the set of all options of size at most $K_b$}
\For{$l \in \mathcal{L}$}
    \State $\Delta_l \leftarrow V_{T(\ear)}(S_0) - V_{T(\ear)}(S_0\backslash\{l\})$
\EndFor
\State $S \leftarrow \{\}$, $I \leftarrow \{\}$, $d \leftarrow 0$ \Comment{$I$ is the set of items included in $S$}
\For{$i \in \mathcal{O}$}
    \State $\bar{x}_i \leftarrow \mathbb{E}_\rv \left[e^{q_i^\rv + \beta_p \cdot \sum_{l \in i} \Delta_l}\right]$
\EndFor
\For{$i \in \text{arg\_reverse\_sort}(\bar{x})$}
\If{($i \cap I = \emptyset$) and ($\#i=1$ or $d<K_s$)}
    \State $S \leftarrow S \cup \{i\}$
    \For{$l \in i$}
        \State $I \leftarrow I \cup \{l\}$
    \EndFor
    \If{$\#i > 1$}
        \State $d \leftarrow d + 1$
    \EndIf
\EndIf
\EndFor
\State \Return $S$
\end{algorithmic}
\end{algorithm}

\subsection{More on Industry-Driven Case Study} \label{apx:case_study}

\begin{figure}[h]
    \centering
    \begin{subfigure}{0.6\linewidth}
        \centering
        \includegraphics[width=\linewidth]{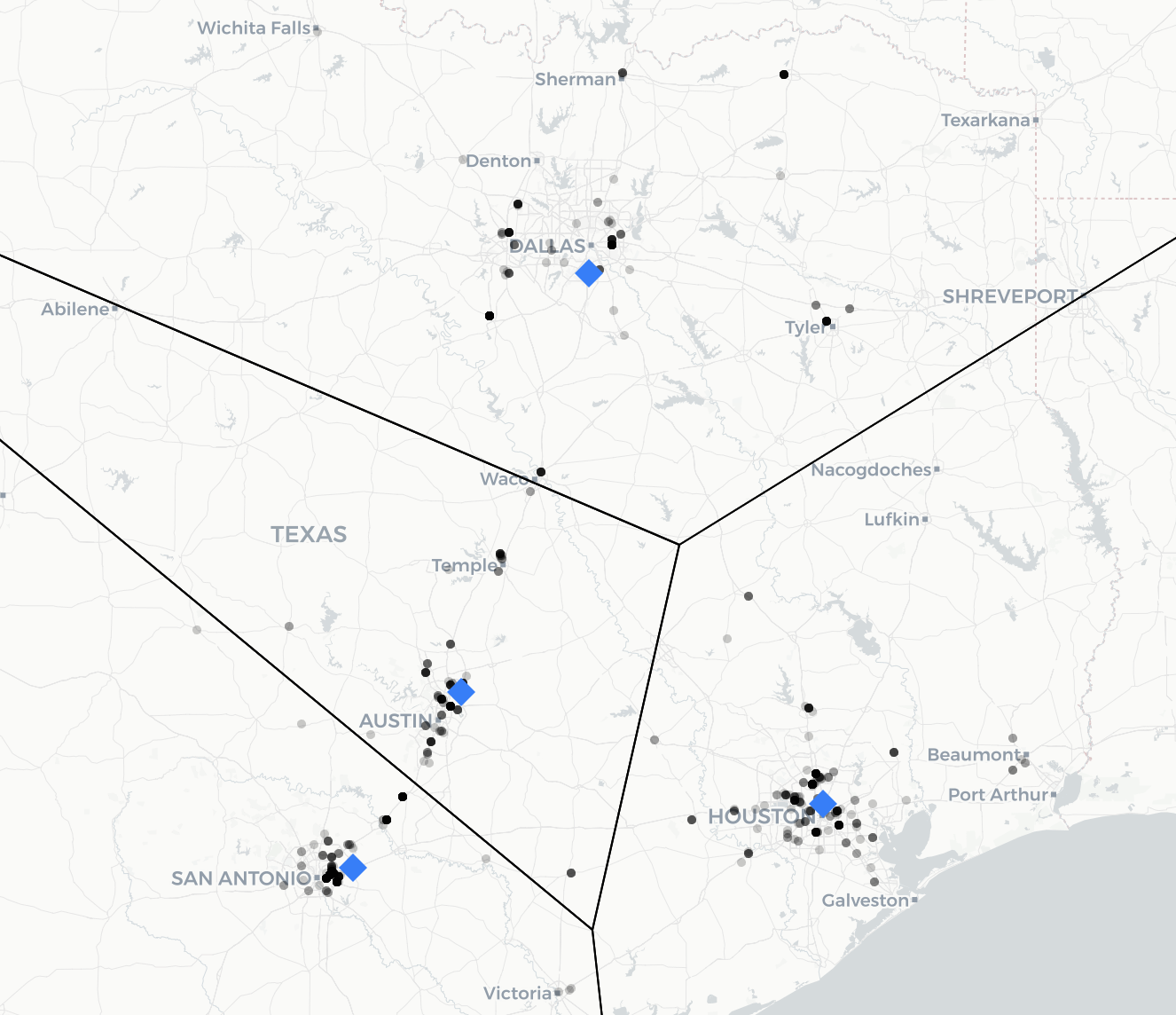}
    \end{subfigure}
    \caption{Distribution of the drop-off locations of historical loads (black) and centroids of the K-means clustering (blue).}
    \label{fig:voronoi}
\end{figure}

We give additional information on the experiments ran in Section \ref{section:ExperimentalResults}. Figure \ref{fig:voronoi} shows the distribution of drop-off locations and the centroids of the 4 markets obtained using K-means algorithm. The black solid lines split the space into 4 markets: San Antonio, Austin, Dallas, and Houston.

In Section \ref{section:ExperimentalResults}, we use a MNL model with parameters fitted to a logit model by our industry partner. Table \ref{tab:simulation_results_logit} shows that running the simulation with a sequence of binary decisions (logit model) yields similar results to the simulation with a single multiple choice decision (MNL model). Since the probability of accepting each option is higher with the logit model, the sequential binary decision simulation yields more optimistic results in terms of cost per mile and unmet deadline across all frameworks, but the discussion in Section \ref{section:ExperimentalResults} remains unchanged.

\begin{table}
\centering
{
\small
\begin{tabular}{|l|c|c||c|c|c|}
\hline
Framework & \makecell{Bunlding \\ method} & Pricing & \makecell{Cost/loaded mile \\ (\$/mi.)}& \makecell{Avg. empty miles \\ (mi.)} & Unmet deadlines\\
\hline
No bundle & - & - & $2.7284\pm0.0097$ & $31.6096\pm0.4478$ & $5.19\%\pm0.29$\\
\hline
\multirow{4}{*}{Rolling horizon}
& \multirow{2}{*}{Empty miles} 
& Linear & $2.9693\pm0.0112$ & $28.2238\pm0.3896$ & $9.61\%\pm0.39$\\
    \cline{3-6}
& & Custom & $2.8284\pm0.0091$ & $27.9261\pm0.3765$ & $3.32\%\pm0.24$\\
\cline{2-6}
& \multirow{2}{*}{Greedy Alg.} 
& Linear & $2.6847\pm0.0093$ & $28.8958\pm0.3673$ & $4.39\%\pm0.27$\\
    \cline{3-6}
& & Custom & $2.6795\pm0.0096$ & $29.1673\pm0.3759$ & $2.97\%\pm0.22$\\
\hline
\multirow{2}{*}{Personalized}
& \multirow{2}{*}{-}
& Linear & $2.6497\pm0.0071$ & $\bm{23.2409\pm0.2475}$ & $4.83\%\pm0.28$\\
    \cline{3-6}
& & Custom & $\bm{2.6042\pm0.0082}$ & $24.4197\pm0.2724$ & $\bm{2.21\%\pm0.19}$\\
\hline
\end{tabular}
}
\caption{Two-year simulation with sequential decisions: operational results.}
\label{tab:simulation_results_logit}
\end{table}

\begin{figure}[H]
    \centering
    \begin{subfigure}{0.45\linewidth}
        \centering
        \includegraphics[width=\linewidth]{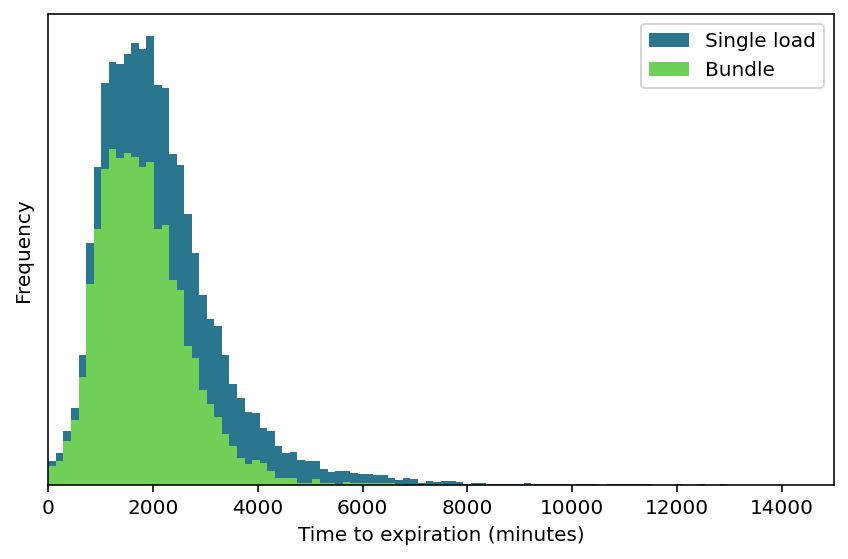}
        \caption{Linear pricing.}\label{fig:acceptance_linear}
    \end{subfigure}
    \begin{subfigure}{0.45\linewidth}
        \centering
        \includegraphics[width=\linewidth]{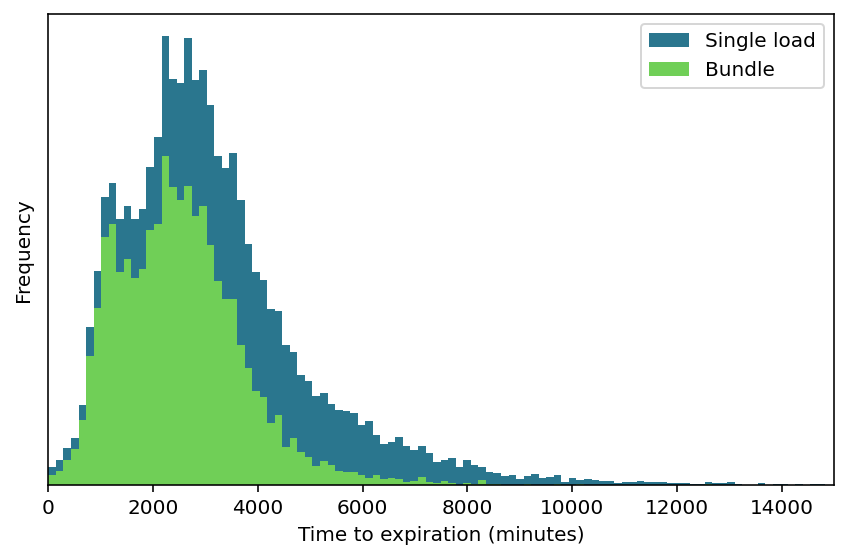}
        \caption{Custom pricing.}\label{fig:acceptance_premium}
    \end{subfigure}
    \caption{Distribution of remaining time before expiration of loads booked by a carrier for personalized bundles.}
    \label{fig:acceptance_distrib}
\end{figure}

Figure \ref{fig:acceptance_distrib} illustrates the cannibalism effect on personalized bundles. In figure \ref{fig:acceptance_linear}, underpriced bundles negatively impact both the early booking of bundles and single loads. These underpriced bundles fail to attract interest of carriers, and as a result, individual loads remain overshadowed by these bundles for longer. This leads to an accumulation of loads on the platform and options being accepted close to the deadline. On the other hand, custom pricing (figure \ref{fig:acceptance_premium}) allows individual loads and bundles to be accepted earlier, reducing the average number of loads on the platform from 123 to 95 and allowing it to better handle loads supplied shortly before their deadline.

\subsection{Proofs of Section \ref{section:DP}}

\begin{lemma}\label{lemma:price_proba_bijection}
Under the MNL model, the price vector $p$ is uniquely determined by the probability vector of choice $\rho$. In particular, the price of option $i$ is given by:
\begin{equation}
p_i(\rho) = \frac{1}{\beta_p}\left(\ln{\rho_i} - \ln{\rho_0} -q_i\right) \label{eq:price}
\end{equation}
\end{lemma}

\begin{cproof}{Lemma}{\ref{lemma:price_proba_bijection}}
Under the MNL model, the probability that option $i$ is selected given a price vector $p$ is:
\begin{equation} \label{eq:rho}
    \rho_i = \frac{e^{q_i + \beta_p p_i}}{1 + \sum_{j=1}^N e^{q_j + \beta_p p_j}}
\end{equation}
and the probability that no option is selected is:
$$\rho_0 = \frac{1}{1 + \sum_{j=1}^N e^{q_j + \beta_p p_j}}$$
As $\rho_0 > 0$, we can express the ratio between the two probabilities as
$$\frac{\rho_i}{\rho_0} = e^{q_i + \beta_p p_i}$$
Hence, by taking the logarithm to the above equation, we can write $p_i$ as a function of $\rho$ as given in \eqref{eq:price}.
\end{cproof}

\begin{lemma}\label{lemma:joint_convexity}
Let $\phi(x_1, x_2) = x_1 \left[\log x_1 - \log x_2\right]$. Then, $\phi$ is jointly convex in ($x_1$, $x_2$) in the region $(0,1)^2$.
\end{lemma}

\begin{cproof}{Lemma}{\ref{lemma:joint_convexity}}
The Hessian of $\phi$ is $H_\phi(x) = \left[\begin{matrix}
\frac{1}{x_1} & \frac{1}{x_2}\\
\frac{1}{x_2} & \frac{x_1}{x_2^2}
\end{matrix}\right]$. Then for any $y \in \R^2$ we have
$$y^\top H_\phi(x) y = \frac{y_1^2}{x_1} + \frac{y_1 y_2}{x_2} + \frac{y_1 y_2}{x_2} + \frac{y_2^2 x_1}{x_2^2} = \left(\frac{y_1}{\sqrt{x_1}}+\frac{y_2\sqrt{x_1}}{x_2}\right)^2$$
which is strictly positive for any $(x_1, x_2) \in (0,1)^2$. Thus the Hessian of $\phi$ is positive semi-definite and $\phi$ is jointly convex.
\end{cproof}

\begin{cproof}{Lemma}{\ref{lemma:inner_problem_solution}}
For readability we drop the $\omega$ superscript in this proof.
Using Lemma \ref{lemma:price_proba_bijection}, we can express the expected revenue function as a function of the probability vector $\rho$:
\begin{align*}
    r(\rho) &\triangleq \sum_{i \in S} \rho_i \cdot \left(p_i(\rho)-\Delta_i\right)\\
    &= \frac{1}{\beta_p}\sum_{i \in S} \rho_i \cdot \left(\ln{\rho_i} - \ln{\rho_0} -q_i-\beta_p\Delta_i\right)
\end{align*}
and write \eqref{eq:inner_problem} as a maximization problem over $\rho$:
\begin{align*}
    \max_{\rho \in (0,1)^{|S|}} r(\rho)
\end{align*}
For all $i \in S$, $\rho_i \left(\ln{\rho_i} - \ln{\rho_0}\right)$ is jointly convex in $(\rho_i, \rho_0)$ according to Lemma \ref{lemma:joint_convexity} and $\rho_i \left(-q_i-\beta_p\Delta_i\right)$ is lineary in $\rho_i$. It follows that the $r$ is concave in $\rho$ as a sum of convex functions multiplied by the negative scalar $\beta_p$.\\
As $\rho_0 = 1 - \sum_{j \in S} \rho_j$, the partial derivative of $p_i$ with respect to $\rho_j$ is given by:
$$\frac{\partial p_i}{\partial \rho_j}(\rho) = \frac{1}{\beta_p}\left(\frac{1}{\rho_0} + \frac{1}{\rho_i} \mathds{1}_{\{i=j\}}\right)$$
and the partial derivative of the expected revenue with respect to $\rho_j$ is:
\begin{align*}
    \frac{\partial r}{\partial \rho_j}(\rho) &= p_j(\rho) - \Delta_j + \sum_{i \in S} \rho_i \frac{\partial p_i}{\partial \rho_j}(\rho)\\
    &= \frac{1}{\beta_p}\left(\ln{\rho_j} - \ln{\rho_0} -q_j\right) - \Delta_j + \sum_{i \in S} \frac{\rho_i}{\beta_p}\left(\frac{1}{\rho_0} + \frac{1}{\rho_i} \mathds{1}_{\{i=j\}}\right)\\
    &= \frac{1}{\beta_p}\left(\ln{\rho_j} - \ln{\rho_0} -q_j\right) - \Delta_j + \frac{1}{\beta_p}\left(1 + \frac{\sum_{i \in S} \rho_i}{\rho_0}\right)\\
    &= \frac{1}{\beta_p}\left(\ln{\rho_j} - \ln{\rho_0} -q_j\right) - \Delta_j + \frac{1}{\beta_p \rho_0}\\
\end{align*}
Since $r$ is concave in $\rho$ we can use the first-order condition and equate the partial derivative of $r$ at $\rho^*$ and obtain:
\begin{equation}
\rho^*_j = \rho^*_0 e^{q_j + \beta_p \Delta_j - \frac{1}{\rho^*_0}} \label{eq:rhostar}
\end{equation}
Since $\rho_0^* > 0$, we can define $m = \frac{1}{\rho^*_0}$ and reformulate $\rho^*_0 = 1-\sum_{i \in S} \rho^*_i$ as:
\begin{align*}
    \frac{1}{m} &= 1-\sum_{i \in S} \frac{1}{m} e^{q_i + \beta_p \Delta_i - m}\\
    e^m &= me^m -\sum_{i \in S} e^{q_i + \beta_p \Delta_i}\\
    (m-1)e^{m-1} &= e^{-1} \sum_{i \in S} e^{q_i + \beta_p \Delta_i}\\
    m-1 &= W\left(\sum_{i \in S} e^{q_i + \beta_p \Delta_i - 1}\right)\\
    m &= 1+\Gamma
\end{align*}
Thus $\rho_0^* = 1/(1+\Gamma)$. Substituting \eqref{eq:rhostar} into \eqref{eq:price} we obtain:
\begin{align*}
    p_i(\rho^*) &= \frac{1}{\beta_p}\left(\ln{\rho^*_i} - \ln{\rho^*_0} -q_i\right)\\
    &= \frac{1}{\beta_p}\left(\ln{\rho^*_0} + q_i + \beta_p \Delta_i - \frac{1}{\rho_0^*} - \ln{\rho^*_0} -q_i\right)\\
    &= \Delta_i - \frac{1}{\beta_p \rho_0^*}\\
    &= \Delta_i - \frac{1+\Gamma}{\beta_p}
\end{align*}
Therefore the expected revenue with optimal $p^*$ and $\rho^*$ is:
\begin{align*}
    r^* &= \sum_{i \in S} \rho^*_i\left(p_i(\rho^*)-\Delta_i\right)\\
     &= - \frac{1+\Gamma}{\beta_p} \sum_{i \in S} \rho^*_i\\
     &= - \frac{1+\Gamma}{\beta_p} \left(1-\rho^*_0\right)\\
     &= -\frac{\Gamma}{\beta_p}
\end{align*}
Finally, as $e^\Gamma = \Gamma^{-1} \sum_{i \in S} e^{q_i + \beta_p \Delta_i - 1}$, substituting optimal prices into \eqref{eq:rho} gives:
\begin{align*}
    \rho_i^* &=\frac{e^{q_i + \beta_p \Delta_i - 1 - \Gamma}}{1 + \sum_{j \in S} e^{q_j + \beta_p \Delta_j - 1 - \Gamma}}\\
     &= \frac{e^{q_i + \beta_p \Delta_i}}{e^{1+\Gamma} + \sum_{j \in S} e^{q_j + \beta_p \Delta_j}}\\
     &= \frac{e^{q_i + \beta_p \Delta_i}}{\left(1+\Gamma^{-1}\right)\sum_{j \in S} e^{q_j + \beta_p \Delta_j}}\\
     &= \frac{\Gamma}{1+\Gamma} \frac{e^{q_i + \beta_p \Delta_i}}{\sum_{j \in S} e^{q_j + \beta_p \Delta_j}}
\end{align*}
\end{cproof}

\begin{lemma} \label{lemma:nonexpansion_W_exp}
The function $g: x \rightarrow W\left(e^{x}\right)$ is strictly increasing and nonexpansive over $\R$. If the domain of $g$ is restrained to the interval $[-M, M]$ for $M > 0$, then $g$ is a contraction with constant $\left(1+\frac{1}{W(e^M)}\right)^{-1}$ and an expansion with constant $\left(1+\frac{1}{W(e^{-M})}\right)^{-1}$.
\end{lemma}

\begin{cproof}{Lemma}{\ref{lemma:nonexpansion_W_exp}}
Note that $g$ is $C^1$ over $\R$ as the composition of $W$, $C^1$ over $\R^+$, and $\exp$, $C^1$ over $\R$ with value in $\R^+$. As shown in \cite{Lambert_W_Function}, the derivative of the Lambert $W$ function is:
\begin{align*}
    W'(z) &= \frac{W(z)}{z(1+W(z))}
\end{align*}
for all $z \notin \left(-e^{-1}, 0\right)$. Using the identity $W(z)e^{W(z)} = z$, we obtain for all $z > 0$:
\begin{align*}
    W'(z) = \frac{1}{z+e^{W(z)}}
\end{align*}
It follows that the derivative of $g$ is:
\begin{align*}
    g'(x) &= \frac{1}{1+\frac{e^{W(e^x)}}{e^x}} = \frac{1}{1+\frac{1}{W(e^x)}}
\end{align*}
and $g'(x) \in (0,1)$ for all $x \in \R$.
For any closed interval $[x,y] \subseteq \R$, the Mean Value Theorem states that there exists $c \in (x,y)$ such that $g(y) - g(x) = g'(c)(y-x)$. Thus:
\begin{align*}
    |g(y) - g(x)| &= g'(c) \cdot |y-x|\\
    &\leq |y-x|
\end{align*}
showing that $g$ is a nonexpansive mapping over $\R$.\\
Suppose now that the domain of $g$ is restrained to the interval $[-M, M]$, the maximum of $g'$ is attained in $M$. Therefore the Mean Value Theorem now gives for any interval $[x,y] \subseteq [-M, M]$:
\begin{align*}
    |g(y) - g(x)| &= g'(c) \cdot |y-x|\\
    &\leq \left(1+\frac{1}{W(e^M)}\right)^{-1} |y-x|
\end{align*}
and we conclude that $g$ is a contraction with constant $\left(1+\frac{1}{W(e^M)}\right)^{-1}$.\\
Similarly, the minimum of $g'$ is attained in $-M$ so $g$ is an expansion with constant $\left(1+\frac{1}{W(e^{-M})}\right)^{-1}$.
\end{cproof}


\begin{lemma} \label{lemma:convergence_rate_noisy_sequence}
Let $(r_n)_{n\in\N}$ and $(\epsilon_n)_{n\in\N}$ satisying $r_n \in [1,n]$ for all $n \in \N$ and $\epsilon_n = O(\log^k(n)/n^2)$ for $k \in \N$. Let $(y_n)_{n\in\N}$ be a sequence such that:
$$y_{n} = \left(1-\frac{r_n}{n}\right)y_{n-1} + \epsilon_n$$
Then $|y_n| = O\left(\frac{\ln^{k+1}(n)}{n}\right)$.
\end{lemma}
\begin{cproof}{Lemma}{\ref{lemma:convergence_rate_noisy_sequence}}
We first show that $\ln(n+1) \geq \ln(n) + \frac{1}{n+1}$ for all $n \in \N$. Let $f(x) = \ln(x+1)-\ln(x)-\frac{1}{x+1}$. Then $f'(x) = \frac{-1}{x(x+1)^2} \leq 0$ for all $x > 0$. Since $\lim_{x \to \infty} f(x) = 0$ we conclude that $f(x) \geq 0$ for all $x \geq 0$. In particular, $\ln(n+1) \geq \ln(n) + \frac{1}{n+1}$ for all $n \in \N$.\\
Since $\epsilon_n = O(\ln^k(n)/n^2)$, there exists $M_\epsilon > 0$ and $n_\epsilon \geq 2$ such that $|\epsilon_n| \leq M_\epsilon\frac{\ln^k(n)}{n^2}$ for all $n \geq n_\epsilon$.\\
Suppose that, for a fixed $n \geq n_\epsilon$, there exists $M_y > M_\epsilon$ such that $|y_n| \leq M_y \frac{\ln^{k+1} n}{n}$. Then:
\begin{align*}
    |y_{n+1}| &\leq \left(1-\frac{r_{n+1}}{n+1}\right)M_y \frac{\ln^{k+1} n}{n} + M_\epsilon\frac{\ln^k(n+1)}{(n+1)^2}\\
    &\numsym{1}{\leq} \left(\frac{n}{n+1}\right)M_y \frac{\ln^{k+1} n}{n} + M_\epsilon\frac{\ln^k(n+1)}{(n+1)^2}\\
    &\leq \frac{\ln^{k} (n+1)}{n+1}\left(\left(M_y-M_\epsilon\right)\ln(n+1) + M_\epsilon\left(\ln(n)+\frac{1}{n+1}\right)\right)\\
    &\numsym{2}{\leq} \frac{\ln^{k} (n+1)}{n+1}\left(\left(M_y-M_\epsilon\right)\ln(n+1) + M_\epsilon\ln(n+1)\right)\\
    &= M_y \frac{\ln^{k+1}(n+1)}{n+1}
\end{align*}
where (1) holds since $r_{n+1} \in [1,n+1]$ and (2) follows from $\ln(n+1) \geq \ln(n) + \frac{1}{n+1}$. Since $|y_{n_\epsilon}| = M_y \frac{\ln^{k+1} n_\epsilon}{n_\epsilon}$ with $M_y = |y_{n_\epsilon}|\frac{n_\epsilon}{\ln^{k+1} n_\epsilon}$, we obtain $|y_n| \leq M_y \frac{\ln^{k+1} n}{n}$ for all $n \geq {n_\epsilon}$ and we conclude that $|y_n| = O\left(\frac{\ln^{k+1}(n)}{n}\right)$.
\end{cproof}


\begin{lemma} \label{lemma:asymptotic_costs}
Let $\bar{S} \in \mathcal{P}(\mathcal{L})$ and, for all $S \subseteq \bar{S}$, let $(z^S_t)_{t \in \N}$ be a sequence such that:
\begin{equation}\label{eq:revenue_recurrence}
     z^S_{t+1} = z^S_t - \frac{\arate}{\beta_p} \E_{\rv}\left[ W\left( \sum_{i \in S} e^{q_i^\rv + \beta_p \left(z^S_t - z_t^{S\backslash\{i\}}\right)-1}\right) \right], \quad \forall t \in \N
\end{equation}
and $(v^S_t)_{t \in \N}$ be a sequence defined as:
$$v^S_t = \frac{-1}{\beta_p}\left( |S|(\ln (\arate t) - 1) + 
\sum_{i \in S} \ln \E_{\rv}\left[e^{q_i^\rv} \right]\right), \quad \forall t \in \N$$
Then $z^S_t - v^S_t = O\left(\frac{\ln^{|S|}(t)}{t}\right)$ for all $S \subseteq \bar{S}$.
\end{lemma}

\begin{cproof}{Lemma}{\ref{lemma:asymptotic_costs}}
Let $(\epsilon_t)_{t \in \N}$ be such that $\epsilon_t^S = z_t^S - v_t^S$ for all $t \in \N$ and for all $S \subseteq \bar{S}$.
Let's consider a set $S \subseteq \bar{S}$ with cardinality $k$. Let the hypothesis of induction be that $\lim_{t \to \infty} \epsilon_t^{S\backslash\{i\}} = O\left(\frac{\ln^{k-1}(t)}{t^2}\right)$ for all $i \in S$. For notation clarity we drop the $S$ superscript, replace the $S\backslash\{i\}$ superscript with $-i$, and let $\E$ represent the expectation taken over $\rv$ unless stated otherwise. To demonstrate the result, we show that:
\begin{enumerate}
    \item the mappings $F_t: \epsilon_t \to \epsilon_{t+1}$ are nonexpansive with unique fixed point $\epsilon_t^*$
    \item the sequence $(\epsilon_t^*)_{t \in \N}$ converges to 0
    \item the sequence $(\epsilon_t)_{t \in \N}$ is bounded
    \item the mappings $F_t$, with domain restricted to the range of $(\epsilon_t)_{t \in \N}$, are contractions with parameter $\gamma_t$
    \item the sequence $(\epsilon_t)_{t \in \N}$ converges to 0
    \item $(\epsilon_t)_{t \in \N}$ is $O\left(\frac{\ln^{k}(t)}{t}\right)$
\end{enumerate}
Note that:
\begin{align*}
    z_t - z_t^{-i} &= \frac{-1}{\beta_p}\left( \ln (\arate t) - 1 + \ln \E\left[e^{q_i^\omega} \right] \right) + \epsilon_t - \epsilon_t^{-i}\\
    v_t - v_{t+1} &= \frac{k}{\beta_p} \ln \left(\frac{t+1}{t}\right)
\end{align*}
We can express $z_{t+1}$ as a function of $v_{t+1}$ and $\epsilon_{t+1}$ as follows:
\begin{align*}
    z_{t+1} &= z_t - \frac{\arate}{\beta_p} \E\left[ W\left( \sum_{i \in S} e^{q_i^\rv + \beta_p (z_t - z_t^{-i}) - 1}\right)\right]\\
    &= v_t + \epsilon_t - \frac{\arate}{\beta_p} \E\left[ W\left( \sum_{i \in S} e^{q_i^\rv - \ln(\arate t) - \ln \E_{\rv'}\left[e^{q_i^{\rv'}} \right] +\beta_p \left( \epsilon_t - \epsilon_t^{-i} \right)}\right)\right]\\
    &= v_{t+1} + \frac{k}{\beta_p} \ln \left(\frac{t+1}{t}\right) + \epsilon_t - \frac{\arate}{\beta_p} \E\left[ W\left( \sum_{i \in S} e^{q_i^\rv - \ln(\arate t) - \ln \E_{\rv'}\left[e^{q_i^{\rv'}} \right] +\beta_p \left( \epsilon_t - \epsilon_t^{-i} \right)}\right)\right]\\
    &= v_{t+1} + \frac{k}{\beta_p} \ln \left(\frac{t+1}{t}\right) + \epsilon_t - \frac{\arate}{\beta_p} \E\left[  W\left(\frac{e^{\beta_p \epsilon}}{\arate t} \sum_{i \in S} \frac{e^{q_i^\rv - \beta_p \epsilon_t^{-i}}}{\E_{\rv'}\left[e^{q_i^{\rv'}} \right]}\right)\right]\\
    &= v_{t+1} + \epsilon_{t+1}
\end{align*}
For all $t \in \N^*$, lets define:
\begin{align*}
    &h_t^\omega = \sum_{i \in S} \frac{e^{q_i^\omega - \beta_p \epsilon_t^{-i}}}{\E\left[e^{q_i^\rv} \right]}
    &F_t:\epsilon \to \epsilon + \frac{k}{\beta_p} \ln \left(\frac{t+1}{t}\right) - \frac{\arate}{\beta_p} \E\left[  W\left(\frac{h_t^\rv}{\arate t} e^{\beta_p \epsilon} \right)\right]
\end{align*}
Note that $F_t - Id_\R$ is strictly decreasing in $\epsilon$, $\lim_{\epsilon \to -\infty} F_t(\epsilon) - \epsilon = \infty$ and $\lim_{\epsilon \to \infty} F_t(\epsilon) - \epsilon = \frac{k}{\beta_p} \ln \left(\frac{t+1}{t}\right) < 0$. Therefore $F_t$ has a unique fixed point $\epsilon_t^*$.\\
We claim that $F_t$ is nonexpansive for all $t \in \N^*$. Let's fix $t \in \N$. Using Lemma \ref{lemma:nonexpansion_W_exp}, we know that the mapping $z \to W(e^z)$ is strictly increasing and nonexpansive. If follows that for any $b\geq a$, we have $\beta_p a \geq \beta_p b$ and:
\begin{align*}
    W\left(\frac{h_t^\omega}{\arate t} e^{\beta_p a} \right) - W\left(\frac{h_t^\omega}{\arate t} e^{\beta_p b} \right) &= W\left(e^{\beta_p a + \ln \frac{h_t^\omega}{\arate t}} \right) - W\left(e^{\beta_p b + \ln \frac{h_t^\omega}{\arate t}} \right)\\
    &\leq \beta_p a - \beta_p b
\end{align*}
Taking the expectation on both sides and multiplying by $-\arate/\beta_p$ we obtain:
\begin{equation} \label{eq:ab_contraction}
     \frac{\arate}{\beta_p} \E\left[  W\left(\frac{h_t^\rv}{t} e^{\beta_p b} \right)\right] - \frac{\arate}{\beta_p} \E\left[  W\left(\frac{h_t^\rv}{t} e^{\beta_p a} \right)\right] \leq \arate \left(b - a\right) \leq b-a
\end{equation}
Since the sign of both sides is identical we obtain:
\begin{align*}
    \left|F_t(b) - F_t(a)\right| &= \left| b - a  - \frac{\arate}{\beta_p} \E\left[  W\left(\frac{h_t^\rv}{\arate t} e^{\beta_p b} \right)\right] + \frac{\arate}{\beta_p} \E\left[  W\left(\frac{h_t^\rv}{\arate t} e^{\beta_p a} \right)\right] \right|\\
    &\leq \left|b - a\right|
\end{align*}
and $F_t$ is indeed a nonexpansive mapping.

We have that $(h_t^\omega)_{t \in \N^*}$ is bounded since $\epsilon_t^{-i} \to 0$ for all $i \in S$ and as $\Omega$ is finite the Dominated convergence theorem applies and we obtain:
\begin{align*}
    \lim_{t \to \infty} \E\left[h_t^\rv\right] &= \E\left[ \lim_{t \to \infty} h_t^\rv\right]\\
     &= \E\left[ \sum_{i \in S} \frac{e^{q_i^\rv}}{\E_{\rv'}\left[e^{q_i^{\rv'}} \right]} \right]\\
     &= k
\end{align*}
The Taylor series of $W$ gives $W(z) \underset{y \to 0}{=} z + O(z^2)$ so we have:
\begin{align*}
F_t(\epsilon) - \epsilon &= \frac{k}{\beta_p}\frac{1}{t} - \frac{\arate \E\left[h_t^\rv\right]}{\beta_p} \frac{e^{\beta_p \epsilon}}{\arate t} + O\left(\frac{1}{t^2}\right)\\
&= \frac{k}{\beta_p}\frac{1}{t} \left(1 - e^{\beta_p \epsilon} \right) + o\left(\frac{1}{t}\right)
\end{align*}
and we conclude that $\lim_{t \to \infty} \epsilon_t^* = 0$.\\
Let $\delta > 0$. Since $\epsilon_t^*$ converges to $0$, there exists $N \in \N$ such that $|\epsilon_t^*| < \delta$ for all $t \geq N$.\\
We'll now show that $(\epsilon_t)_{t \in \N}$ is bounded. Suppose by contradiction that $(\epsilon_t)_{t \in \N}$ diverges. For any $t \geq N$, two case arise:
\begin{itemize}
    \item if $\epsilon_t \leq \epsilon_t^*$, then $\epsilon_{t+1} \geq \epsilon_t$
    \item if $\epsilon_t > \epsilon_t^*$, then $\epsilon_t > \epsilon_{t+1} \geq \epsilon_t+\frac{k}{\beta_p} \ln \left(\frac{t+1}{t}\right) > \epsilon_t^*+\frac{k}{\beta_p} \ln \left(\frac{t+1}{t}\right) > -\delta+\frac{k}{\beta_p} \ln \left(\frac{N+1}{N}\right)$
\end{itemize}
Hence, since $(\epsilon_t)_{t\in \N}$ diverges, it must be an alternating sequence. However for all $t > N$, $\epsilon_t$ is lower bounded by $ -\delta+\frac{k}{\beta_p} \ln \left(\frac{N+1}{N}\right)$, and the nonexpansiveness of $F_t$ implies that, for any $t$ such that $\epsilon_t \leq \epsilon_t^*$, we have:
\begin{align*}
    |\epsilon_{t+1}-\epsilon_t^*| &= |F_t(\epsilon_t) - F_t(\epsilon_t^*)|\\
    &\leq |\epsilon_t - \epsilon_t^*|\\
    &\leq |\epsilon_t| + |\epsilon_t^*|\\
    &\leq 2\delta - \frac{k}{\beta_p} \ln \left(\frac{N+1}{N}\right)
\end{align*}
so terms above $\epsilon_t^*$ are upper bounded, which contradict that $(\epsilon_t)_{t\in \N^*}$ diverges. Therefore there exist $M \in \R^+$ such that $|\epsilon_t| \leq M$ for all $t \in \N^*$.\\
Let's define:
$$r_t^\omega = \left(1+\frac{1}{W\left(\frac{h_t^\omega}{\arate t}e^{\beta_p M}\right)}\right)^{-1}$$
With the domain of $F_t$ restained to $[-M, M]$, we can use the last part of Lemma \ref{lemma:nonexpansion_W_exp} and for any $-M \leq a \leq b \leq M$, we have $\beta_p a \geq \beta_p b$ so:
\begin{align*}
    W\left(\frac{h_t^\omega}{\arate t} e^{\beta_p a} \right) - W\left(\frac{h_t^\omega}{\arate t} e^{\beta_p b} \right) &= W\left(e^{\beta_p a + \ln \frac{h_t^\omega}{\arate t}} \right) - W\left(e^{\beta_p b + \ln \frac{h_t^\omega}{\arate t}} \right)\\
    &\geq r_t^\omega\left(\beta_p a - \beta_p b\right)
\end{align*}
Taking the expectation on both sides and multiplying by $-\arate/\beta_p$ we obtain:
\begin{equation*}
     \frac{\arate}{\beta_p} \E\left[  W\left(\frac{h_t^\rv}{\arate t} e^{\beta_p b} \right)\right] - \frac{\arate}{\beta_p} \E\left[  W\left(\frac{h_t^\rv}{\arate t} e^{\beta_p a} \right)\right] \geq \arate \E\left[r_t^\rv\right]\left(b-a\right)
\end{equation*}
Combined with \eqref{eq:ab_contraction} and the fact that the sign of both sides of the equation above is identical we now obtain:
\begin{align*}
    \left|F_t(b) - F_t(a)\right| &= \left| b - a  - \frac{\arate}{\beta_p} \E\left[  W\left(\frac{h_t^\rv}{\arate t} e^{\beta_p b} \right)\right] + \frac{\arate}{\beta_p} \E\left[  W\left(\frac{h_t^\rv}{\arate t} e^{\beta_p a} \right)\right] \right|\\
    &\leq \left(1-\arate \E\left[r_t^\rv\right]\right)\left|b - a\right|
\end{align*}
showing that $F_t$ is a contraction with constant $\gamma_t = 1-\arate \E\left[r_t^\rv\right] \in (0,1)$.\\
Note that, since for all $\omega$, $h_t^\omega$ converges to a positive constant, and $\E\left[h_t^\rv\right]$ converges to $k$, then for $t$ large enough, $\E\left[r_t^\rv\right] \geq \frac{k}{2 \arate}e^{\beta_p M} \cdot \frac{1}{t}$. It follows that $\sum_{t=1}^\infty \E\left[r_t^\rv\right] = \infty$. At the same time, since $\log(1-z)\leq -z$ for all $z \in [0,1]$ we have for any $n \in \N$:
\begin{align*}
    \ln \prod_{t=1}^\infty \gamma_t &\leq \ln \prod_{t=1}^n \gamma_t\\
    &= \sum_{t=1}^n \ln \left(1-\arate \E\left[r_t^\rv\right]\right)\\
    &\leq -\arate \sum_{t=1}^n \E\left[r_t^\rv\right]
\end{align*}
This is true for any $n \in \N$. By taking the limit as $n$ tends to $\infty$, we obtain $\ln \prod_{t=1}^\infty \gamma_t = -\infty$, and thus $\prod_{t=1}^\infty \gamma_t = 0$.\\
For any $t$, two cases arise:
\begin{itemize}
    \item if $\epsilon_t \leq \epsilon_t^*$, then $\epsilon_{t+1} \geq \epsilon_{t}$ and:
    \begin{align*}
        |\epsilon_{t+1}| &\leq |F_t(\epsilon_t) -  F_t(\epsilon_t^*)| + |\epsilon_t^*|\\
        &\leq |\epsilon_t -  \epsilon_t^*| + |\epsilon_t^*|\\
        &\leq |\epsilon_t| + 2|\epsilon_t^*|
    \end{align*}
    \item if $\epsilon_t > \epsilon_t^*$, then $\epsilon_{t+1} < \epsilon_{t}$ and:
    \begin{align*}
        \epsilon_{t+1} &\geq \epsilon_t+\frac{k}{\beta_p} \ln \left(\frac{t+1}{t}\right)\\
        &> \epsilon_t^*+\frac{k}{\beta_p} \ln \left(\frac{t+1}{t}\right)
    \end{align*}
\end{itemize}
Let $t^*$ be the smallest integer greater than $N$ and satisfying $\left|\frac{k}{\beta_p} \ln \left(\frac{t^*+1}{t^*}\right)\right| \leq \delta$. Then for any $t \geq t^*$ the sequence $(\epsilon_t)$ decreases by at most $\delta$ at each step and increases if $\epsilon_t \leq \epsilon^*_t$. Thus we obtain $\epsilon_t \geq \epsilon_t^* - \delta \geq -2\delta$ for all $t \geq t^*$.\\
Let $I_\delta = [-2\delta, 4\delta]$. Suppose that $\epsilon_t \in I_\delta$ for $t \geq t^*$. If $\epsilon_t \geq \epsilon_t^*$ then $-2\delta \leq \epsilon_{t+1} \leq \epsilon_t$ and $\epsilon_{t+1} \in I_\delta$. Else if $\epsilon_t \leq \epsilon_t^*$ then $|\epsilon_t|\leq 2\delta$ so $\epsilon_t \leq \epsilon_{t+1} \leq |\epsilon_{t+1}| \leq |\epsilon_t| + 2 |\epsilon_t^*| \leq 4\delta$ and $\epsilon_{t+1} \in I_\delta$. Hence we have that $\epsilon_t \in I$ implies $\epsilon_{t+1} \in I_\delta$.\\
Suppose by contradiction that $\epsilon_t \geq 4 \delta$ for all $t \geq t^*$. Then:
\begin{align*}
    \epsilon_{t+1} &\leq \epsilon_t^* + \gamma_t (\epsilon_t - \epsilon_t^*)\\
    &\leq \delta + \gamma_t (\epsilon_t - \delta)\\
    &\leq \delta + \gamma_t \gamma_{t-1} (\epsilon_{t-1} - \delta)\\
    &\hdots\\
    &\leq \delta + (\epsilon_{t^*} - \delta) \prod_{j=t^*}^t \gamma_j
\end{align*}

Thus $\lim_{t \to \infty} \epsilon_t \leq \delta$ which contradicts $\epsilon_t \geq 4 \delta$ for all $t$. Therefore there must exist a $N' \in \N$ such that $\epsilon_t \in I_\delta$ for all $t \geq N'$.\\
In conclusion, for any $\delta > 0$, there exists $N' \in \N$ such that $\epsilon_t \in I_\delta$ for all $t \geq N'$. Hence we have that $\lim_{t \to \infty} \epsilon_t = 0$.\\

Let's now address the convergence rate of $\epsilon_t$. As seen previously, the Taylor series of $\ln$ and $W$ give:
\begin{align*}
\epsilon_{t+1} &= \epsilon_t + \frac{k}{\beta_p}\frac{1}{t} - \frac{ \E\left[h_t^\rv\right]}{\beta_p} \frac{e^{\beta_p \epsilon}}{t} + O\left(\frac{1}{t^2}\right)
\end{align*}
From our hypothesis of induction, $\epsilon_t^{-i} = O\left(\frac{\ln^{k-1}(t)}{t}\right)$. Thus we obtain $\E\left[h_t^\rv\right] = k + O\left(\frac{\ln^{k-1}(t)}{t}\right)$ and, combined with $e^{\beta_p \epsilon_t} = 1 + \beta_p \epsilon_t + \beta_p \epsilon_t \eta_t$ where $\eta_t = o(1)$, we obtain:
\begin{align*}
\epsilon_{t+1} &= \epsilon_t + \frac{k}{\beta_p}\frac{1}{t}\left(1-e^{\beta_p \epsilon}\right) + O\left(\frac{\ln^{k-1}(t)}{t^2}\right)\\
&= \epsilon_t + \frac{k}{t}\left(- \epsilon_t - \epsilon_t \eta_t \right) + O\left(\frac{\ln^{k-1}(t)}{t^2}\right)\\
&= \epsilon_t \left(1 - \frac{k(1+\eta_t)}{t}\right) + O\left(\frac{\ln^{k-1}(t)}{t^2}\right)
\end{align*}
For $t$ large enough, we have $k(1+\eta_t) \in [1,t]$. Indeed, if $k=1$ then $\epsilon_t^{-i} = 0$ so $\E\left[h_t^\omega\right] = k$ and $\eta_t = 0$. Otherwise, $k \geq 2$ and since $\eta_t = o(1)$ there exists $t^*$ such that $\eta_t \in \left[\frac{1}{k}-1,\frac{t}{k}-1\right]$ and $k(1+\eta_t) \in [1,t]$ for all $t \geq t^*$. Thus Lemma \ref{lemma:convergence_rate_noisy_sequence} applies and we obtain $\epsilon_t =  O\left(\frac{\ln^{k}(t)}{t}\right)$.\\

The base case ($k=0$) is trivially satisfied as $z_t^\emptyset = v_t^\emptyset = 0$ for all $t \in \N$. We conclude by induction that $\epsilon_t^S =  O\left(\frac{\ln^{|S|}(t)}{t}\right)$ for all $S \in \mathcal{P}(\mathcal{L})$.
\end{cproof}

\begin{cproof}{Theorem}{\ref{theorem:asymptotic_costs}}
We remark that the value function satisfy the recurrence in \eqref{eq:revenue_recurrence} with $z_t^{S} = V_t(S)$ for all $t \in \N$ and all $S \subseteq \bar{S}$. The first result follows from Lemma \ref{lemma:asymptotic_costs}. We can then use this result to express for any set $S \subseteq \bar{S}$ the marginal value of option $i$ as $t$ tends to $+\infty$:
\begin{align*}
    \Delta_i V^*_t(S) &= \frac{-1}{\beta_p} \left(\ln(\arate t) - 1 + \kappa_i \right) + O\left(\frac{\ln^{|S|}(t)}{t}\right)
\end{align*}
Therefore, for all option $i$ and all type $\omega$ we have $e^{q_i^\omega + \beta_p \Delta_i V^*_t(S)} = \frac{1}{\arate t}e^{q_i^\omega - \kappa_i + 1}+O\left(\frac{\ln^{|S|}(t)}{t^2}\right)$. Since $W(z) = z + O(z^2)$ as $z \to 0$, we obtain $\Gamma_t^\omega(S) = \frac{1}{\arate t}\sum_{i \in S}e^{q_i^\omega - \kappa_i} + O\left(\frac{\ln^{|S|}(t)}{t^2}\right)$. Using the result of Lemma \ref{lemma:inner_problem_solution} we have:
\begin{align*}
    p_{t,i}^{*,\omega} &= \Delta_i V^*_t(S) - \frac{1+\Gamma_t^\omega(S)}{\beta_p}\\
    &= \frac{-1}{\beta_p} \left(\ln(\arate t) + \kappa_i \right) + O\left(\frac{\ln^{|S|}(t)}{t}\right)
\end{align*}
and
\begin{align*}
    \rho_{t,i}^{*,\omega} &= \frac{\Gamma_t^\omega(S)}{1+\Gamma_t^\omega(S)} \frac{e^{q_i^\omega + \beta_p \Delta_i V^*_t(S)}}{\sum_{j \in S} e^{q_j^\omega + \beta_p \Delta_j V^*_t(S)}}\\
    &= e^{q_i^\omega + \beta_p \Delta_i V^*_t(S)} \cdot \left(e^{-1}+O\left(\frac{\ln^{|S|}(t)}{t}\right)\right)\\
    &= \frac{1}{\arate t}e^{q_i^\omega - \kappa_i}+O\left(\frac{\ln^{|S|}(t)}{t^2}\right)
\end{align*}
Since $T(\ear)$ and $\ear$ grow proportionally, we obtain the desired result.
\end{cproof}

\begin{cproof}{Theorem}{\ref{theorem:revenue_utility}}
Let $\Gamma_t^\omega(S) \triangleq W\left(\sum_{i \in S} e^{q_i^\omega + \beta_p \Delta_i V^*_t(S) - 1}\right)$. At any time period $t \geq 1$ of problem \eqref{eq:bellman_equation}, we have for all $\omega \in \Omega$:
\begin{align*}
    \sum_{i \in S} e^{u_{t,i}^{\omega,*}} &= \sum_{i \in S} e^{q_i^\omega + \beta_p p_{t,i}^*}\\
    &\numsym{1}{=} \sum_{i \in S} e^{q_i^\omega + \beta_p \Delta_i V^*_t(S) - 1 - \Gamma_t^\omega(S)}\\
    &= e^{-\Gamma_t^\omega(S)} \sum_{i \in S} e^{q_i^\omega + \beta_p \Delta_i V^*_t(S) - 1}\\
    &\numsym{2}{=} \Gamma_t^\omega(S)
\end{align*}
where (1) follows from Lemma \ref{lemma:inner_problem_solution} and (2) from $W(z) = z\cdot e^{-W(z)}$. Taking the expectation over $\omega$ on both sides we obtain:
$$\sum_{i \in S} \mathbb{E}_\rv\left[e^{u_{t,i}^{\rv,*}}\right] = \mathbb{E}_\rv\left[\Gamma_t^\rv(S)\right]$$
Therefore, since $\arate>0$ and $\sum_{i \in S} \xi_i = \sum_{i \in S'} \xi_i$ for any $S' \in \mathcal{P}(\mathcal{L})$, we have for any $T \geq 1$:
\begin{align*}
\argmax_{S \in \mathcal{P}(\mathcal{L})} V^*_{T(\ear)}(S) = \argmax_{S \in \mathcal{P}(\mathcal{L})} \sum_{i \in S} \xi_i + \arate \sum_{t=1}^{T(\ear)} \mathbb{E}_\rv\left[\Gamma_t^\rv(S)\right] = \argmax_{S \in \mathcal{P}(\mathcal{L})} \ln \sum_{t=1}^{T(\ear)} \sum_{i \in S} \mathbb{E}_\rv\left[e^{u_{t,i}^{\rv,*}}\right]
\end{align*}
\end{cproof}

\subsection{Proofs of Section \ref{section:bounds}}

\begin{cproof}{Proposition}{\ref{prop:upper_bound}}
Let the hypothesis of induction be that, for a given $t$, $V^*_t(S) \leq V^U_{t}(S)$ for all $S \subseteq \bar{S}$. Then for any $S \subseteq \bar{S}$:
\begin{align*}
    V^*_{t+1}(S) &= \displaystyle \E_\rv\left[ \max_{p \in \R^{|S|}} (1-\arate + \arate f_0(p)) V_t^{*}(S) + \arate \sum_{i \in S} f_i^\rv(p, S) \left(p_i + V_t^{*}(S\backslash\{i\})\right)\right]\\
    &\numsym{1}{\leq} \displaystyle \E_\rv\left[ \max_{p \in \R^{|S|}} (1-\arate + \arate f_0(p)) V^U_t(S) + \arate \sum_{i \in S} f_i^\rv(p, S) \left(p_i + V^U_t(S\backslash\{i\})\right)\right]\\
    &\numsym{2}{=} V^U_t(S) + \arate \cdot \E_\rv\left[ \max_{p \in \R^{|S|}} \sum_{i \in S} f_i^\rv(p, S) \left(p_i + V^U_t(S\backslash\{i\}) - V^U_t(S)\right)\right]\\
    &\numsym{3}{=} \sum_{i \in S}r_{t,i} + \arate \cdot \E_\rv\left[ \max_{p \in \R^{|S|}} \sum_{i \in S} f_i^\rv(p, S) \left(p_i - r_{t,i}\right)\right]\\
    &\numsym{4}{\leq} \sum_{i \in S}r_{t,i} + \arate \cdot \E_\rv\left[ \max_{p \in \R^{|S|}} \sum_{i \in S} f_i^\rv(p, \{i\}) \left(p_i - r_{t,i}\right)\right]\\
    &\numsym{5}{=} \sum_{i \in S}r_{t,i} + \arate \cdot \E_\rv\left[ \sum_{i \in S} \max_{p \in \R} f_i^\rv(p, \{i\}) \left(p - r_{t,i}\right)\right]\\
    &\numsym{6}{=} \sum_{i \in S} r_{t+1,i}\\
    &= V^U_{t+1}(S)
\end{align*}
where (1) uses the hypothesis of induction, (2) substitutes $f_0(p)$ with $1-\sum_{i \in S} f_i(p)$, (3) uses the definition of $V^U_t$, (4) follows from the decreasing monotonicity of $f_i^\omega(p,S)$ with respect to $S$ for any price vector $p$, and that any optimal $p$ satisfies $p_i \geq r_{t,i}$ (Lemma $\ref{lemma:inner_problem_solution}$), (5) results from the separability of the maximization problem, and (6) uses the linearity of the expectation and the definition of $r_{t+1,i}$.\\
The base case is immediate as $V^U_0(S) = \sum_{i \in S} \xi_i = V^*_0(S)$, and we conclude by induction that $V^U_{t}(S) \geq V^*_{t}(S)$ for all $t$ and for all $S \subseteq \bar{S}$.\\

We know from Theorem $\ref{theorem:asymptotic_costs}$ applied to any single option set $\{i\}$ that $\lim_{t \to \infty} r_{t,i} - v_t^{\{i\}} = 0$. Note that $V^U_t = \sum_{i \in S} r_{t,i}$ and $v^S_t = \sum_{i \in S} v_t^{\{i\}}$. Therefore $\lim_{t \to \infty} V^U_{t}(S) - v^S_t = 0$. Since Theorem $\ref{theorem:asymptotic_costs}$ also states that $\lim_{t \to \infty} V^*_{t}(S) - v^S_t = 0$ we obtain $\lim_{t \to \infty} V^*_{t}(S) -  V^U_{t}(S) = 0$.
\end{cproof}

\begin{cproof}{Proposition}{\ref{prop:prices_upper_bound}}
Let $g_t(p) = f_i^\omega(p,\{i\})(p-r_{t,i})$ for $\omega \in \Omega$ and $i \in \bar{S}$ fixed. We first note that $\max_{p \in \R} g_t(p) \geq 0$ as choosing $p = r_{t,i}$ gives $g_t(p) = 0$. It follows from \eqref{eq:u_recursion} that $r_{t,i}$ is non-decreasing in $t$.\\
From $\textit{p.3}$ of definition \ref{def:choice_funciton}, $\lim_{p \to \infty} g_t(p) = 0$ so $g_t$ admits a maximum. Hence the maximum of $g_t$ is achieved for $p^*$ satisfying $g_t'(p^*) = 0$:
$$\frac{\partial}{\partial p}f_i^\omega(p^*,\{i\})(p^*-r_{t,i}) + f_i^\omega(p^*,\{i\}) = 0$$ or equivalently:
$$p^* + \frac{f_i^\omega(p^*,\{i\})}{\frac{\partial}{\partial p}f_i^\omega(p^*,\{i\})} = r_{t,i}$$
Since the partial derivative of the left-hand side with respect to $p^*$ is non-negative ($\textit{p.5}$ of definition \ref{def:choice_funciton}) we obtain $\frac{\partial r_{t,i}}{\partial p^*} \geq 0$. Thus the fact that $r_{t,i}$ is non-decreasing in $t$ implies that $\tau^U_{t,i}$ is also non-decreasing in $t$.
\end{cproof}

\begin{cproof}{Proposition}{\ref{prop:lower_bound}}
Let the hypothesis of induction be that, for a given $t$, $V^L_{t}(S) \leq V^*_t(S)$ for all $S \subseteq \bar{S}$. From Lemma $\ref{lemma:inner_problem_solution}$ we know that $\hat{p}_{t,i}^\omega \geq l_{t,i}$. Then for any $S \subseteq \bar{S}$:
\begin{align*}
    V^*_{t+1}(S) &= \displaystyle \E_\rv\left[ \max_{p \in \R^{|S|}} (1-\arate + \arate f_0(p)) V_t^{*}(S) + \arate \sum_{i \in S} f_i^\rv(p, S) \left(p_i + V_t^{*}(S\backslash\{i\})\right)\right]\\
    &\numsym{1}{\geq} \displaystyle \E_\rv\left[ \max_{p \in \R^{|S|}} (1-\arate + \arate f_0(p)) V^L_t(S) + \arate \sum_{i \in S} f_i^\rv(p, S) \left(p_i + V^L_t(S\backslash\{i\})\right)\right]\\
    &\numsym{2}{=} V^L_t(S) + \arate \cdot \E_\rv\left[ \max_{p \in \R^{|S|}} \sum_{i \in S} f_i^\rv(p, S) \left(p_i + V^L_t(S\backslash\{i\}) - V^L_t(S)\right)\right]\\
    &\numsym{3}{=} \sum_{i \in S}l_{t,i} + \arate \cdot \E_\rv\left[ \max_{p \in \R^{|S|}} \sum_{i \in S} f_i^\rv(p, S) \left(p_i - l_{t,i}\right)\right]\\
    &\numsym{4}{\geq} \sum_{i \in S}l_{t,i} + \arate \cdot \E_\rv\left[ \sum_{i \in S} f_i^\rv(\hat{p}_t^\rv, S) \left(\hat{p}_{t,i}^\rv - l_{t,i}\right)\right]\\
    &\numsym{5}{\geq} \sum_{i \in S}l_{t,i} + \arate \cdot \E_\rv\left[ \sum_{i \in S} f_i^\rv(\hat{p}_t^\rv, \bar{S}) \left(\hat{p}_{t,i}^\rv - l_{t,i}\right)\right]\\
    &\numsym{6}{=} \sum_{i \in S} l_{t+1,i}\\
    &= V^L_{t+1}(S)
\end{align*}
where (1) uses the hypothesis of induction, (2) substitutes $f_0(p)$ with $1-\sum_{i \in S} f_i(p)$, (3) uses the definition of $V^L_t$, (4) results from the suboptimality of $\hat{p}_t^\omega$, (5) follows from $f_i^\omega(p, S) \geq f_i^\omega(p,\bar{S})$ for any price vector $p$ and that $\hat{p}_{t,i}^\omega \geq l_{t,i}$, and (6) uses the linearity of the expectation and the definition of $l_{t+1,i}$.\\
The base case is immediate as $V^L_0(S) = \sum_{i \in S} \xi_i = V^*_0(S)$, and we conclude by induction that $V^L_{t}(S) \leq V^*_{t}(S)$ for all $t$ and for all $S \subseteq \bar{S}$.
\end{cproof}

\begin{lemma}\label{lemma:convexity_tilde_rho}
For any vector $p \in \R^N$, the function $\tilde{\rho}_i^\omega(p, \cdot)$ is convex and non-increasing, where $\tilde{\rho}_i^\omega$ is defined in \eqref{eq:tilde_rho_def}.
\end{lemma}
\begin{cproof}{Lemma}{\ref{lemma:convexity_tilde_rho}}
Let's fix $p \in \R^{|S|}$. Let:
\begin{align*}
    &f:\;\begin{array}{c @{} c @{} l} 
        \R^+ &\to&\; [0,1]\\
        x &\mapsto&\; 1/(1+x)
    \end{array}
    &g:\;\begin{array}{c @{} c @{} l} 
        [0,1]^N &\to&\; \R^+\\
        a &\mapsto&\; 1 + e^{q_i^\omega + \beta_p p_i} + \sum_{j \in S \backslash\{i\}} a_j e^{q_j^\omega + \beta_p p_j}
    \end{array}
\end{align*}
We first note that is $f$ convex and decreasing on its domain of definition since for any $x \geq 0$, $f'(x) = -1/(1+x)^2 < 0$ and $f''(x) = 2/(1+x)^3 > 0$. As $g$ is linear and non-decreasing in $a$, we obtain that $f \circ g (a)$ is convex and non-increasing. It follows that $\tilde{\rho}_i^\omega(p, \cdot) = e^{q_i^\omega + \beta_p p_i} \cdot f \circ g(\cdot)$ is convex and non-increasing.
\end{cproof}

\begin{lemma}\label{lemma:telescopic_inequality}
Let $\tau$ be a homogeneous and non-increasing price trajectory and $A^1, A^2 \in \R^{T \times N}$ be such that $A^1 \leq A^2$ and $A^1_{T} = A^2_{T}$. Then:
$$\sum_{i} A^1_{0,i} \xi_i + \sum_{i} \sum_{t=1}^T \left( A^1_{t,i} - A^1_{t-1,i}\right) \tau_{t,i} \geq \sum_{i} A^2_{0,i} \xi_i + \sum_{i} \sum_{t=1}^T \left( A^2_{t,i} - A^2_{t-1,i}\right) \tau_{t,i}$$
\end{lemma}

\begin{cproof}{Lemma}{\ref{lemma:telescopic_inequality}}
Since $A^1 \leq A^2$ and that $\tau_{T,i} \geq \tau_{T-1,i} \geq \hdots \geq \tau_{1,i} \geq \xi_i$, we have $\left(A^2_{t,i}-A^1_{t,i}\right)\left(\tau_{t,i}-\tau_{t+1,i}\right) \leq 0$ for all $(i,t)$. It follows that for any $i$:
\begin{align*}
    &A^1_{0,i} \xi_i + \sum_{t=1}^T \left( A^1_{t,i} - A^1_{t-1,i}\right) \tau_{t,i}\\
    &\geq \left[A^1_{0,i} + (A^2_{0,i}-A^1_{0,i})\right] \xi_i + \sum_{t=1}^{T-1} \left[- (A^2_{t-1,i}-A^1_{t-1,i}) + A^1_{t,i} - A^1_{t-1,i} +  (A^2_{t,i}-A^1_{t,i})\right] \tau_{t,i}\\
    &\quad \quad \quad + \left[ - (A^2_{T-1,i}-A^1_{T-1,i}) +  A^1_{T,i} - A^1_{T-1,i} \right] \tau_{T,i}\\
    &= A^2_{0,i} \xi_i + \sum_{t=1}^T \left( A^2_{t,i} - A^2_{t-1,i}\right) \tau_{t,i}
\end{align*}
where the inequality holds since we added $\sum_{t=1}^{T-1} \left(A^2_{t,i}-A^1_{t,i}\right)\left(\tau_{t,i}-\tau_{t+1,i}\right) \leq 0$, and the last equality follows from the cancellation of terms and that $A^1_{T,i} = A^2_{T,i}$ for all $i$. Summing over $i$, we obtain the desired result.
\end{cproof}

\begin{cproof}{Theorem}{\ref{theorem:DFA}}
Let's first show that $A^* \leq \hat{A}$. Assume that for a given $t$, $A^*_{t,i} \leq \hat{A}_{t,i}$ for all $i \in \bar{S}$. We then have for all $i \in \bar{S}$:
\begin{align*}
    A^*_{t-1,i} &= A^*_{t,i} \cdot \left(1 - \arate \E_{\rv,A}\left[\tilde{f}_i^\rv\left(\tau_t, A_{t}\right)\right]\right)\\
    &\numsym{1}{\leq} A^*_{t,i} \cdot \left(1 - \arate \E_\rv \left[\tilde{f}_i^\rv\left(\tau_t,\E_A\left[ A_{t}\right]\right)\right]\right)\\
    &= A^*_{t,i} \cdot \left(1 - \arate \E_\rv \left[\tilde{f}_i^\rv\left(\tau_t,A_t^*\right)\right]\right)\\
    &\numsym{2}{\leq} \hat{A}_{t,i} \cdot \left(1 - \arate \E_\rv \left[\tilde{f}_i^\rv\left(\tau_t,\hat{A}_t\right)\right]\right)\\
    &= \hat{A}_{t-1,i}
\end{align*}
where the first equality comes from \eqref{eq:availability_recursion}.\\
Here, (1) uses Jensen's inequality and the convexity of $\tilde{f}_i^\omega$ with respect to $a$, and (2) results from the fact that $\tilde{f}_i^\omega$ is non-increasing with respect to $A$ and the hypothesis of induction. Since $A^*_T = \hat{A}_T = \mathds{1}_{\R^N}$, the induction is initialized and we conclude by induction that $A^* \leq \hat{A}$.\\
We remark that \eqref{eq:availability_recursion} can be reformulated:
\begin{equation} \label{eq:availability_telescoping}
     A^*_{t,i}\cdot \arate \E_{\rv,A}\left[\tilde{f}_i^\rv(\tau_t, A_t)\right] = A^*_{t,i} - A^*_{t-1,i}
\end{equation}
Similarly, we obtain from the recursive definition of $\hat{A}$:
\begin{equation} \label{eq:estimate_availability_telescoping}
     \hat{A}_{t,i}\cdot \arate \E_\rv\left[\tilde{f}_i^\rv(\tau_t, \hat{A}_t)\right] = \hat{A}_{t,i} - \hat{A}_{t-1,i}
\end{equation}
Then:
\begin{align*}
    V_T^*(\bar{S}) &\geq V_T(\bar{S}, \tau)\\
    &= \sum_{i \in \bar{S}} A^*_{0,i} \xi_i + \sum_{i \in \bar{S}} \sum_{t=1}^T A^*_{t,i} \cdot \arate \E_{\rv,A}\left[\tilde{f}_i^\rv(\tau_t, A_t)\right] \tau_{t,i}\\
    &\numsym{1}{=} \sum_{i \in \bar{S}} A^*_{0,i} \xi_i + \sum_{i \in \bar{S}} \sum_{t=1}^T \left( A^*_{t,i} - A^*_{t-1,i}\right) \tau_{t,i}\\
    &\numsym{2}{\geq} \sum_{i \in \bar{S}} \hat{A}_{0,i} \xi_i + \sum_{i \in \bar{S}} \sum_{t=1}^T \left( \hat{A}_{t,i} - \hat{A}_{t-1,i}\right) \tau_{t,i}\\
    &\numsym{3}{=} \sum_{i \in \bar{S}} \hat{A}_{0,i} \xi_i + \sum_{i \in \bar{S}} \sum_{t=1}^T \hat{A}_{t,i} \cdot \arate \E_\rv \left[\tilde{f}_i^\rv(\tau_t, \hat{A}_t)\right] \tau_{t,i}\\
    &= V_T^\text{DFA}(\bar{S}, \tau)
\end{align*}
where (1) uses equation $\eqref{eq:availability_telescoping}$, (2) follows from Lemma \ref{lemma:telescopic_inequality} combined with the fact that $\hat{A} \geq A^*$, $A^*_{T,i} = 1 = \hat{A}_{T,i}$ for all $i \in \bar{S}$, and that $\tau_{T,i} \geq \tau_{T-1,i} \geq \hdots \geq \tau_{1,i} \geq \xi_i$, and (4) uses equation $\eqref{eq:estimate_availability_telescoping}$.
\end{cproof}

\begin{cproof}{Theorem}{\ref{theorem:performance_bound}}
Let $S \in \mathcal{P}\left(\mathcal{L}\right)$ be of cardinality $N$, $\tau$ be as defined in \eqref{eq:tau^U}, and $\delta_{T,i} \triangleq \arate \E_{\rv}\left[\tilde{\rho}_i^\rv(\tau_{T}, \mathds{1}_{\R^N})\right]$ be the probability that item $i$ is accepted at the first period with horizon $T$. Let $\hat{A}(T)$ and $\hat{A}(T-1)$ be the probabilities of availability associated to $V_T^\text{DFA}(S, \tau)$ and $V_{T-1}^\text{DFA}(S, \tau)$, respectively. Finally, let $\Tilde{A}(T)$ be such that $\hat{A}_{i,t}(T) = \Tilde{A}_{i,t}(T) \left(1-\delta_{T,i}\right)$ for all $i$ and $t \leq T-1$.\\
We first note that $\hat{A}_{t}(T) \leq \hat{A}_{t}(T-1)$ for all $t \leq T-1$. Indeed, $\hat{A}_{T-1}(T) \leq \hat{A}_{T-1}(T-1) = 1$, and we obtain from \eqref{eq:estimated_availability_recursion} and the fact that $\tilde{\rho}_i^\omega$ is non-increasing in $a$ (Lemma \ref{lemma:convexity_tilde_rho}) that $\hat{A}_{t}(T) \leq \hat{A}_{t}(T-1)$ for a given $t$ implies $\hat{A}_{t-1}(T) \leq \hat{A}_{t-1}(T-1)$. We conclude by induction that $\hat{A}_{t}(T) \leq \hat{A}_{t}(T-1)$ for all $t \leq T-1$.\\
Substituting $\Tilde{A}(T)$ in \eqref{eq:estimate_availability_telescoping} we obtain:
$$\Tilde{A}_{t,i}(T) \cdot \arate\E_\rv\left[\tilde{\rho}_i^\rv(\tau_t, \hat{A}_t(T))\right] = \Tilde{A}_{t,i}(T) - \Tilde{A}_{t-1,i}(T)$$
With an induction argument similar to the one above and the fact that $\hat{A}_{t}(T) \leq \hat{A}_{t}(T-1)$ we obtain $\tilde{A}_{t}(T) \leq \hat{A}_{t}(T-1)$ for all $t \leq T-1$. In addition $\Tilde{A}_{T-1,i}(T) = 1 = \hat{A}_{T-1,i}(T-1)$ for all $i \in S$. As $\tau$ is homogeneous and non-increasing (Proposition \ref{prop:prices_upper_bound}), Lemma \ref{lemma:telescopic_inequality} gives for all $i\in S$:
$$\tilde{A}_{0,i}(T) \xi_i + \sum_{t=1}^{T-1} \tilde{A}_{t,i}(T) \cdot \arate \E_\rv \left[\tilde{\rho}_i^\rv(\tau_t, \hat{A}_t(T))\right] \tau_{t,i} \geq \hat{A}_{0,i}(T-1) \xi_i + \sum_{t=1}^{T-1} \hat{A}_{t,i}(T-1) \cdot \arate \E_\rv \left[\tilde{\rho}_i^\rv(\tau_t, \hat{A}_t(T-1))\right] \tau_{t,i}$$
We note from its definition in \eqref{eq:def_DFA} that $V_{T}^\text{DFA}(S, \tau)$ is separable in $i$. Let $V_{T,i}^\text{DFA}(S, \tau)$ be the contribution of option $i$ to $V_{T}^\text{DFA}(S, \tau)$. We can lower bound $V_{T,i}^\text{DFA}(S, \tau)$ as:
\begin{align*}
    V_{T,i}^\text{DFA}(S, \tau)
    &= \hat{A}_{0,i}(T) \xi_i + \sum_{t=1}^T \hat{A}_{t,i}(T) \cdot \arate \E_{\rv}\left[ \tilde{\rho}_i^\rv(\tau_{t}, \hat{A}_t(T))\right] \tau_{t,i}\\
    &= \hat{A}_{0,i}(T) \xi_i + \sum_{t=1}^{T-1} \hat{A}_{t,i}(T) \cdot \arate \E_{\rv}\left[ \tilde{\rho}_i^\rv(\tau_{t}, \hat{A}_t(T))\right] \tau_{t,i} + \delta_{T,i} \cdot \tau_{T,i}\\
    &= \left(1-\delta_{T,i}\right)\left(\Tilde{A}_{0,i}(T) \xi_i + \sum_{t=1}^{T-1} \tilde{A}_{t,i}(T) \cdot \arate \E_{\rv}\left[ \tilde{\rho}_i^\rv(\tau_{t}, \hat{A}_t(T))\right] \tau_{t,i}\right) + \delta_{T,i} \cdot \tau_{T,i}\\
    &\geq \left(1-\delta_{T,i}\right)\left(\hat{A}_{0,i}(T-1) \xi_i + \sum_{t=1}^{T-1} \hat{A}_{t,i}(T-1) \cdot \arate \E_\rv \left[\tilde{\rho}_i^\rv(\tau_t, \hat{A}_t(T-1))\right] \tau_{t,i}\right) + \delta_{T,i} \cdot \tau_{T,i}\\
    &= \left(1-\delta_{T,i}\right)V_{T-1,i}^\text{DFA}(S, \tau) + \delta_{T,i} \cdot \tau_{T,i}
\end{align*}
Hence $V_{T,i}^\text{DFA}(S, \tau) \geq y_{T,i}$ where $(y_{t,i})_{t \in \N}$ is the sequence satisfying:
\begin{align*}
    y_{t,i} &= \left(1-\delta_{t,i}\right) \cdot y_{t-1,i} + \delta_{t,i} \cdot \tau_{t,i}, \quad \forall t \in \N^*\\
    y_{0,i} &= \xi_i
\end{align*}
and we have  $V_{T}^\text{DFA}(S, \tau) = \sum_{i \in S}V_{T,i}^\text{DFA}(S, \tau) \geq Y_T \triangleq \sum_{i \in S} y_{T,i}$.
Since the price trajectory $(\tau_{t,i})_{t \in \N^*}$ is obtained considering a option $i$ only, Theorem \ref{theorem:asymptotic_costs} gives $\tau_{t,i} = \frac{-1}{\beta_p}\left(\ln(\arate t) + \kappa_i\right) + O\left(\frac{\ln^N(t)}{t}\right)$ so we obtain $e^{q_i^\omega + \beta_p \tau_{t,i}} = \frac{1}{\arate t} e^{q_i^\omega - \kappa_i} + O\left(\frac{\ln^N(t)}{t^2}\right)$. Substituting the price in \eqref{eq:rho_definition} with $\tau_{t,i}$ we obtain:
\begin{align*}
\delta_{t,i} \triangleq \arate \E_\omega\left[\tilde{\rho}_i^\omega(\tau_{t},\mathds{1}_{\R^N})\right] &= \frac{\E_\omega\left[e^{q_i^\omega - \kappa_i}\right]}{t} + O\left(\frac{\ln^N(t)}{t^2}\right)\\
&= \frac{1}{t} + O\left(\frac{\ln^N(t)}{t^2}\right)
\end{align*}
and similarly:
$$\arate \E_{\rv}\left[ \tilde{\rho}_i^\rv(\tau_{t}, 0_{\R^N})\right] = \frac{1}{t} +O\left(\frac{\ln^N(t)}{t^2}\right)$$
In addition, we have for all $\omega \in \Omega$:
\begin{align*}
    \tilde{\rho}_i^\omega(\tau_{t},0_{\R^N}) -\tilde{\rho}_i^\omega(\tau_{t},\mathds{1}_{\R^N})
    &= \frac{e^{q_i^\omega + \beta_p \tau_{t,i}}}{1+e^{q_i^\omega + \beta_p \tau_{t,i}}} - \frac{e^{q_i^\omega + \beta_p \tau_{t,i}}}{1+\sum_{j \in S} e^{q_j^\omega + \beta_p \tau_{t,j}}}\\
    &= \frac{e^{q_i^\omega + \beta_p \tau_{t,i}}\left(\sum_{j \in S\backslash\{i\}} e^{q_j^\omega + \beta_p \tau_{t,j}}\right)}{\left(1+e^{q_i^\omega + \beta_p \tau_{t,i}}\right)\left(1+\sum_{j \in S} e^{q_j^\omega + \beta_p \tau_{t,j}}\right)}\\
    &= O\left(\frac{1}{t^2}\right)
\end{align*}
Thus we have:
\begin{align*}
    &\delta_{t,i} = \frac{1}{t} + \epsilon_t^1\\
    &\arate \E_{\rv}\left[ \tilde{\rho}_i^\rv(\tau_{t}, 0_{\R^N})\right] = \frac{1}{t} + \epsilon_t^0\\
    &\arate \E_{\rv}\left[ \tilde{\rho}_i^\rv(\tau_{t}, 0_{\R^N})\right] - \delta_{t,i} = \epsilon_t^d
\end{align*}
where $\epsilon_t^1 = O\left(\frac{\ln^N(t)}{t^2}\right)$, $\epsilon_t^0 = O\left(\frac{\ln^N(t)}{t^2}\right)$, and $\epsilon_t^d = O\left(\frac{1}{t^2}\right)$.

Let $(x_{t,i})_{t\in \N}$ be the sequence such that $x_{t,i} =  r_{t,i} - y_{t,i}$ for all $t$, where $r_{t,i}$ is the contribution of option $i$ in the upper bound $V^U_T$. Since $\tau_t$ is the optimal solution to problem \eqref{eq:u_recursion}, we obtain:
\begin{align*}
    r_{t,i} &= r_{t-1,i} + \arate \E_{\rv}\left[ \tilde{\rho}_i^\rv(\tau_{t}, \{i\}) \left(\tau_{t,i}-r_{t-1,i}\right)\right]\\
    &=r_{t-1,i} + \arate \E_{\rv}\left[ \tilde{\rho}_i^\rv(\tau_{t}, 0_{\R^N})\right] \left(\tau_{t,i}-r_{t-1,i}\right)
\end{align*}
It follows that:
\begin{align*}
    x_{t,i} &= r_{t-1,i} + \arate \E_{\rv}\left[ \tilde{\rho}_i^\rv(\tau_{t}, 0_{\R^N})\right] \left(\tau_{t,i}-r_{t-1,i}\right) - \left(\left(1-\delta_{t,i}\right) \cdot y_{t-1,i} + \delta_{t,i} \cdot \tau_{t,i}\right)\\
    &= r_{t-1,i} - y_{t-1,i} - \arate \E_{\rv}\left[ \tilde{\rho}_i^\rv(\tau_{t}, 0_{\R^N})\right] r_{t-1,i} + \delta_{t,i} y_{t-1,i} + \tau_{t,i}\left(\arate \E_{\rv}\left[ \tilde{\rho}_i^\rv(\tau_{t}, 0_{\R^N})\right] - \delta_{t,i}\right)\\
    &= x_{t-1,i}\left(1-\frac{1}{t}\right) - \epsilon_t^0 r_{t-1,i} + \epsilon_t^1 y_{t-1,i} + \epsilon_t^d \tau_{t,i}\\
    &= x_{t-1,i}\left(1-\frac{1}{t}\right) - \left(\epsilon_t^1 + \epsilon_t^d\right) r_{t-1,i} + \epsilon_t^1 y_{t-1,i} + \epsilon_t^d \tau_{t,i}\\
    &= x_{t-1,i}\left(1-\frac{1}{t} - \epsilon_t^1\right) + \epsilon_t^d \left(\tau_{t,i} - r_{t-1,i}\right)\\
    &= x_{t-1,i} \left(1-\delta_{t,i}\right) + \epsilon_t^d \left(\tau_{t,i} - r_{t-1,i}\right)
\end{align*}
Since $\delta_{t,i} \in [0,1]$ for all $t \in \N$, we have:
$$x_{T,i} = x_{0,i} \prod_{t=1}^T (1-\delta_{t,i}) + \sum_{t=1}^T \left(\epsilon_t^d \left(\tau_{t,i} - r_{t-1,i}\right) \cdot \prod_{n=t+1}^T (1-\delta_{n,i})\right)$$
Let $C = \sum_{t=1}^\infty |\epsilon_t^1| < \infty$ since $\epsilon_t^1 = O\left(\frac{\ln^N(t)}{t^2}\right)$. We have for any $T \geq t' \geq 0$:
\begin{align*}
    \ln \prod_{t=t'+1}^T (1-\delta_{t,i}) &= \sum_{t=t'+1}^T \ln(1-\delta_{t,i})\\
    &\numsym{1}{\leq} -\sum_{t=t'+1}^T \delta_{t,i}\\
    &= -\sum_{t=t'+1}^T \left(\frac{1}{t} + \epsilon_t^1\right)\\
    &= -\left(\sum_{t=1}^T \frac{1}{t} - \sum_{t=1}^{t'} \frac{1}{t}\right) - \sum_{t=t'+1}^T \epsilon_t^1\\
    &\numsym{2}{\leq} -\left(\ln(T) - \sum_{t=1}^{t'} \frac{1}{t}\right) + C
\end{align*}
where (1) holds since $\ln(1-z) \leq -z$ for all $z \in [0,1]$ and (2) result from $\sum_{n=1}^K \frac{1}{n} \geq \ln(K)$.
Taking the exponential on both sides gives:
\begin{align*}
    \prod_{t=t'+1}^T (1-\delta_{t,i}) &\leq \frac{\exp\left(\sum_{t=1}^{t'} \frac{1}{t} + C\right)}{T} = O\left(\frac{1}{T}\right)
\end{align*}
Hence we obtain $x_{0,i} \cdot \prod_{t=1}^T(1-\delta_{t,i}) = e^C/T = O(1/T)$. From Theorem \ref{theorem:asymptotic_costs} applied to the set containing option $i$ only, we have $\lim_{t \to \infty}\tau_{t,i} - r_{t-1,i} = -1/\beta_p$. Thus, as $\epsilon_t^d = O\left(\frac{1}{t^2}\right)$, there exists $t_0 \in \N^*$ and $M > 0$ such that for all $t \geq t_0$, $\epsilon_t^d \left(\tau_{t,i} - r_{t-1,i}\right) \leq \frac{M}{t^2}$. Hence we have:
\begin{align*}
    \sum_{t=1}^T \left(\epsilon_t^d \left(\tau_{t,i} - r_{t-1,i}\right) \prod_{n=t+1}^T (1-\delta_{n,i})\right) &\leq \underbrace{\sum_{t=1}^{t_0-1} \left(\epsilon_t^d \left(\tau_{t,i} - r_{t-1,i}\right) \prod_{n=t+1}^T (1-\delta_{n,i})\right)}_{(a)} + \underbrace{M \sum_{t=t_0}^{T} \frac{\prod_{n=t+1}^T (1-\delta_{n,i})}{t^2}}_{(b)}
\end{align*}
Since $\prod_{n=t+1}^T (1-\delta_{n,i}) = O(1/T)$, (a) is $O(1/T)$ as a finite sum of $O(1/T)$ elements. For (b), we observe that:
\begin{align*}
    \sum_{t=t_0}^{T} \frac{\prod_{n=t+1}^T (1-\delta_{n,i})}{t^2} &\leq \sum_{t=t_0}^{T} \frac{\exp\left(\sum_{n=1}^{t} \frac{1}{n} + C\right)}{T \cdot t^2}\\
    &\numsym{1}{\leq} \sum_{t=t_0}^{T} \frac{\exp\left(1 + \ln(t) + C\right)}{T \cdot t^2}\\
    &= \frac{e^{1+C}}{T} \sum_{t=t_0}^{T} \frac{1}{t}\\
    &\leq \frac{1+\ln(T)}{T} e^{1+C}
\end{align*}
where (1) follows from $\sum_{n=1}^t \frac{1}{n} \leq \ln(t) + 1$ for all $t \in \N^*$. Thus we have that (b) is $O\left(\frac{\ln(T)}{T}\right)$.\\
Combining the above results we obtain that $x_{T,i} = O\left(\frac{\ln(T)}{T}\right)$. Summing over $i \in S$ yields:
\begin{align*}
    V^U_{T} - Y_{T} = \sum_{i} x_{T,i} = O\left(\frac{\ln(T)}{T}\right)
\end{align*}
From Proposition \ref{prop:upper_bound} and Theorem \ref{theorem:DFA} we have:
$$Y_T \leq V_T^\text{DFA}(S, \tau) \leq V_T^* \leq U_T$$
It follows that:
\begin{align*}
    V_T^\text{DFA}(S) - V_T^*(S) &= O\left(\frac{\ln(T)}{T}\right)\\
    V_T^U(S) - V_T^*(S) &= O\left(\frac{\ln(T)}{T}\right)
\end{align*}
and since Theorem \ref{theorem:asymptotic_costs} gives $V_T^*(S) = \Theta\left(\ln(T)\right)$, we finally obtain:
\begin{align*}
   \frac{V_T^\text{DFA}(S) - V_T^*(S)}{V_T^*(S)} &= O\left(\frac{1}{T}\right)\\
   \frac{V_T^U(S) - V_T^*(S)}{V_T^*(S)} &= O\left(\frac{1}{T}\right)
\end{align*}
which concludes the proof.
\end{cproof}

\end{document}